\documentclass[reqno,pdftex,11pt]{amsart}
\usepackage[top=25truemm,bottom=25truemm,left=25truemm,right=25truemm]{geometry}

\usepackage{amsmath,amsthm,amssymb,amsfonts,amscd,mathtools,color,tikz,tikz-cd}
\usepackage{graphicx} %To submit on arXiv, replace `dvipdfmx' to `pdftex'.
\usetikzlibrary{arrows.meta}
\numberwithin{equation}{section}
\theoremstyle{plain} 
	\newtheorem{thm}{Theorem}[section]
	\newtheorem*{thm*}{Theorem}
	\newtheorem{cor}[thm]{Corollary}
	\newtheorem{lem}[thm]{Lemma}
	
	\newtheorem{prop}[thm]{Proposition}
	
	\newtheorem*{conj*}{Conjecture}
	
\theoremstyle{definition}
	\newtheorem{defn}[thm]{Definition}%[section]

\theoremstyle{remark}
	\newtheorem{rem}[thm]{Remark}
	
	\newtheorem*{pf}{Proof}
\def\EE{{\mathbb E}}
\def\CC{{\mathbb C}}
\def\DD{{\mathbb D}}

\def\HH{{\mathbb H}}

\def\LL{{\mathbb L}}
\def\MM{{\mathbb M}}

\def\PP{{\mathbb P}}

\def\RR{{\mathbb R}}
\def\SS{{\mathbb S}}

\def\XX{{\mathbb X}}
\def\ZZ{{\mathbb Z}}
\def\A{{\mathcal A}}
\def\B{{\mathcal B}}
\def\C{{\mathcal C}}
\def\D{{\mathcal D}}
\def\E{{\mathcal E}}
\def\F{{\mathcal F}}

\def\H{{\mathcal H}}
\def\I{{\mathcal I}}
\def\J{{\mathcal J}}

\def\L{{\mathcal L}}
\def\M{{\mathcal M}}
\def\N{{\mathcal N}}
\def\O{{\mathcal O}}
\def\P{{\mathcal P}}

\def\R{{\mathcal R}}
\def\S{{\mathcal S}}
\def\T{{\mathcal T}}

\def\W{{\mathcal W}}
\def\X{{\mathcal X}}

\def\p{{\partial}}

\def\Aut{{\rm Aut}}

\def\FEC{{\rm FEC}}
\def\Fuk{{\rm Fuk}}
\def\Hom{{\rm Hom}}

\def\Pic{{\rm Pic}}

\def\Stab{{\rm Stab}}
\def\RHom{{\rm \mathbf{R}Hom}}
\def\Cone{{\rm Cone}}
\def\inv{{\rm inv}}

\def\id{{\rm id}}
\begin{document}
\title{Frobenius manifold for the nodal quiver\\ \footnotesize In Memory of Professor Wolfgang Ebeling}
\date{\today}
\author{Atsuki Nakago} 
\address{Department of Mathematics, Graduate School of Science, Osaka University, 
Toyonaka Osaka, 560-0043, Japan}
\email{u400778f@ecs.osaka-u.ac.jp}
\author{Yuuki Shiraishi}
\address{School of Economics and Management, University of Hyogo, 
Nishiku Kobe Hyogo, 651-2197, Japan}
\email{s912y025@guh.u-hyogo.ac.jp}
\author{Atsushi Takahashi}
\address{Department of Mathematics, Graduate School of Science, Osaka University, 
Toyonaka Osaka, 560-0043, Japan}
\email{takahashi@math.sci.osaka-u.ac.jp}

% \maketitle
%%%%%%%%%%%%%%%%%%%%%%%%%%%%%%%%%%%%%%%%%%%%%%%%%%%%%%%%%%%%%%%%
\begin{abstract}
Starting from the Weierstrass elliptic function, we study the associated Frobenius structure, incorporating the perspective of derived categories, particularly that of homological mirror symmetry.
Given a deformation of the Weierstrass elliptic function, we construct a primitive form normalized to be compatible with the period map for integral cycles, and obtain a Frobenius structure 
whose Frobenius potential is defined over the rational numbers.
We also construct a Frobenius structure using elliptic Weyl group invariants (as opposed to Jacobi group invariants), and establish an isomorphism between these two Frobenius structures.
We further examine the relationship between the degree of the Lyashko--Looijenga map modulo the modular group and the number of full exceptional collections up to the braid group action and translations, 
as well as the associated Gamma-integral structure.
\end{abstract}
\maketitle
%%%%%%%%%%%%%%%%%%%%%%%%%%%%%%%%%%%%%%%%%%%%%%%%%%%%%%%%%%%%%%%%
\section{Introduction}

A Frobenius manifold is a complex manifold equipped with a structure of a graded Frobenius algebra on its tangent sheaf satisfying certain flatness conditions. 
By mirror symmetry one expects an isomorphism between a Frobenius manifold arising from the Gromov--Witten theory and another Frobenius manifold arising from deformation theory, 
thereby establishing deep connections between mathematical objects of different origins.
The notion of a Frobenius manifold originates from the flat structure by K.~Saito \cite{Sa1} in his study of primitive forms, introduced as a framework to describe the period integrals and the monodromy associated to 
universal unfoldings of singularities in a unified manner.
Subsequently, this differential-geometric structure was axiomatized by B.~Dubrovin \cite{D} in the context of integrable systems arising from two-dimensional topological field theories.
It is also important that the flat structure first appeared in the study of the invariant ring of a finite reflection group and its orbit space \cite{Sa4, SYS} and 
its regular subspace coincides with the domain for period mappings of (the Gelfand--Leray form associated to) a primitive form on cycles in the Milnor fiber.
This approach has been extended to affine Weyl groups \cite{DZ}, elliptic Weyl groups \cite{Sat1} and Jacobi groups \cite{Be1, Be2}. 
Moreover, for several classes it has been verified that Frobenius structures obtained in this way are isomorphic to those arising from primitive forms (see \cite{SatTa, IST}).
In particular, the period mapping associated with a primitive form plays an important role in proving the isomorphism of Frobenius structures between that arising from deformation theory and that arising from the corresponding Weyl group invariant theory \cite{ShTa}.
In this paper, since we also discuss Weyl group invariants, we have emphasized some of their background. 
However, what is important here is that there are essentially three different ways to construct a Frobenius structure.

In contrast, there are also three types of realizations of derived categories: those based on coherent sheaves over algebraic varieties, those based on Floer homology
of symplectic manifolds, and those based on representations of algebras. 
Several variants of the homological mirror symmetry conjecture proposed by Kontsevich \cite{Kon} claim equivalences among those categories mentioned above
and categorical notions such as full exceptional collections, spherical objects, and auto-equivalence groups have been interpreted as counterparts of monodromy phenomena in the geometric side.
Moreover, at this intersection of mirror symmetry and categorical methods, T.~Bridgeland \cite{Br} introduced the notion of stability conditions for triangulated categories,  
consisting of tuples of a group homomorphism from the Grothendieck group to $\CC$ and an $\RR$-graded family of full additive subcategories
satisfying certain conditions, and showed that the set of stability conditions has a structure of a complex manifold. 
Since the manifold can be regarded as a formulation of (the universal covering of) the complexified K\"ahler moduli space, it is natural to expect that 
if the category is smooth and proper then it is equipped with a Frobenius structure (in fact, more strongly, with a $tt^*$-geometry).
Under this expectation, especially for the cases of the smooth algebraic varieties, 
it is conjectured by Dubrovin~\cite{D3} that the Frobenius structure above is semi-simple if and only if 
the bounded derived category of coherent sheaves has a full exceptional collection,
and that the monodromy data of the semi-simple Frobenius structure are compatible with numerical data of the full exceptional collection.
Dubrovin's conjecture has evolved into Gamma conjecture II (Galkin--Golyshev--Iritani~\cite{GGI}) and
the refined Dubrovin conjecture (Cotti--Dubrovin--Guzzetti~\cite{CDG}).

Thus, while both classical and homological mirror symmetry aim at understanding deformation theory and monodromy, no unified framework directly connecting the two has yet been established. 
In particular, making precise the construction that derives classical mirror symmetry from homological mirror symmetry remains a major open problem in current research on mirror symmetry.
Therefore, the goal here is, setting aside for the moment the direct categorical construction of a Frobenius structure on the space of stability conditions, 
to focus instead on the background and expectations coming from mirror symmetry, and to verify all possible consistencies, 
in a special case where most quantities can be computed by means of elliptic function theory and the elliptic Weyl group
(In the case of affine ADE types, we refer the reader to \cite{IST, STW, Ta1, OST}).

In this paper, starting from the Weierstrass $\wp$-function, we study the associated Frobenius structure, incorporating the perspective of derived categories, in particular that of homological mirror symmetry.
Given a deformation of the Weierstrass $\wp$-function, we construct a primitive form normalized to be compatible with the period map for integral cycles and obtain a Frobenius structure 
whose Frobenius potential is defined over rational numbers.
We also construct a Frobenius structure using elliptic Weyl group invariants (as opposed to Jacobi group invariants), and establish an isomorphism between these two Frobenius structures.
We further examine the relationship between the degree of the Lyashko--Looijenga map modulo the modular group and the number of full exceptional collections up to the braid group action and translations, 
as well as the associated Gamma-integral structure.

The contents of this paper are as follows:
First, in Section~\ref{sec : preliminaries}, we organize Weierstrass $\wp$-function and its basic identities.
We also prepare a minimal set of tools for the nodal quiver and its derived category, which will be used later.
This serves as a preliminary bridge connecting deformation theory, Weyl group invariant theory, and category theory.

In Section~\ref{sec: category}, we describe a triangulated category $\D$, which is one of the main subjects of consideration, 
in three ways and analyze their internal structures based on the categorical equivalence among them.
Specifically, we consider the partially wrapped Fukaya category $\P\W(\SS_1^1,\MM,\ell)$ of the marked bordered torus with two marked points on one boundary $(\SS_1^1,\MM)$ and the line field $\ell$ giving a constant foliation, 
the derived category $\D^{b}{\rm coh}(\A_{\bf E})$ of the Auslander sheaf $\A_{\bf E}$ associated to the nodal cubic curve ${\bf E}$, 
and the derived category $\D^{b}{\rm mod}(\CC Q/I)$ of finite-dimensional representations of the nodal quiver $(Q,I)$.
The equivalence $\D^{b}{\rm coh}(\A_{\bf E})\cong \P\W(\SS_1^1,\MM,\ell)$ is the homological mirror symmetry in this context, 
but what is more important is the triadic relation that includes the equivalence with $\D^{b}{\rm mod}(\CC Q/I)$.
When discussing properties of triangulated categories that do not depend on a particular construction, 
we shall, for simplicity, denote any one of these triangulated categories by $\D$. 
Since concrete descriptions are more accessible, it is often taken to be $\D=\D^{b}{\rm mod}(\CC Q/I)$.
Based on these equivalences, we determine the left action on $\FEC(\D)/\ZZ^{3}$, the set of full exceptional collections up to translations, by the subgroup $B(\D)$ of 
the auto-equivalence group $\Aut(\D)$ generated by twist functors associated to exceptional cycles (including spherical objects).
We show that this action can be identified with the mutation action of braid groups and is transitive.
Furthermore, for the space $\Stab(\D)$ of stability conditions for $\D=\P\W(\SS_1^1,\MM,\ell)$, using the results of Haiden-Katzarkov-Kontsevich \cite{HKK} and Takeda \cite{Tak} we see that
the universal covering of $M=\CC\times\CC^*\times\HH$ is isomorphic to $\Stab(\D)$.
Here $M$ is the parameter space in deformation theory, which will be described later, on which a Frobenius structure is constructed via a primitive form.
In addition, by comparing the Serre dimension on $\Stab(\D)$ with the global dimension function, we show that there is no Serre-invariant stability condition, which will be an important observation for the later consistency with the Lyashko-Looijenga map.

In the first half of Section~\ref{sec:Frob_Str}, we consider the universal unfolding $F(z;{\bf s})$ with ${\bf s}=(s_1,s_2,\tau)\in M$ of the Weierstrass $\wp$-function which gives $(\SS_1^1,\MM)$. 
The critical points of $F(z;{\bf s})$ are the three half-periodic points and do not collide with each other.
Therefore, the three corresponding critical values are distinct on $M$ and serve as global canonical coordinates.
It is shown that, by using the basic identities in Section~\ref{sec : preliminaries}, $2\pi\sqrt{-1}dz$ gives a primitive form with minimal exponent $0$.
The line field $\ell$ coincides with the one obtained from the primitive form.
As a result, we obtain a globally semi-simple Frobenius structure on $M$ of rank $3$ and dimension $1$, 
whose potential and intersection form can be expressed concretely using Eisenstein series.
Furthermore, from classical elliptic function theory, we see that the Lyashko-Looijenga map, a map sending points in $M$ to corresponding tuples of three critical values in
the configuration space ${\rm Conf}(\CC,3)$ of three points in $\CC$, induces the isomorphism ${\rm LL}^{alg}: M/{\rm SL}(2;\ZZ)\longrightarrow {\rm Conf}(\CC,3)$.
It is important that the degree of ${\rm LL}^{alg}$ agrees with the number of orbits ($=1$) of the braid group action on $\FEC(\D)/\ZZ^{3}$.
From the viewpoint of period mappings, if we take a basis of the relative homology group in the domain $M^{reg}$, the complement of the zeros of the discriminant, 
and consider the period integral along the path associated to the critical value, its derivative gives a flat section of the Levi-Civita connection with respect to the intersection form $g$. 
Indeed, we can take the coordinates $(x,\tau)$ corresponding to two particular period integrals and an additional function $\phi$, determined uniquely up to a constant
and naturally coincides with the coordinates of a suitable subspace on which the elliptic Weyl group acts, so that ${d\phi,dx,d\tau}$ make $g$ a constant matrix.

In the second half, we then develop the theory of the elliptic Weyl group invariants associated to $\D=\D^{b}{\rm mod}(\CC Q/I)$.
Namely, we first derive the elliptic root system (of type $A_1^{(1,1)*}$) from $\D$, especially the elliptic Weyl group $W$ as a reflection group on the Grothendieck group $K_0(\D)$.
Since the radical of this root system is of rank $2$, we enlarge the root lattice so that the extension $\widetilde{I}$ of the intersection form $I$ is non-degenerate.
The hyperbolic extension $\widetilde{W}$ of the elliptic Weyl group $W$ is defined as the reflection group with respect to $\widetilde{I}$.   
We shall explain the actions of $\widetilde{W}$ and ${\rm SL}(2;\ZZ)$ on a suitable subspace $\widetilde{\EE}$ of the space of group homomorphisms from
the enlarged root lattice to $\CC$.
This construction matches the deformation theory of the Weierstrass $\wp$-function, namely, in relation to the fundamental group and period mappings.
We introduce the $\widetilde{W}$-invariants $y_1,y_2$ and the $\widetilde{W}$-anti-invariant $\mathcal{J}$ using theta functions and the Weierstrass $\widetilde\wp$ function.
These are generators in the Chevalley-type theorem (Proposition \ref{prop: Chevalley}), and the intersection form, the unit vector field, and the Euler vector field turn out to be 
identical to those from the deformation theory above. Therefore, the two Frobenius structures from different constructions are globally isomorphic by the uniqueness theorem (Proposition~\ref{unique}).
The advantage of the elliptic Weyl group $\widetilde{W}$ rather than the Jacobi group is that the assertion about the $\widetilde{W}$-anti-invariant in the Chevalley-type theorem is correctly stated, 
namely, the discriminant is expressed as the square $\mathcal{J}^2$ of the $\widetilde{W}$-anti-invariant $\mathcal{J}$ and 
is a monic cubic polynomial in the $\widetilde{W}$-invariant $y_1$ of maximal degree, and hence is consistent with the Lyashko-Looijenga map.

In the final Section, Section~\ref{sec:Gamma}, we introduce the analogue of Iritani's gamma integral structure for a "smooth proper non-commutative algebraic variety'' $\D=\D^{b}{\rm coh}(\A_{\bf E})$
and compare this with the monodromy data on the geometric side. More precisely, we define the characteristic class map ${\rm ch}_\Gamma$ from $K_0(\D)$ 
to the Hochschild homology group $H\!H_\bullet(\D)$ as a matrix with respect to specific bases.
By evaluating the asymptotic expansion of the exponential periods along relative cycles at appropriate points and directions, 
and using the fact that their derivatives coincide with the flat solutions of the first structure connection, 
we find that the Stokes matrix is equal to the Euler matrix of the full exceptional collection, and the central connection matrix is given by $(2\pi)^{-1/2}{\rm ch}_\Gamma$.
This embodies the correspondence proposed in the (refined) Dubrovin conjecture and the Gamma Conjecture II 
(the relationship between semi-simplicity of a Frobenius structure and full exceptional collections in a derived category) 
in the case of the non-commutative resolution $\D^{b}{\rm coh}(\A_{\bf E})$ of nodal cubic ${\bf E}$.

\bigskip
\noindent
{\bf Acknowledgements.}
A.T. would like to express his heartfelt gratitude to Professor Wolfgang Ebeling, his long-standing research collaborator and co-author. His dedication greatly shaped the depth of our work. 
It is with great respect and sorrow that we honor his memory through this paper.
A.T. is supported by JSPS KAKENHI Grant Number JP21H04994.
Y.S. is supported by JSPS KAKENHI Grant Number 19K14531 and 23K03111.

%%%%%%%%%%%%%%%%%%%%%%%%%%%%%%%%%%%%%%%%%%%%%%%%%%%%%%%%%%%%%%%%%
\section{Preliminaries}\label{sec : preliminaries}
\subsection{Weierstrass elliptic functions}
Let $\wp(z;\omega_1,\omega_2)$ be the Weierstrass's $\wp$-function 
\[
\wp(z;\omega_1,\omega_2)
:=\frac{1}{z^{2}}
+\!\!\sum_{\substack{(m,n)\in\ZZ^2\\(m,n)\ne(0,0)}}\!
\left(\frac{1}{(z-m\omega_1-n\omega_2)^2}-\frac{1}{(m\omega_1+n\omega_2)^2}\right),
\]
and set
\[
\widetilde\wp(z;\tau):=\frac{1}{(2\pi\sqrt{-1})^2}\wp(z;1,\tau),
\]
which satisfies
\begin{subequations}\label{eq: cubic}
\begin{gather}
\frac{1}{(2\pi\sqrt{-1})^2}\left(\frac{\p\widetilde \wp(z;\tau)}{\p z}\right)^2
=4\widetilde \wp(z;\tau)^3-\frac{1}{12}E_4(\tau)\widetilde \wp(z;\tau)+\frac{1}{216}E_6(\tau),\label{eq: cubic a}\\
\frac{1}{(2\pi\sqrt{-1})^2}\frac{\p^2\widetilde \wp(z;\tau)}{\p z^2}=6\widetilde \wp(z;\tau)^2-\frac{1}{24}E_4(\tau).\label{eq: cubic b}
\end{gather}
\end{subequations}
Here $E_4(\tau)$ and $E_6(\tau)$ are Eisenstein series defined as 
\[
E_4(\tau)=1+240\sum_{n=1}\sigma_3(n){\bf e}[n\tau],\quad E_6(\tau)=1-504\sum_{n=1}\sigma_5(n){\bf e}[n\tau], 
\]
where $\sigma_{2k-1}(n):=\sum_{d|n}d^{2k-1}$ and ${\bf e}[-]:=e^{2\pi\sqrt{-1}(-)}$.

Similarly, let $\zeta(z;\omega_1,\omega_2)$ be the Weierstrass's $\zeta$-function and set
\[
\widetilde\zeta(z;\tau):=\frac{1}{(2\pi\sqrt{-1})^2}\zeta(z;1,\tau),
\]
which satisfies $-\p\widetilde\zeta(z;\tau)/\p z=\widetilde\wp(z;\tau)$ and 
\begin{equation}\label{eq: zeta}
\widetilde\zeta(z+1;\tau)-\widetilde\zeta(z;\tau)=-\frac{1}{12}E_2(\tau),\ 
\widetilde\zeta(z+\tau;\tau)-\widetilde\zeta(z;\tau)=-\frac{1}{12}E_2(\tau)\tau-\frac{1}{2\pi\sqrt{-1}},
\end{equation}
due to Legendre's relation where 
\[
E_2(\tau):=1-24\sum_{n=1}^\infty \sigma_1(n){\bf e}[n\tau].
\]
It is well known that 
\begin{subequations}\label{eq: E d-ring}
{\Small
\begin{gather}
\frac{1}{2\pi\sqrt{-1}}\frac{\p E_2(\tau)}{\p \tau}=\frac{1}{12}\left(E_{2}(\tau)^{2}-E_{4}(\tau)\right),\label{eq: E d-ring a}\\ 
\frac{1}{2\pi\sqrt{-1}}\frac{\p E_{4}(\tau)}{\p \tau}=\frac{1}{3}\left(E_{2}(\tau)E_{4}(\tau)-E_{6}(\tau)\right),\label{eq: E d-ring b}\\
\frac{1}{2\pi\sqrt{-1}}\frac{\p E_{6}(\tau)}{\p \tau}=\frac{1}{2}\left(E_{2}(\tau)E_{6}(\tau)-E_{4}(\tau)\right).\label{eq: E d-ring c}
\end{gather}}
\end{subequations}
\begin{lem}
{\small
\begin{equation}\label{lem: E_2}
\frac{1}{(2\pi\sqrt{-1})^3}\frac{\p^3 E_2(\tau)}{\p \tau^3}=E_2(\tau)\frac{1}{(2\pi\sqrt{-1})^2}\frac{\p^2 E_2(\tau)}{\p \tau^2}-\frac{3}{2}\left(\frac{1}{2\pi\sqrt{-1}}\frac{\p E_2(\tau)}{\p \tau}\right)^2.
\end{equation}}
\end{lem}
\begin{pf}
It follows from \eqref{eq: E d-ring} by direct calculation.
\qed\end{pf}

Note that equations \eqref{eq: cubic a} and \eqref{eq: cubic b} are rewritten in terms of $\wp(z;\tau)$ and $E_2(\tau)$ as 
\begin{subequations}\label{eq: cubic-tilde}
{\Small
\begin{multline}\label{eq: cubic-tilde a}
\frac{1}{(2\pi\sqrt{-1})^2}\left(\frac{\p}{\p z}\left(\widetilde\wp(z;\tau)-\frac{1}{12}E_2(\tau)\right)\right)^2\\
=4\left(\widetilde\wp(z;\tau)-\frac{1}{12}E_2(\tau)\right)^3
+E_2(\tau)\left(\widetilde\wp(z;\tau)-\frac{1}{12}E_2(\tau)\right)^2\\
+\frac{1}{2\pi\sqrt{-1}}\frac{\p E_2(\tau)}{\p\tau} \left(\widetilde\wp(z;\tau)-\frac{1}{12}E_2(\tau)\right)
+\frac{1}{6}\frac{1}{(2\pi\sqrt{-1})^2}\frac{\p^2 E_2(\tau)}{\p \tau^2}
\end{multline}}
{\Small
\begin{multline}\label{eq: cubic-tilde b}
\frac{1}{(2\pi\sqrt{-1})^2}\frac{\p^2}{\p z^2}\left(\widetilde\wp(z;\tau)-\frac{1}{12}E_2(\tau)\right)\\
=6\left(\widetilde\wp(z;\tau)-\frac{1}{12}E_2(\tau)\right)^2
+E_2(\tau)\left(\widetilde\wp(z;\tau)-\frac{1}{12}E_2(\tau)\right)
+\frac{1}{2}\frac{1}{2\pi\sqrt{-1}}\frac{\p E_2(\tau)}{\p\tau},
\end{multline}}
\end{subequations}

\begin{lem}
We have
{\Small 
\begin{multline}\label{eq: key identity-1}
\frac{1}{2\pi\sqrt{-1}}\frac{\p}{\p\tau}\widetilde\zeta(z;\tau)\\
=-\frac{1}{2}\frac{1}{(2\pi\sqrt{-1})^2}\frac{\p}{\p z}\widetilde\wp(z;\tau)
+\frac{1}{12}E_2(\tau)\widetilde\zeta(z;\tau)-\frac{1}{144}E_4(\tau)z
+\left(-\widetilde\zeta(z;\tau)-\frac{1}{12}E_2(\tau)z\right)\widetilde\wp(z;\tau),
\end{multline}}
equivalently, 
{\Small 
\begin{multline}\label{eq: key identity-tilde-1}
\frac{1}{2\pi\sqrt{-1}}\frac{\p}{\p\tau}\left(-\widetilde\zeta(z;\tau)-\frac{1}{12}E_2(\tau)z\right)\\
=-\left(-\widetilde\zeta(z;\tau)-\frac{1}{12}E_2(\tau)z\right)\left(\widetilde\wp(z;\tau)-\frac{1}{12}E_2(\tau)\right)+\frac{1}{2}
\frac{1}{(2\pi\sqrt{-1})^2}\frac{\p}{\p z}\left(\widetilde\wp(z;\tau)-\frac{1}{12}E_2(\tau)\right).
\end{multline}}
\end{lem}
\begin{pf}
By \eqref{eq: zeta}, it is easy to see that the LHS minus RHS is an elliptic function.
We see that at $z=0$ its principal part and its constant term vanish by direct calculation.
\qed\end{pf}
\begin{cor}
We have
{\small
\begin{equation}\label{eq: key identity}
\frac{1}{2\pi\sqrt{-1}}\frac{\p\widetilde \wp(z;\tau)}{\p \tau}
=2\widetilde \wp(z;\tau)^2+\frac{1}{6}E_2(\tau)\widetilde \wp(z;\tau)-\frac{1}{36}E_4(\tau)-\left(-\widetilde\zeta(z;\tau)-\frac{1}{12}E_2(\tau)z\right)\frac{\p \widetilde\wp(z;\tau)}{\p z},
\end{equation}}
equivalently, 
{\Small
\begin{multline}\label{eq: key identity-tilde}
\frac{1}{2\pi\sqrt{-1}}\frac{\p}{\p \tau}\left(\widetilde\wp(z;\tau)-\frac{1}{12}E_2(\tau)\right)\\
=2\left(\widetilde\wp(z;\tau)-\frac{1}{12}E_2(\tau)\right)^2+\frac{1}{2}E_2(\tau)\left(\widetilde\wp(z;\tau)-\frac{1}{12}E_2(\tau)\right)
+\frac{1}{4}\frac{1}{2\pi\sqrt{-1}}\frac{\p E_2(\tau)}{\p\tau}\\
 -\left(-\widetilde\zeta(z;\tau)-\frac{1}{12}E_2(\tau)z\right)\cdot\frac{\p }{\p z}\left(\widetilde\wp(z;\tau)-\frac{1}{12}E_2(\tau)\right).
\end{multline}}
\end{cor}
\begin{pf}
By differentiating both sides of \eqref{eq: key identity-1} and \eqref{eq: key identity-tilde-1} with respect to $z$, we obtain the statements.
\qed\end{pf}

\subsection{Nodal quiver and its representations}
\begin{defn}
Let $(Q, I)$ be the quiver with relations given by the following:
\begin{equation}\label{eq: nodal quiver}
Q:
\begin{tikzcd}
\underset{1}{\bullet} \ar[r, shift left, "{\alpha_1}"] \ar[r, shift
right, "{\beta_1}"'] & \underset{2}{\bullet} \ar[r, shift left,
"{\alpha_2}"] \ar[r, shift right, "{\beta_2}"'] & \underset{3}{\bullet}
\end{tikzcd},\quad
I= \langle\alpha_1\beta_2,\ \beta_1\alpha_2\rangle,
\end{equation}
which we call the {\it nodal quiver}. 
\end{defn}
Let $\CC Q/I$ be the path algebra of the nodal quiver, which is a finite dimensional basic $\CC$-algebra of global dimension $2$. 
It is also {\it gentle} in the sense of {\cite[Section~1.1]{AS}}. Namely, the algebra $\CC Q/I$ satisfies the following conditions:
\begin{enumerate}
\item For each vertex $v$ of $Q$, there are at most two arrows starting at $v$ and at most two arrows ending at $v$.
\item For each arrow $\alpha$ in $Q$, there is at most one arrow $\beta$ such that $\alpha\beta \notin I$, and at most one arrow $\gamma$ such that $\gamma\alpha \notin I$.
\item For each arrow $\alpha$ in $Q$, there is at most one arrow $\beta$ such that $\alpha\beta \in I$, and at most one arrow $\gamma$ such that $\gamma\alpha \in I$.
\item The ideal $I$ is generated by paths of length $2$.
\end{enumerate}

Denote by ${\rm mod}(\CC Q/I)$ the abelian category of finite dimensional $\CC Q/I$-modules and 
by $\D^{b}{\rm mod} (\CC Q/I)$ the bounded derived category of ${\rm mod}(\CC Q/I)$.
The category $\D^{b}{\rm mod} (\CC Q/I)$ has a structure of a $\CC$-linear triangulated category, 
whose translation functor is denoted by $[1]$. For $p\in\ZZ$ the $p$-times composition of $[1]$ 
will be denoted by $[p]$.

To simplify the notation, we write $\D^{b}{\rm mod} (\CC Q/I)$ as $\D$ from now on. 
Let $ \D^{b}{\rm mod}(\CC)$ be the bounded derived category of complexes of finite dimensional $\CC$-modules.
Derived categories of an abelian category have canonical differential graded (dg) enhancements.
Therefore, we have the $\RHom$ functor
\[
\RHom_{\D}(-,-):\D^{op}\times \D\longrightarrow \D^{b}{\rm mod}(\CC),
\]
with the property that for all $X,Y\in\D$ the cohomology group of the complex $\RHom_\D(X,Y)$ (called the $\RHom$ complex) satisfies
$H^p(\RHom_{\D}(X,Y))=\Hom_\D(X,Y[p])$ for all $p\in\ZZ$.
We also have functorial mapping cones for all morphisms in $\D$.
 
\begin{defn}
Let $\T$ be a $\CC$-linear enhanced triangulated category.
Denote by $[1]$ its translation functor and by $\RHom_\T$ its $\RHom$ functor.
\begin{enumerate}
\item
An object $E \in \T$ is called {\it indecomposable} if $E$ is nonzero and $E \cong X \oplus Y$ 
for some $X,Y\in\T$ then $X$ or $Y$ is the zero object.
\item
An object $E \in \T$ is called {\it exceptional} if $\RHom_\T(E,E) =\CC\cdot {\rm id}_E$ in $\D^{b}{\rm mod}(\CC)$.
If an object is exceptional then it is indecomposable.
\item 
A sequence $\E = (E_1, \dots, E_\mu)$ of exceptional objects $E_1, \dots, E_\mu\in\T$ is called 
an {\it exceptional collection} if $\RHom_\T(E_i, E_j) \cong  0$ in $ \D^{b}{\rm mod}(\CC)$ if  $i > j$.
\item 
An exceptional collection $\E$ is called {\it full} if the smallest full triangulated subcategory 
$\langle E_1,\dots, E_\mu\rangle$ of $\T$ containing $E_1,\dots, E_\mu$ is equivalent to $\T$.
\item 
Two full exceptional collections $\E=(E_1, \dots, E_\mu)$, $\E'=(E'_1, \dots, E'_\mu)$ in $\T$ 
are said to be {\it isomorphic} if $E_i \cong E'_i$ for all $i = 1, \dots, \mu$.
\item
A full exceptional collection $\E= (E_1, \dots, E_\mu)$ is called {\it strong} 
if $\RHom_\T(E_i, E_j)$ is isomorphic in $ \D^{b}{\rm mod}(\CC)$ to a complex concentrated in degree zero,
equivalently, $\Hom_\T(E_i,E_j[p])=0$ for all $p\ne 0$.
\end{enumerate}
\end{defn}

Let $S(1), S(2), S(3)$ be the simple $\CC Q/I$-modules given by
\[
S(1):
\begin{tikzcd}
\CC \ar[r, shift left] \ar[r, shift right] & 0 \ar[r, shift left] \ar[r,
shift right] & 0
\end{tikzcd},\quad
S(2):
\begin{tikzcd}
0 \ar[r, shift left] \ar[r, shift right] & \CC \ar[r, shift left] \ar[r,
shift right] & 0
\end{tikzcd},\quad
S(3):
\begin{tikzcd}
0 \ar[r, shift left] \ar[r, shift right] & 0 \ar[r, shift left] \ar[r,
shift right] & \CC
\end{tikzcd},
\]
and $P(1), P(2), P(3)$ the projective $\CC Q/I$-modules defined by
\[
P(3):
\begin{tikzcd}
0 \ar[r, shift left] \ar[r, shift right] & 0 \ar[r, shift left] \ar[r,
shift right] & \CC
\end{tikzcd},
\quad
P(2):
\begin{tikzcd}[ampersand replacement=\&]
0 \ar[r, shift left] \ar[r, shift right,] \& \CC \ar[r, shift left,
"{\begin{psmallmatrix} 1\\0 \end{psmallmatrix}}"] \ar[r, shift right,
"{\begin{psmallmatrix} 0\\1 \end{psmallmatrix}}"'] \& \CC^2
\end{tikzcd},\quad
P(1):
\begin{tikzcd}[ampersand replacement=\&]
\CC \ar[r, shift left, "{\begin{psmallmatrix} 1\\0 \end{psmallmatrix}}"]
\ar[r, shift right, "{\begin{psmallmatrix} 0\\1 \end{psmallmatrix}}"']
\& \CC^2 \ar[r, shift left, "{\begin{psmallmatrix} 1&0\\0&0
\end{psmallmatrix}}"] \ar[r, shift right, "{\begin{psmallmatrix}
0&0\\0&1 \end{psmallmatrix}}"'] \& \CC^2
\end{tikzcd}.
\]
These $\CC Q/I$-modules $S(1)$, $S(2)$, $S(3)$, $P(1)$, $P(2)$, $P(3)$ 
are exceptional in $\D$. 
In $\D$, the sequence $(S(1),S(2),S(3))$ is a full exceptional collection, 
which is not strong, and the sequence $(P(3),P(2),P(1))$ is a full strong exceptional collection.

\begin{defn}
Let $\T$ be a $\CC$-linear triangulated category with the translation functor $[1]$.
\begin{enumerate}
\item
Consider the free abelian group $F$ with generators $\{[X]\,|\,X\in\T\}$ and the subgroup $F_0$
of $F$ generated by $[X]-[Y]+[Z]$ for all exact triangles $X\longrightarrow Y\longrightarrow Z\longrightarrow X[1]$ in $\T$.
The {\it Grothendieck group} $K_0(\T)$ is defined as the quotient group $F/F_0$.
\item
Assume that $\T$ is of {\it finite type}, namely, $\displaystyle\sum_{p\in\ZZ} \dim_\CC\Hom_\T(X,Y[p])<\infty $ for all $X,Y\in\T$.
The {\it Euler form} $\chi_\T:K_0(\T)\times K_0(\T)\longrightarrow \ZZ$ is the $\ZZ$-bilinear form defined as 
\[
\chi_\T([X],[Y]):=\sum_{p\in\ZZ}(-1)^p \dim_\CC \Hom_\T(X,Y[p]),\quad X,Y\in\D,
\]
\item 
An auto-equivalence $\Phi$ of $\T$ is an exact functor $\Phi:\T\longrightarrow\T$ 
which is also an equivalence of $\T$.
Denote by $\Aut(\T)$ the group of auto-equivalences of $\T$.
\item Assume that $\T$ is of {\it finite type}.
A {\it Serre functor} $\S$ is an auto-equivalence of $\T$ with bi-functorial isomorphisms $\Hom_\T(X,Y) \cong\Hom_\T(Y,\S_\T(X))^*$ 
for all $X,Y\in\T$.
\end{enumerate}
\end{defn}
The notion of Serre functors is introduced by Bondal-Kapranov~\cite[Section~3]{BK} and several properties are given there.
If a Serre functor exists, it is unique up to canonical isomorphism which commutes with the above bi-functorial isomorphisms.
Therefore, it is usually referred to as {\it the} Serre functor of $\T$ and is denoted by $\S_\T$.
It is also known that the Serre functor $\S_\T$ commutes with any auto-equivalence of $\T$, namely, $\S_\T\circ \Phi\cong \Phi\circ \S_\T$ for any $\Phi\in \Aut(\T)$. 

The key properties used in this paper are listed below:
\begin{prop}
Let $\T$ be a $\CC$-linear triangulated category of finite type 
with a full exceptional collection $(E_1,\dots,E_\mu)$. 
\begin{enumerate}
\item The Serre functor $\S_\T$ exists. 
\item The subset $\{[E_1],\dots,[E_\mu]\}$ of $K_0(\T)$ forms a $\ZZ$-basis of $K_0(\T)$.
\item Let $[\S_\T]$ be the automorphism of $K_0(\T)$ induced by $\S_\T$.
We have $[\S_\T]=\chi_\T^{-1}\chi_\T^T$.
\end{enumerate}
\end{prop}
\begin{pf}
The first statement follows from \cite[Corollary 3.5]{BK}. 
For each $X\in\T$ we have the following collection of exact triangles
\[
\begin{tikzcd}
0=:X_{\mu+1} \ar[rr] & & X_{\mu}\ar[ld] & \dots   & X_{2}\ar[rr]  & & X_{1}:=X\ar[ld]\\
 & Y_{\mu}\ar[lu,dashed] & & \dots  & & Y_1\ar[lu,dashed] & 
\end{tikzcd}
\]
with $Y_i\in\langle E_i\rangle$. Therefore, 
$K_0(\T)=\displaystyle\bigoplus_{i=1}^\mu K_0(\langle E_i\rangle)=\displaystyle\bigoplus_{i=1}^\mu \ZZ[E_i]$.
By definitions of $\S_\T$ and the Euler form $\chi_\T$, 
we have $\chi_\T([E_i],[E_j])=\chi_\T([E_j],[\S_\T]([E_i]))$.
It follows that $[\S_\T]=\chi_\T^{-1}\chi_\T^T$. 
\qed\end{pf}

In particular, with respect to the $\ZZ$-bases $\{[S(1)], [S(2)], [S(3)]\}$ and $\{[P(3)], [P(2)], [P(1)]\}$ of $K_0(\D)$,
the Euler form $\chi_\D$ is represented by the matrices
\[
\begin{pmatrix}
1 & -2 & 2\\
0 & 1 & -2\\
0 & 0 & 1
\end{pmatrix},\quad
\begin{pmatrix}
1 & 2 & 2\\
0 & 1 & 2\\
0 & 0 & 1
\end{pmatrix},
\]
respectively.

\begin{lem}[{\cite[Lemma~2.3]{Su}}]\label{lem: sung}
The $\CC Q/I$-modules $E_+$ and $E_-$ defined by
\begin{equation}\label{eq: E+-}
E_+=
\begin{tikzcd}
\CC \ar[r, shift left, "1"] \ar[r, shift right, "0"'] & \CC \ar[r, shift
left, "1"] \ar[r, shift right, "0"'] & \CC
\end{tikzcd},\quad
E_-=
\begin{tikzcd}
\CC \ar[r, shift left, "0"] \ar[r, shift right, "1"'] & \CC \ar[r, shift
left, "0"] \ar[r, shift right, "1"'] & \CC
\end{tikzcd}.
\end{equation}
are exceptional, of projective dimension two, and satisfy $\S_{\D}(E_\pm)=E_\mp[2]$.
\qed\end{lem}

%%%%%%%%%%%%%%%%%%%%%%%%%%%%%%%%%%%%%%%%%%%%%%%%%%%%%%%%%%%%%%%%
\section{Marked bordered surface associated to $F$ and Homological Mirror Symmetry}\label{sec: category}
\subsection{Holomorphic function $F$}
Consider the complex manifolds $M=\CC\times \CC^\ast\times\HH$ with the coordinate system $(s_1,s_2,\tau)$ and 
$\CC\times M$ with the coordinate system $(z;s_1,s_2,\tau)$. Sometimes we denote $2\pi\sqrt{-1}\tau$ by $s_3$ and 
$(z;s_1,s_2,\tau)$ by $(z;{\bf s})$. 

Let $\X$ be the quotient of $\CC\times M$ under the group action 
$(z;s_1,s_2,\tau)\mapsto (z+m+n\tau;s_1,s_2,\tau)$, $m,n\in\ZZ$, and 
$p$ the natural projection $\X\longrightarrow M$;
\[
p: \X\longrightarrow M,\quad [(z;s_1,s_2,\tau)]\mapsto (s_1,s_2,\tau).
\]
Let $F: \X\longrightarrow \PP^1$ be a holomorphic map given by
\[
F=F(z;{\bf s})=F(z;s_1,s_2,\tau):=s_2^2\widetilde\wp(z;\tau)+s_1,
\]
and $\X^o$ the space obtained by removing from $\X$ the zero section of $p$, namely,
\[
\X^o:=\X\setminus \{[(0;s_1,s_2,\tau)]\in \X~|~(s_1,s_2,\tau)\in M\},
\]
and we denote the restrictions of $p$ and $F$ to $\X^o$ by the same symbols $p$ and $F$, respectively.
Note that $\X^o=\X\setminus F^{-1}(\infty)$, in other words, $F$ defines a holomorphic function on $\X^o$.

We consider the real blow-up of the point at infinity $\infty$ in $\PP^1=\CC\cup\{\infty\}$, 
which we denote by $\DD_\infty=\CC\cup \{e^{2\pi\sqrt{-1}\theta}\cdot\infty\,|\,0\leq \theta <1\}$.
Define a compact oriented smooth manifold with boundary $\overline{\X}$ as the fiber product of $\X$ and $\DD_\infty$.
Let $\overline{p}:\overline{\X}\longrightarrow M$ be the composition of the natural map $\overline{\X}\longrightarrow \X$ and the map $p:\X\longrightarrow M$.
\[
\begin{tikzcd}
\overline{\X} \arrow[dd, bend right=50, "\overline{p}"'] \arrow[d] \arrow[r, "\overline{F}"] & \DD_\infty \arrow[d]\\
\X\arrow[r, "F"] \arrow[d, "p"']& \PP^1\\
M & 
\end{tikzcd}
\]
Note that for each point ${\bf s}\in M$ the fiber $\overline{p}^{-1}({\bf s})$ is the torus with a single boundary component and 
the restriction of $\overline{F}$ to $\overline{p}^{-1}({\bf s})$
is a double cover of the disk $\DD_\infty$ branched over $3$ points.

\subsection{Homological mirror symmetry}
\begin{defn}
Let $(\SS_1^1, \MM, \ell)$ be a triple of a real surface $\SS_1^1$ of genus $1$ with $1$ boundary component,  
a set $\MM$ of two marked points on the boundary $\p\SS_1^1$ and the line field $\ell\in\Gamma(\SS_1^1,\PP(T\SS_1^1))$ giving a constant foliation on $\SS_1^1$.  
Denote by $\P\W(\SS_1^1, \MM, \ell)$ the partially wrapped Fukaya category associated to the tuple $(\SS_1^1, \MM, \ell)$.
\end{defn}
Note that for each ${\bf s}\in M$ we can identify $\SS_1^1$ and $\MM$ with $\overline{p}^{-1}({\bf s})$ and $\overline{F}^{-1}(-1\cdot\infty)$.

Set
\begin{equation}\label{eq: e}
\widetilde{e}_1(\tau):=\widetilde\wp\left(\frac{1}{2};\tau\right),\quad \widetilde{e}_2(\tau):=\widetilde\wp\left(\frac{1+\tau}{2};\tau\right),\quad \widetilde{e}_3(\tau):=\widetilde\wp\left(\frac{\tau}{2};\tau\right).
\end{equation}
Then we have $\widetilde{e}_2(\sqrt{-1})=0$, $\widetilde{e}_3(\sqrt{-1})\in \RR_{>0}$ and $\widetilde{e}_1(\sqrt{-1})=-\widetilde{e}_3(\sqrt{-1})$.

Choose sufficiently small $\varepsilon, \theta\in\RR_{>0}$ and fix a base point ${\bf s}_0:=(-1, \varepsilon e^{\pi\sqrt{-1}\theta},\sqrt{-1})\in M$.
Consider an ordered system $({\bf p}_3, {\bf p}_2, {\bf p}_1)$ of paths ${\bf p}_i: [0,1]\longrightarrow \DD_\infty$ with ${\bf p}_i(t)=\varepsilon^2 e^{2\pi\sqrt{-1}\theta}\widetilde{e}_i(\sqrt{-1})-1+\frac{t}{1-t}$ and set $L_i=\overline{F}_{{\bf s}_0}^{-1}({\bf p}_i)\subset \SS_1^1$. Note that $L_i$ defines a class in $H_1(\SS_1^1,\p\SS_1^1\setminus\MM;\ZZ)$.
For clarity, we provide a picture in Figure~\ref{fig: 1} below.

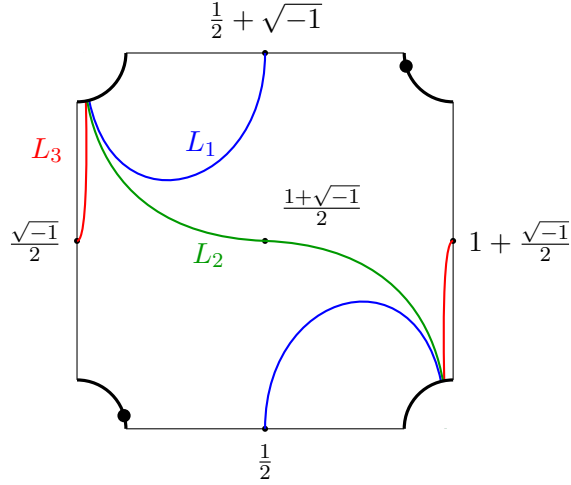
\begin{figure}[htbp]
  \centering
	\begin{tikzpicture}[scale=5,>=Stealth,
		myred/.style={line width=1.2pt,red!70!black},
		myblue/.style={line width=1.2pt,blue!70!black},
		mygreen/.style={line width=1.2pt,green!60!black},
		pole/.style={black, font=\scriptsize},
		tick/.style={fill=black, circle, inner sep=0.8pt}
	]
	% ===== 枠とラベル =====
	%\draw[black, line width=0.7pt] (0,0) rectangle (1,1);
	\draw (0.13,0) -- (0.87,0);
	\draw (0,0.13) -- (0,0.87);
	\draw (0.13,1) -- (0.87,1);
	\draw (1,0.13) -- (1,0.87);
	\coordinate (A) at (0.02,1);
	\coordinate (B) at (0,0);
	\coordinate (C) at (0.98,0);
	\coordinate (D) at (1,1);

	% 周期（極）
	\node[pole, anchor=north east] at (0,0) {};
	\node[pole, anchor=north west] at (1,0) {};
	\node[pole, anchor=south east] at (0,1) {};
	\node[pole, anchor=south west] at (1,1) {};

	% 半周期（臨界点）
	\coordinate (E) at (0.0,0.5);   % 1/2
	\coordinate (F) at (0.5,1);   % i/2
	\coordinate (G) at (1,0.5);   % 1/2
	\coordinate (H) at (0.5,0);   % i/2
	\coordinate (I) at (0.5,0.5);   % (1+i)/2

	\node[left=2pt] at (E) {$\tfrac{\sqrt{-1}}{2}$};
	\node[above right=2pt] at (I) {$\tfrac{1+\sqrt{-1}}{2}$};
	\node[right=2pt] at (G) {$1+\tfrac{\sqrt{-1}}{2}$};
	\node[above=2pt] at (F) {$\tfrac12+\sqrt{-1}$};
	\node[below=2pt] at (H) {$\tfrac12$};

	\node[tick] at (E) {};
	\node[tick] at (F) {};
	\node[tick] at (G) {};
	\node[tick] at (H) {};
	\node[tick] at (I) {}; 

	% ===== 逆像 =====

	% L_3
	\draw[red, thick]
		(E) .. controls +(0.04,0) and +(0,-0.0) .. (A);
	\draw[red, thick]
		(G) .. controls +(-0.04,0) and +(0,0) .. (C);

	% L_1
	\draw[blue, thick] 
		(F) .. controls +(-0.00,-0.4) and +(0,-0.5) .. (A);
	\draw[blue, thick] 
		(H) .. controls +(0.00,0.4) and +(0,0.5) .. (C);

	% L_2
	\draw[mygreen, thick] 
		(I) .. controls +(-0.53,0.02) and +(0,0) .. (A);
	\draw[mygreen, thick] 
		(I) .. controls +(0.53,-0.02) and +(0,0) .. (C);

	\node[red, anchor=north west]     at (-0.15,0.8) {$L_3$};
	\node[mygreen, anchor=south]         at (0.35,0.4) {$L_2$};
	\node[blue, anchor=south east] at (0.4,0.7) {$L_1$};

	% 四隅ごとにブローアップ（半径 r）
	\def\r{0.13}

	% 左下
	%\fill[white] (0,0) -- (\r,0) arc[start angle=0,end angle=90,radius=\r] -- cycle;
	%\draw (0,0) -- (\r,0) arc[start angle=0,end angle=90,radius=\r] -- cycle;
	\draw[black, very thick] (0.13,0) arc (0:90:0.13);
	% 右下
	\fill[white] (1,0) -- (1-\r,0) arc[start angle=180,end angle=90,radius=\r] -- cycle;
	%\draw (1,0) -- (1-\r,0) arc[start angle=180,end angle=90,radius=\r] -- cycle;
	\draw[black, very thick] (0.87,0) arc (180:90:0.13);

	% 左上
	\fill[white] (0,1) -- (0,1-\r) arc[start angle=270,end angle=365,radius=\r] -- cycle;
	%\draw (0,1) -- (0,1-\r) arc[start angle=270,end angle=360,radius=\r] -- cycle;
	\draw[black, very thick] (0,0.87) arc (270:360:0.13);

	% 右上
	%\fill[white] (1,1) -- (1,1-\r) arc[start angle=270,end angle=180,radius=\r] -- cycle;
	%\draw (1,1) -- (1,1-\r) arc[start angle=270,end angle=180,radius=\r] -- cycle;
	\draw[black, very thick] (0.87,1) arc (180:270:0.13);

	\coordinate (M1) at (0.125,0.035);
	\coordinate (M2) at (0.875,0.965);

	\fill (M1) circle (0.5pt);
	\fill (M2) circle (0.5pt);

	\end{tikzpicture}
	\caption{This shows how $L_3,L_2,L_1$ are arranged in $\SS_1^1$. The two points indicated by black dots on the boundary are the marked points.}
  \label{fig: 1}
\end{figure}

\begin{prop}[{\cite[Section~2.3]{LP}} and also {\cite[Section~3.4]{HKK}}]\label{HMS1}
Let the notation be as above.
There exists an equivalence of triangulated categories
\begin{equation}
\P\W(\SS_1^1, \MM, \ell)\cong \D^{b}{\rm mod} (\CC Q/I)
\end{equation}
which sends $L_i$ to $P(i)$ for $i=1,2,3$.
\end{prop}

\begin{pf}
Consider the sequence of arcs $\L:=(L_{3}, L_{2},L_{1})$ in $\SS_1^1$ (see Figure \ref{fig: 1}) 
and associate an $A_\infty$-category $\Fuk(\L)$ whose objects are $L_{3}$, $L_{2}$ and $L_{1}$.
Since $\ell$ gives the constant foliation, gradings of the arcs and morphisms between them can be chosen to be zero and
hence the space of morphisms $\Fuk(\L)(L_{i},L_{j})$ is the $\CC$-vector space with a basis given 
by counterclockwise paths inside $\p\SS_1^1\setminus \MM\cong S^1\setminus \{2\text{ pts}\}$ connecting $L_{i}$ to $L_{j}$ and the identity if $i=j$.
All $A_\infty$-operations except the product $m_{2}$ vanish 
since each polygon cut out of $\SS_1^1$ by the arcs contains exactly one element in $\MM$. 
The composition maps $m_{2}: \Fuk(\L)(L_{j},L_{k})\times \Fuk(\L)(L_{i},L_{j})\longrightarrow \Fuk(\L)(L_{i},L_{k})$ are defined 
through the concatenations of boundary paths.
The category $\P\W(\SS_1^1, \MM, \ell)$ is defined as the category of twisted complexes over $\Fuk(\L)$, 
which is independent of the particular choice of $\L$ as a triangulated category.
In particular, $\L$ becomes a full strong exceptional collection in $\P\W(\SS_1^1, \MM, \ell)$.
Calculating morphisms and their compositions according to the above rules, 
it turns out that the functor
\[
\P\W(\SS_1^1, \MM, \ell)\longrightarrow \D^{b}{\rm mod} (\CC Q/I), \quad L \mapsto 
{\rm Hom}_{\P\W(\SS_1^1, \MM, \ell)}(L_3\oplus L_2\oplus L_1,L)
\]
is an equivalence of triangulated categories.
\qed
\end{pf}

\begin{defn}[{\cite[Section~2, Section~7]{BD}}]
Let ${\bf E}\subset \PP^{2}$ be a nodal Weierstra{\ss} cubic curve over $\CC$ given by the equation $z_0 z_2^{2}=z_1^{3}+z_0z_1^{2}$ 
and $\mathcal{I}\subset \O_{\bf E}$ the ideal sheaf of the singular point $[0:0:1]\in {\bf E}$. 
Denote by $\nu: \PP^{1}\rightarrow {\bf E}$ the normalization of ${\bf E}$ where we choose coordinates on $\PP^{1}$ in such a way that
$\nu^{-1}([0:0:1])=\{\infty:=[1:0], 0:=[0:1]\}$.   
The Auslander sheaf of ${\bf E}$ is the sheaf of $\O_{\bf E}$-algebras $\A_{\bf E}:=\E nd_{\bf E}(\I \oplus \O_{\bf E})$.
Denote by ${\rm coh}(\A_{\bf E})$ the abelian category of coherent left $\A_{\bf E}$-modules
and by $\D^{b}{\rm coh}(\A_{\bf E})$ the bounded derived category of ${\rm coh}(\A_{\bf E})$.
\end{defn}

According to \cite[Section~7]{BD}, the derived category $\D^{b}{\rm coh}(\A_{\bf E})$ plays the role of the subcategory of compact objects of 
a categorical resolution of $\D^{b}{\rm coh}({\bf E})$ in the sense of Kuznetsov-Lunts~\cite[Definition~1.3]{KL}:
\begin{equation}\label{eq: cat resol}
\begin{tikzcd}[row sep=0.1cm]
& {\rm Perf}({\bf E}) \ar[ld, "\pi^*"'] \ar[dd, phantom, "\subseteq" sloped]\\
\D^{b}{\rm coh}(\A_{\bf E}) \ar[rd, "\pi_*"']\\
& \D^{b}{\rm coh}({\bf E})
\end{tikzcd}
\end{equation}
where the pair $(\pi^*, \pi_*)$ forms an adjoint pair of exact functors, which is given by
\[ 
\pi^* := (\I \oplus \O_{\bf E}) \overset{\mathbf{L}}{\otimes}_{\bf E} -,\quad
\pi_* := \mathbf{R}\H om_{\A_{\bf E}}(\I \oplus \O_{\bf E}, -).
\]
Furthermore, as shown in \cite[Theorem~2]{BD}, the functor $\pi^*$ is fully faithful.

Define $\M_3,\M_2,\M_1\in {\rm coh}(\A_{\bf E})$ by
\[
\M_3:=
\begin{pmatrix}
0\\
\O_{[0:0:1]}
\end{pmatrix},\quad 
\M_2:=
\begin{pmatrix}
\nu_{*}\O_{\PP^{1}}(-1)\\
\nu_{*}\O_{\PP^{1}}(-1)
\end{pmatrix},\quad
\M_1:=
\begin{pmatrix}
\nu_{*}\O_{\PP^{1}}\\
\nu_{*}\O_{\PP^{1}}
\end{pmatrix}.
\]

\begin{prop}[{\cite[Section~7]{BD}}]\label{prop: FSEC}
There exists an equivalence of triangulated categories
\begin{equation}
\D^{b}{\rm coh}(\A_{\bf E})\cong \D^{b}{\rm mod} (\CC Q/I)
\end{equation}
which sends $\M_i$ to $P(i)$ for $i=1,2,3$.
\end{prop}
\begin{pf}
The sequence $(\M_3,\M_2,\M_1)$ is a full strong exceptional collection in $\D^{b}{\rm coh}(\A_{\bf E})$,
whose morphism spaces are given by
\begin{align*}
&{\rm Hom}_{\D^{b}{\rm coh}(\A_{\bf E})} (\M_3,\M_2)=\CC x_{1,\infty}\oplus \CC x_{1,0}, \\
&{\rm Hom}_{\D^{b}{\rm coh}(\A_{\bf E})} (\M_2,\M_1)=\CC x_{2,\infty}\oplus \CC x_{2,0},\\
&{\rm Hom}_{\D^{b}{\rm coh}(\A_{\bf E})} (\M_3,\M_1)
=\CC x_{2,\infty}\circ x_{1,\infty}\oplus \CC x_{2,0}\circ x_{1,0},
\end{align*}
such that $x_{2,\infty}\circ x_{1,0}=x_{2,0}\circ x_{1,\infty}=0$.
In particular, the endomorphism algebra ${\rm End}_{\D^{b}{\rm coh}(\A_{\bf E})}(\M_3\oplus \M_2\oplus \M_1)$ 
is isomorphic to the path algebra $\CC Q/I$ of the nodal quiver $(Q,I)$ and the functor
\[
\D^{b}{\rm coh}(\A_{\bf E})\cong \D^{b}{\rm mod} (\CC Q/I),\quad \M\mapsto 
{\rm Hom}_{\D^{b}{\rm coh}(\A_{\bf E})} (\M_3\oplus\M_2\oplus\M_1,\M)
\]
is an equivalence of triangulated categories.
\qed\end{pf}

\begin{cor}\label{HMS2}
There exist triangulated equivalences:
\begin{equation}
\D^{b}{\rm coh}(\A_{\bf E})\cong \D^{b}{\rm mod} (\CC Q/I)\cong \P\W(\SS_1^1, \MM, \ell),\quad 
\M_i\mapsto P(i)\mapsto L_i.
\end{equation}
\qed\end{cor}

We will continue to denote $\D^{b}{\rm mod} (\CC Q/I)$ by $\D$ to simplify the notation,
though due to the categorical equivalence discussed above, one may use $\D^{b}{\rm coh}(\A_{\bf E})$ or $\P\W(\SS_1^1, \MM, \ell)$ as $\D$.

\subsection{Braid group action on $\FEC(\D)$}\label{subsec: Braid group action}
Denote by $\FEC(\D)$ the set of isomorphism classes of full exceptional collections in $\D$ and 
by $\Aut(\D)$ the group of auto-equivalences of $\D$.
Note that the group $\Aut(\D)$ acts on $\FEC(\D)$ from the left by
\[
\Phi \cdot \E := (\Phi(E_1),\Phi(E_2),\Phi(E_3)), \quad \Phi \in \Aut(\D),\ \E=(E_1,E_2,E_3)\in\FEC(\D).
\]

The Artin {\it braid group} $B_3$ on $3$-strands is a group presented by the following generators and relations: 
\begin{equation*}
B_3:=\langle {\boldsymbol\sigma}_1,{\boldsymbol\sigma}_2\,|\,{\boldsymbol\sigma}_{1}{\boldsymbol\sigma}_{2}{\boldsymbol\sigma}_{1}={\boldsymbol\sigma}_{2}{\boldsymbol\sigma}_{1}{\boldsymbol\sigma}_{2}\rangle.
\end{equation*}
The groups $B_3$ and $\ZZ^3$ act on $\FEC(\D)$ from the right by mutations \cite{BP}:
\begin{gather*}
\E \cdot {\boldsymbol\sigma}_1 := (\L_{E_1}E_2,E_1,E_3),\quad \E \cdot {\boldsymbol\sigma}_2 := (E_1,\L_{E_2}E_3,E_2),\\
\E \cdot {\boldsymbol\sigma}_1^{-1} := (E_2,\R_{E_2}E_1,E_3),\quad \E \cdot {\boldsymbol\sigma}_2^{-1} := (E_1,E_3,\R_{E_3}E_2),\\
\E\cdot (n_1,n_2,n_3):=(E_1[n_1],E_2[n_2],E_3[n_3]),\quad (n_1,n_2,n_3)\in\ZZ^3,
\end{gather*}
where $\L_{E_i}E_{i+1},~\R_{E_{i+1}}E_{i}$ for $i=1,2$ are defined by the following exact triangles
\begin{gather*}
\L_{E_{i}}E_{i+1}\longrightarrow \RHom_\D(E_{i},E_{i+1})\otimes E_{i}\longrightarrow E_{i+1}\longrightarrow \L_{E_{i}}E_{i+1}[1],\\
\R_{E_{i+1}}E_{i}[-1]\longrightarrow E_{i}\longrightarrow \RHom_\D(E_{i},E_{i+1})^*\otimes E_{i+1}\longrightarrow \R_{E_{i+1}}E_{i}.
\end{gather*}
Note that on $\FEC(\D)$ the left action of $\Aut(\D)$ and the right action of $B_3$ commute, namely, $(\Phi\cdot\E)\cdot {\boldsymbol\sigma}=\Phi\cdot(\E\cdot {\boldsymbol\sigma})$ for all $\Phi\in\Aut(\D)$ and ${\boldsymbol\sigma}\in B_3$.
If $\E\in\FEC(\D)$ is strong, the actions of $\Aut(\D)$ and $B_3$ preserve the property of being strong.

\begin{prop}\label{prop: fsec}
Let $\E=(E_1,E_2,E_3)$ be a full strongly exceptional collection in $\D$. 
The endomorphism algebra ${\rm End}_\D(\E):={\rm Hom}_\D(A,A)$, $A:=E_1\oplus E_2\oplus E_3$ is isomorphic to the path algebra $\CC Q/I$ of the nodal quiver $(Q,I)$.
\end{prop}
\begin{pf}
Since $\E$ is a full strongly exceptional collection, the endomorphism algebra ${\rm End}_\D(\E)$ is derived equivalent to the path algebra $\CC Q/I$ of the nodal quiver $(Q,I)$. 
Since $\CC Q/I$ is a finite-dimensional gentle algebra, ${\rm End}_\D(\E)$ is also a gentle algebra by \cite[Corollary~1.2]{SZ}.

Suppose that $[E_i]=l_i[P(3)]+m_i[P(2)]+n_i[P(1)]$ for some $l_i,m_i,n_i\in\ZZ$ in the Grothendieck group $K_0(\D)$.
Since $E_i$ is exceptional, we have $\chi_\D([E_i],[E_i])=(l_i+m_i+n_i)^2=1$ and hence $l_i+m_i+n_i=\pm 1$. 
Since the full exceptional collection $(E_1,E_2,E_3)$ is strong, $\chi_\D([E_i],[E_j])=\dim\Hom_\D(E_i,E_j)\geq 0$ and $\chi_\D([E_j],[E_i])=0$ if $i<j$. Therefore, for $i<j$ we have
\begin{eqnarray*}
\dim\Hom_\D(E_i,E_j) &=& \chi_\D([E_i],[E_j])+\chi_\D([E_j],[E_i])\\
&=& (l_il_j+m_im_j+n_in_j+2l_im_j+2l_in_j+2m_in_j)\\
& & +(l_il_j+m_im_j+n_in_j+2m_il_j+2n_il_j+2n_im_j)\\
&=& 2(l_i+m_i+n_i)(l_j+m_j+n_j)=2.
\end{eqnarray*}
Since ${\rm End}_\D(\E)$ is gentle, this calculation implies that ${\rm End}_\D(\E)$ is isomorphic to $\CC Q/I$.
\qed
\end{pf}

\begin{defn}[{\cite[Definition~2.14]{ST}}]
An object $S\in \D$ is \textit{$m$-spherical} if the following holds
\begin{enumerate}
\item $\RHom_\D(S,S)\cong\CC\oplus\CC[-m]$,
\item For the Serre functor $\S_\D$ of $\D$, $\S_\D(S)\cong S[m]$.
\end{enumerate}
\end{defn}
A spherical object gives an auto-equivalence as follows.
\begin{prop}[{\cite[Theorem~2.10]{ST}}]
Let $S\in\D$ be an $m$-spherical object. There exists an auto-equivalence $T_{S}\in\Aut(\D)$ given by the following exact triangle
\[
T_{S}(X)[-1]\longrightarrow \RHom_\D(S,X)\otimes S\longrightarrow X\longrightarrow T_{S}(X),\quad X\in\D.
\]
In particular, $T_S$ induces an automorphism $[T_S]\in\Aut_\ZZ(K_0(\D),\chi_\D)$ satisfying 
\[
[T_{S}(X)]=[X]-\chi_\D([S],[X])[S],\ X\in \D.
\]
\qed
\end{prop}

\begin{prop}
Let $\E=(E_1,E_2,E_3)$ be the full strongly exceptional collection $(P(3),P(2),P(1))$ and set 
\[
c_i:=c_{i,\infty}x_{1,\infty}+c_{i,0} x_{1,0}\in \Hom_\D(E_{i}, E_{i+1}),\quad i=1,2,
\]
for $(c_{i,\infty},c_{i,0})\in\CC^2\setminus ((\{0\}\times \CC)\cup(\CC\times \{0\}))$.
An object $S_{i,c_i}\in\D$ defined by the following exact triangle
\begin{equation}\label{1-sph}
E_{i}\overset{c_i}{\longrightarrow} E_{i+1}\longrightarrow S_{i,c_i}\longrightarrow E_{i}[1]
\end{equation}
is $1$-spherical. Moreover, we have 
\begin{equation}\label{eq: hom s1-s2}
\RHom_\D(S_{1,c_1},S_{2,c_2})\cong \CC.
\end{equation}
\end{prop}
\begin{pf}
It follows from direct calculation using the hom-functors.
\qed
\end{pf}
Note that $S_{i,c_i}\cong S_{i,c_i'}$ if $c_i'=c\cdot c_i$ for some $c\in\CC^*$.

Set $S_i := S_{i,1}$ by choosing $(c_{i,\infty}, c_{i,0}) = (1, 1)$.
By \eqref{eq: hom s1-s2}, we have
\begin{cor}\label{cor: left B_3}
The correspondence ${\boldsymbol\sigma}_i\mapsto T_{S_{i}}^{-1}$, $i=1,2$ yields a surjective group homomorphism
\[
B_3\longrightarrow \langle T_{S_{1}}^{-1}, T_{S_{2}}^{-1} \rangle\subset \Aut(\D).
\]
\end{cor}
\begin{pf}
It is a direct consequence of Seidel--Thomas \cite[Theorem 2.17]{ST}. 
\qed
\end{pf}
The right $B_3$-action on $\FEC(\D)$ can be represented by left $\Aut(\D)$-action 
in the following way.
\begin{prop}\label{prop: left right B_3}
Let $\E=(E_1,E_2,E_3)$ be the full strongly exceptional collection $(P(3),P(2),P(1))$.
In $\FEC(\D)$, we have 
\[
T_{S_{i,c_i}}\cdot\E =\E\cdot{\boldsymbol\sigma}_i^{-1}.
\]
\end{prop}
\begin{pf}
We give a proof only for the case when $i=1$ since the case when $i=2$ can be shown in a similar manner.
Note that $\E\cdot{\boldsymbol\sigma}_1^{-1}=(E_2,\R_{E_2}E_1,E_3)$. 
By applying the functor $\Hom_\D(-,E_k)$ to (\ref{1-sph}), we have 
\[
\RHom_\D(S_{1,c_1},E_k)\cong
\begin{cases}
\CC[-1] & \text{if}\quad k=1,2,\\
0 & \text{if}\quad k=3.
\end{cases}
\]
Therefore, 
$T_{S_{1,c_1}}(E_1) \cong \Cone(S_{1,c_1}[-1] \longrightarrow E_1)\cong E_2$ and $T_{S_{1,c_1}}(E_3)\cong E_3$.
On the other hand, \eqref{1-sph} and the isomorphism $\R_{E_{2}} E_1\cong\Cone(E_1\to E_{2}^{\oplus 2})$ show that $\RHom_\D(T_{S_{1,c_1}}(E_{2}),\R_{E_{2}} E_1)\cong\CC$.
Since both $T_{S_{1,c_1}}\cdot\E$ and $\E\cdot{\boldsymbol\sigma}_1^{-1}$ are full exceptional collections, the above computation yields $T_{S_{1,c_1}}(E_{2})\cong\R_{E_{2}} E_1$. 
Thus we have $T_{S_{1,c_1}}\cdot\E=\E\cdot{\boldsymbol\sigma}_1^{-1}$ in $\FEC(\D)$.
\qed
\end{pf}

Define elements $\alpha, \delta_1,\delta_2\in K_0(\D)$ by  
\begin{gather*}
\alpha:=[E_+]=[P(3)]-[P(2)]+[P(1)],\\ 
\delta_1:=[S_{1}]=[P(2)]-[P(3)],\quad
\delta_2:=[S_{2}]=[P(1)]-[P(2)].
\end{gather*}
With respect to  the $\ZZ$-basis $\{\alpha,\delta_1,\delta_2\}$ of $K_0(\D)$, the Euler form $\chi_\D$ is represented by
\[
\chi_\D=
\begin{pmatrix}
1 & 0 & 0\\
0 & 0 & 1\\
0 & -1 & 0
\end{pmatrix}.
\]
The definition of the Serre functor implies that $[\S_\D^2]=\id_{K_0(\D)}$ since 
\[
[\S_\D]=\chi_\D^{-1}\chi_\D^T=
\begin{pmatrix}
1 & 0 & 0\\
0 & 0 & -1\\
0 & 1 & 0
\end{pmatrix}
\begin{pmatrix}
1 & 0 & 0\\
0 & 0 & -1\\
0 & 1 & 0
\end{pmatrix}
=
\begin{pmatrix}
1 & 0 & 0\\
0 & -1 & 0\\
0 & 0 & -1
\end{pmatrix}.
\]
With respect to  the basis $\{\alpha,\delta_1,\delta_2\}$ of $K_0(\D)$, 
the automorphisms induced by $T_{S_{1}}^{-1}$ and $T_{S_{2}}^{-1}$ 
are represented by the following matrices 
\begin{equation}\label{eq: SL2}
[T_{S_{1}}^{-1}]=
\begin{pmatrix}
1 & 0 & 0\\
0 & 1 & 1\\
0 & 0 & 1
\end{pmatrix},
\quad
[T_{S_{2}}^{-1}]=
\begin{pmatrix}
1 & 0 & 0\\
0 & 1 & 0\\
0 & -1 & 1
\end{pmatrix}.
\end{equation}
Hence, restricting these automorphisms 
to $\ZZ\delta_1\oplus\ZZ\delta_2$, we have a group homomorphism 
\[
\langle T_{S_{1}}^{-1}, T_{S_{2}}^{-1} \rangle\longrightarrow 
\Aut_\ZZ(K_0(\D),\chi_\D)\longrightarrow 
\Aut_\ZZ(\ZZ\delta_1\oplus\ZZ\delta_2,\chi_\D|_{\ZZ\delta_1\oplus\ZZ\delta_2})={\rm SL}(2;\ZZ),
\]
which is surjective.
Thus we have the following commutative diagram of short exact sequences:
\[
\begin{tikzcd}
\{1\} \arrow[r] & \langle ({\boldsymbol\sigma}_{1}{\boldsymbol\sigma}_{2})^6 \rangle \arrow[r] \arrow[d] & B_3 \arrow[r] \arrow[d, "\text{Cor.~\ref{cor: left B_3}}"] & {\rm SL}(2;\ZZ) \arrow[r] \arrow[d, equal] & \{1\}\\
\{1\} \arrow[r] & \langle (T_{S_{1}}T_{S_{2}})^{-6} \rangle \arrow[r] & \langle T_{S_{1}}^{-1},T_{S_{2}}^{-1} \rangle \arrow[r] & {\rm SL}(2;\ZZ) \arrow[r] & \{1\}\\
\end{tikzcd}.
\]
\begin{prop}
The group homomorphism $B_3\longrightarrow \langle T_{S_{1}}^{-1},T_{S_{2}}^{-1}\rangle$
is an isomorphism.
\end{prop}
\begin{pf}
Note that $\E\cdot ({\boldsymbol\sigma}_{2}{\boldsymbol\sigma}_{1})^{3}\cong\S_\D[-2]\cdot \E$ for any $\E\in \FEC(\D)$ (cf.~\cite[Assertion 4.2]{Bo}) and hence $\S_\D(P(k))[-2]\cong (T_{S_{2}}^{-1}T_{S_{1}}^{-1})^3(P(k))\cong (T_{S_{1}}T_{S_{2}})^{-3}(P(k))$. 
Since 
\begin{eqnarray*}
 (T_{S_{1}}T_{S_{2}})^3(S_{1,c_1}) 
&\cong& \S_\D^{-1}(\Cone(\S_\D (T_{S_{1}}T_{S_{2}})^3(P(3))[-2] \overset{c_1'}{\longrightarrow} \S_\D  (T_{S_{1}}T_{S_{2}})^3(P(2))[-2]))[2]\\
&\cong& \S_\D^{-1}(S_{1,c_1'})[2]\cong S_{1,c_1'}[1]\ncong S_{1},
\end{eqnarray*}
where $c_1':=\S_\D(T_{S_{1}}T_{S_{2}})^3 (c_1)[-2]\in \Hom_\D(P(3), P(2))$, the map $\langle ({\boldsymbol\sigma}_{1}{\boldsymbol\sigma}_{2})^6 \rangle \longrightarrow \langle (T_{S_{1}}T_{S_{2}})^{-6}\rangle$ is an isomorphism.
Thus we have $B_3\cong \langle T_{S_{1}}^{-1},T_{S_{2}}^{-1} \rangle$.
\qed\end{pf}

\begin{defn}[{Broomhead--Pauksztello--Ploog~\cite[Definition~4.2]{BPP}}]
For any integer $n\ge 2$, a sequence $(E_1,\dots,E_n)$ of objects of $\D$ is an \textit{exceptional $n$-cycle} if 
\begin{enumerate}
\item $E_i$ is exceptional for all $i=1,\dots, n$,
\item there are integers $k_i$ such that $\S_\D(E_i)\cong E_{i+1}[k_i]$ for all $i=1,\dots, n$ where $E_{n+1}:=E_1$,
\item $\RHom_\D(E_i,E_j)=0$ unless $j=i$ or $i+1$.
\end{enumerate}
\end{defn}
If $(E_1,\dots,E_n)$ is an exceptional $n$-cycle, so is $(E_2,\dots,E_n,E_1)$. 

A spherical object $S$ can be thought of as an exceptional $1$-cycle.
As the following proposition shows, the above definition provides a natural generalization of spherical objects and their twist
functors.
\begin{prop}[{\cite[Theorem~4.5]{BPP}}]
For any integer $n\ge 2$ and an exceptional $n$-cycle $E_*=(E_1,\dots,E_n)$,
the functor $T_{E_*}:\D\longrightarrow \D$ defined by the following exact triangle
\[
T_{E_*}(X)[-1]\longrightarrow \bigoplus_{i=1}^n\RHom_\D(E_{i},X)\otimes E_{i}\longrightarrow X\longrightarrow T_{E_*}(X)
\]
is an equivalence.
\qed
\end{prop}
The above functor $T_{E_*}$ is called the {\it twist functor} associated to the $n$-cycle $E_*=(E_1,\dots,E_n)$, 
which induces an automorphism $[T_{E_*}]\in\Aut_\ZZ(K_0(\D),\chi_\D)$ given by
\[
[T_{E_*}]([X])=[X]-\sum_{i=1}^{n}\chi_\D([E_i],[X])[E_i].
\]

\begin{defn}\label{defn: B(D)}
Define a group $B(\D)$ as follows:
\begin{equation*}
B(\D):=\langle T_{E_{*}}^n\,|\, E_{*}=(E_{1},\dots, E_{n});\ \text{an exceptional}\  \text{$n$-cycle}, n\ge 1  \rangle.
\end{equation*}
\end{defn}
As already mentioned, by an exceptional $1$-cycle we mean an $m$-spherical object for some $m\in\ZZ$.  
For any $\CC$-linear triangulated category $\T$, we denote by ${\rm Sph}_m(\T)$ the set of $m$-spherical objects in $\T$.
\begin{rem}
In \cite{OST}, the group generated by the twist functors associated to $1$-spherical objects has been considered. 
However, in general situations, including the one treated in this paper, it is more natural to consider the group $\B(\D)$ defined above. 
Note that in the case of \cite{OST}, these groups coincide.
\end{rem}

In order to understand $B(\D)$, we begin by classifying exceptional cycles.
\begin{prop}[{See also \cite[Proposition~3.8]{O2}}]\label{prop: exceptional cycle}
Let $n\ge 2$. 
A sequence $E_*=(E_1,\dots,E_n)$ of objects of $\D$ is an exceptional $n$-cycle if and only if 
$n=2$ and $E_*=(E_+,E_-)$ or $(E_-,E_+)$ up to translations.

In particular, we have
\[
T_{(E_{+},E_{-})}\cong \S_\D^{-1}[1].
\]
\end{prop}
\begin{pf}
We begin by proving the sufficiency.
\begin{lem}
The sequence $(E_{+},E_{-})$ and $(E_{-},E_{+})$ are exceptional $2$-cycles.
\end{lem}
\begin{pf}
We only need to show that $(E_{+},E_{-})$ is an exceptional $2$-cycle.
Lemma~\ref{lem: sung} shows that $E_{+},E_{-}$ are exceptional objects satisfying $\S_\D(E_\pm)\cong E_\mp[2]$.
The third condition is automatically satisfied since $n=2$. 
\qed
\end{pf}
\begin{lem}\label{lem: exceptional-serre}
Let $T_{(E_{+},E_{-})}$ be the twist functor associated to the $2$-cycle $(E_{+},E_{-})$. 
We have 
\[
T_{(E_{+},E_{-})}\cong \S_\D^{-1}[1].
\]
\end{lem}
\begin{pf}
By a direct calculation based on the definitions of the twist functor and left mutation functor (see \cite[Remark 2.8]{Ku}), 
it is easy to see that $T_{(E_+,E_-)}(E_{+})\cong E_{-}[-1]\cong \S_{\D}^{-1}(E_{+})[1]$
and $\left. T_{(E_{+},E_{-})}\right|_{E_{+}^{\perp}}=\LL_{\S_{\D}(E_{+})}$.
It turns out that $T_{(E_{+},E_{-})}\cong \S_{\D}^{-1}[1]$ by the argument of \cite[Section 3]{Su}.
\qed
\end{pf}
\begin{lem}
If an exceptional $n$-cycle exists in $\D$, then $n=2$.
\end{lem}
\begin{pf}
Without loss of generality, we may assume that an exceptional $n$-cycle of the form $E_*=(E,\S_\D(E),\dots, \S_\D^{n-1}(E))$ exists.
Since $E$ is exceptional and $[\S_\D^2]=\id_{K_0(\D)}$, we have 
\[
[T_{E_*}]([\S_\D^{j}(E)])=[\S_\D^{j}(E)]-\sum_{i=0}^{n-1}\chi([\S_\D^{i}(E)],[\S_\D^{j}(E)])[\S_\D^{i}(E)]
=[\S_\D^{j}(E)]-\sum_{i=0}^{n-1}[\S_\D^{i}(E)].
\]
Therefore,
\[
[T_{E_*}]([E]+[\S_\D^{1}(E)]+\dots+[\S_\D^{n-1}(E)])=(1-n)([E]+[\S_\D^{1}(E)]+\dots+[\S_\D^{n-1}(E)]).
\]
It follows that $[T_{E_*}]$ is an automorphism only if $1-n=\pm 1$ and hence $n=2$.
\qed\end{pf}

\begin{lem}
Suppose that $(E,\S_\D(E))$ is an exceptional $2$-cycle. Then $E$ is isomorphic to $E_+$ or $E_-$ up to translation.
\end{lem}
\begin{pf}
Since $[\S_\D^2]=\id_{K_0(\D)}$, there exists an integer $k$ such that $\S_\D^2(E) \cong E[2k]$. 
By definition of $T_{(E_{+},E_{-})}$, we have the following exact triangle 
\[
T_{(E_{+},E_{-})}(E)[-1]\longrightarrow \RHom_\D(E_+,E)\otimes E_+\oplus \RHom_\D(E_-,E)\otimes E_-\longrightarrow E\longrightarrow T_{(E_{+},E_{-})}(E).
\]
By Lemma~\ref{lem: exceptional-serre}, we have $\Hom_\D(E,T_{(E_{+},E_{-})}(E))\cong\Hom_\D(E,\S_\D^{-1}(E)[1])$ and 
\[
\Hom_\D(E,\S_\D^{-1}(E)[1])\cong\Hom_\D(\S_\D^{-2}(E)[1],E)^*\cong\Hom_\D(E[1-2k],E)^*=0,
\]
which implies that the above exact triangle splits. 
It then follows from the Krull--Schmidt property that $E\cong E_+$ or $E\cong E_-$ up to translations.
\qed\end{pf}
Thus we have finished the proof of the proposition.
\qed
\end{pf}

Next, we consider spherical objects.
\begin{prop}
Suppose $S$ is an $m$-spherical object in $\D$. Then $m=1$.
\end{prop}
\begin{pf}
By definition of $T_{(E_{+},E_{-})}$ and Lemma~\ref{lem: exceptional-serre}, we have the following exact triangle 
\[
S[-m]\longrightarrow \RHom_\D(E_+,S)\otimes E_+\oplus \RHom_\D(E_-,S)\otimes E_-\longrightarrow S\longrightarrow S[1-m].
\]
If $m\ne 1$, then the morphism $S\longrightarrow S[1-m]$ is zero; then $S\cong E_+$ or $E_-$ up to translations due to the Krull-Schmidt property of $\D$, which contradicts the fact that $\S_\D(E_\pm)\cong E_\mp[2]$. 
\qed\end{pf}

\begin{prop}
We have 
\[
B(\D)=\langle T_{(E_+,E_-)}^2,\ T_{S_1,c_1},\ T_{S_2,c_2}\,|\, c_1,c_2\in\CC^2\setminus ((\{0\}\times \CC)\cup(\CC\times \{0\}))\rangle.
\]
\end{prop}
\begin{pf}
The set ${\rm Sph}_1(\D)$ of $1$-spherical objects in $\D$ is characterized as follows:
\begin{lem}[{\cite[Corollary~2]{BD}} and also {\cite[Lemma~4.1]{O2}}]
Let ${\rm Perf}({\bf E})$ be the derived category of perfect complexes of sheaves on ${\bf E}$. 
The fully faithful functor $\pi^*:{\rm Perf}({\bf E})\to\D$ in the diagram \textnormal{(\ref{eq: cat resol})} induces a bijection between the sets of $1$-spherical objects:
\[
{\rm Sph}_1(\D)=\pi^*({\rm Sph}_1({\rm Perf}({\bf E}))).
\]
\qed\end{lem}
It follows from \cite[Propositions~11]{BD} that for any smooth point 
$\nu([c_{1,\infty}:c_{1,0}])\in{\bf E}$ with $[c_{1,\infty}:c_{1,0}]\in \PP^1\setminus\{\infty,0\}$, 
we have 
\[
\pi^*(\CC_{\nu([c_{1,\infty}:c_{1,0}])}) \cong S_{1,(c_{1,\infty},c_{1,0})}.
\]
It is known that $\Pic^0({\bf E})\cong \CC^*$ and by \cite[Propositions~11 and 13]{BD}
the element $\L_{[c_{1,\infty}':c_{1,0}']}\in\Pic^0({\bf E})$ corresponding to $[c_{1,\infty}':c_{1,0}']\in\PP^1\setminus\{\infty,0\}=\CC^*$
is mapped by $\pi^*$ to the complex of $\CC Q/I$-modules 
\[
P(3)\xrightarrow{
\begin{psmallmatrix}
c_{1,\infty}' x_{1,\infty} \\
c_{1,0}'x_{1,0}
\end{psmallmatrix}
} 
{P(2)\oplus P(2)} \xrightarrow{
\begin{psmallmatrix}
x_{2,\infty} & x_{2,0}
\end{psmallmatrix}
} 
P(1),
\]
where the middle term is of degree zero.
Furthermore, a direct computation shows that
\[
\pi^*(\L_{[c_{1,\infty}':c_{1,0}']})\cong T_{S_{1,c_1}}S_{2,c_2}[-1]
\quad\text{for some } c_1, c_2\in\CC^2\setminus ((\{0\}\times \CC)\cup(\CC\times \{0\})).
\]
\begin{lem}[{\cite[Proposition~4.13, Proposition~4.2]{BuKr}}]
Let the notations be as above.
Consider the subgroup of $\Aut({\rm Perf}({\bf E}))$ given by
\[
\langle T_{\CC_{\nu([c_{1,\infty}:c_{1,0}])}}, T_{\L_{[c_{1,\infty}':c_{1,0}']}}\,|\, [c_{1,\infty}:c_{1,0}], [c_{1,\infty}':c_{1,0}']\in\PP^1\setminus\{\infty,0\}\rangle.
\]
This subgroup acts transitively on the set ${\rm Sph}_1({\rm Perf}({\bf E}))$ up to translations.
\qed\end{lem}
Putting these together, we obtain the statement of the proposition.
\qed\end{pf}

Let $\Aut_{\bf Set}(\FEC(\D)/\ZZ^3)$ be the automorphism group of the
set $\FEC(\D)/\ZZ^3$ (modulo translations)
and consider the following group homomorphism:
\begin{equation}
\rho_{\FEC}: \Aut(\D)\longrightarrow \Aut_{\bf Set}(\FEC(\D)/\ZZ^3).
\end{equation}

\begin{prop}[cf. {\cite[Theorem~C]{O1}}]\label{prop: alg Aut}
Let ${\rm Out}_{\CC\text{-}{\bf alg}}(\CC Q/I)$ be the group of outer
automorphisms of $\CC Q/I$ defined  by
\[
{\rm Out}_{\CC\text{-}{\bf alg}}(\CC Q/I):=\Aut_{\CC\text{-}{\bf
alg}}(\CC Q/I)/{\rm Inn}_{\CC}(\CC Q/I)
\]
where $\Aut_{\CC\text{-}{\bf alg}}(\CC Q/I)$ is the automorphism group
of the $\CC$-algebra $\CC Q/I$ and ${\rm Inn}_{\CC}(\CC Q/I)$ is the
normal subgroup of
$\Aut_{\CC\text{-}{\bf alg}}(\CC Q/I)$ consisting of inner
automorphisms.
\begin{enumerate}
\item
We have
\[
{\rm Out}_{\CC\text{-}{\bf alg}}(\CC Q/I)\cong \left(\CC^*\times
\CC^*\right)\rtimes (\ZZ/2\ZZ).
\]
\item
We have the following short exact sequence of groups
\begin{equation}
\{1\}\longrightarrow {\rm Out}_{\CC\text{-}{\bf alg}}(\CC
Q/I)\times\ZZ[1]\longrightarrow\Aut(\D)
\overset{\rho_{\FEC}}{\longrightarrow} B_3\longrightarrow\{1\},
\end{equation}
where ${\rm Out}_{\CC}(\CC Q/I)$ is naturally identified with the
subgroup
\[
\{\Phi\in \Aut(\D)\,|\,\Phi({\rm mod}(\CC Q/I))={\rm mod}(\CC Q/I),\
\Phi(P(i))\cong P(i), i=1,2,3\}.
\]

\end{enumerate}
\end{prop}
\begin{pf}
(i) Since $\CC Q/I$ is a basic finite dimensional $\CC$-algebra, any
$\CC$-algebra automorphism maps a complete set of primitive idempotents
to a complete set of primitive idempotents and, up to inner
automorphisms, fixes them.
Therefore, in determining the outer automorphism group ${\rm
Out}_{\CC\text{-}{\bf alg}}(\CC Q/I)$,
we may assume that any $\CC$-algebra automorphism $\phi \in
\Aut_{\CC\text{-}{\bf alg}}(\CC Q/I)$ fixes the primitive idempotents
$e_1, e_2, e_3$ associated to vertices $1,2,3$ of $Q$,
which is determined by its action on the arrows $\alpha_1, \alpha_2,
\beta_1, \beta_2$. Hence we may write
\begin{gather*}
\phi(\alpha_1) =\phi^1_{11} \alpha_1 +\phi^1_{12} \beta_1,\quad
\phi(\beta_1) = \phi^1_{21} \alpha_1 + \phi^1_{22} \beta_1,\\
\phi(\alpha_2) = \phi^2_{11} \alpha_2 + \phi^2_{12} \beta_2,\quad
\phi(\beta_2) = \phi^2_{21} \alpha_2 + \phi^2_{22} \beta_2,
\end{gather*}
for some $((\phi^1_{ij}),(\phi^2_{ij})) \in {\rm GL}(2;\CC)^2$, which
must satisfy
\[
\phi^1_{11} \phi^2_{21} = \phi^1_{12} \phi^2_{22} = \phi^1_{21}
\phi^2_{11} = \phi^1_{22} \phi^2_{12} = 0.
\]
Solving these equations yields
\[
\phi^1_{12}=\phi^1_{21}=\phi^2_{12}=\phi^2_{21}=0\quad\text{or}\quad
\phi^1_{11}=\phi^1_{22}=\phi^2_{11}=\phi^2_{22}=0.
\]
Thus, the subgroup of $\Aut_{\CC\text{-}{\bf alg}}(\CC Q/I)$ fixing
$e_1,e_2,e_3$ is generated by:
\begin{itemize}
\item[(a)]
the subgroup generated by the scaling transformations of the arrows
$\alpha_1,\beta_1,\alpha_2,\beta_2$ given by
$\{((\phi^1_{11},\phi^1_{22}),(\phi^2_{11}, \phi^2_{22}))\}\cong
(\CC^*)^2\times (\CC^*)^2$,
\item[(b)] an involution that exchanges $\alpha_1$ with $\beta_1$ and
$\alpha_2$ with $\beta_2$, while preserving the relations.
\end{itemize}
Note that the involution exchanges $\phi^1_{11}$ with $\phi^1_{22}$ and
$\phi^2_{11}$ with $\phi^2_{22}$.
As a result, the group is a nontrivial semidirect product
\[
\left((\CC^*)^2\times (\CC^*)^2\right)\rtimes (\ZZ/2\ZZ).
\]

Any inner automorphism fixing $e_1,e_2,e_3$ is given by the element of
the form
\[
u = \sum_{i=1}^3 \psi_i e_i,\quad \psi_i \in \CC^*,
\]
which acts on the arrows as
\[
\alpha_1 \mapsto \frac{\psi_2}{\psi_1} \alpha_1,\quad
\beta_1 \mapsto \frac{\psi_2}{\psi_1} \beta_1,\quad
\alpha_2 \mapsto \frac{\psi_3}{\psi_2} \alpha_2,\quad
\beta_2 \mapsto \frac{\psi_3}{\psi_2} \beta_2.
\]
It is obvious that this action is trivial if and only if
$\psi_1=\psi_2=\psi_3$. Hence, the image of the subgroup under the
natural map
${\rm Inn}_{\CC\text{-}{\bf alg}}(\CC Q/I)\longrightarrow
\Aut_{\CC\text{-}{\bf alg}}(\CC Q/I)$ is $\CC^*\times \CC^*\subset
(\CC^*)^2\times (\CC^*)^2$,
the product of the diagonals of each $(\CC^*)^2$ factor in
$(\CC^*)^2\times (\CC^*)^2$.

(ii) The statement follows from \cite[Proposition~5.3]{O2}. 
Since $\CC Q/I=P(3)\oplus P(2)\oplus P(1)$ in ${\rm mod}(\CC Q/I)$,
the identification of ${\rm Out}_{\CC\text{-}{\bf alg}}(\CC Q/I)$
follows from \cite[Proposition 2.4]{Ch} and references therein.
Indeed, it is clear that $T_{S_1,c_1}T_{S_1}^{-1}$,
$T_{S_2,c_2}T_{S_2}^{-1}$ belong to ${\rm Ker}(\rho_{\FEC})$ and they
respect $P(3)$, $P(2)$ and $P(1)$, and a direct calculation yields
\begin{gather*}
T_{S_1,c_1}T_{S_1}^{-1}(x_{1,\infty})=x_{1,\infty}, \
T_{S_1,c_1}T_{S_1}^{-1}(x_{1,0})=x_{1,0}, \\
  T_{S_1,c_1}T_{S_1}^{-1}(x_{2,\infty})=c_{1,\infty}x_{2,\infty}, \
T_{S_1,c_1}T_{S_1}^{-1}(x_{2,0})=c_{1,0}x_{2,0}, \\
T_{S_2,c_2}T_{S_2}^{-1}(x_{2,\infty})=c_{2,\infty}x_{1,\infty}, \
T_{S_2,c_2}T_{S_2}^{-1}(x_{1,0})=c_{2,0}x_{1,0}, \\
  T_{S_2,c_2}T_{S_2}^{-1}(x_{2,\infty})=x_{2,\infty}, \
T_{S_2,c_2}T_{S_2}^{-1}(x_{2,0})=x_{2,0},
\end{gather*}
which generate the subgroup $\CC^*\times \CC^*$ of ${\rm
Out}_{\CC\text{-}{\bf alg}}(\CC Q/I)$.
Note also that $(\S_\D[-2])\circ (T_{S_{1}}T_{S_{2}})^{3}$ is of order $2$, 
which maps $E_+$ to $E_-$ and belongs to ${\rm Ker}(\rho_{\FEC})$,
which gives the generator of the subgroup $\ZZ/2\ZZ$ of ${\rm
Out}_{\CC\text{-}{\bf alg}}(\CC Q/I)$.
\qed\end{pf}

We are interested in classifying exceptional objects in $\D$. 
\begin{prop}
The following holds for exceptional objects and complete exceptional sequences.
\begin{enumerate}
\item 
If $E\in \D$ is an exceptional object and $E \ncong E_\pm [k]$ for any $k \in \ZZ$, 
then there exists a full exceptional collection $(E_1,E_2,E_3)$ with $E_1\cong E$.
\item
The right $B_3$-action on $\FEC(\D)/\ZZ^3$ is transitive.
\end{enumerate}
\end{prop}
\begin{pf}
A geometric description based on categorical equivalence $\D\cong \P\W(\SS_1^1, \MM, \ell)$ is useful for the classification of exceptional cases.
\begin{lem}[{\cite[Theorem~4.3]{HKK} and \cite[Lemma~2.1~(1)]{CS}}]
Let $E$ be an exceptional object in $\D$. 
The corresponding object $L$ in $\P\W(\SS_1^1, \MM, \ell)$ under the equivalence $\D\cong \P\W(\SS_1^1, \MM, \ell)$ 
is an arc connecting two connected components of $\p\SS_1^1\setminus \MM$ without self-intersections. 
\qed\end{lem}
\begin{lem}[{\cite[Proposition~3.8]{O2}}]
Let $E$ be an exceptional object in $\D$ and let $L$ be the corresponding arc. 
$L$ is homotopic to a part of the  boundary $\p\SS^1_1$ if and only if $E$ satisfies $(\S_\D[-1])^2(E)\cong E[2]$. 
Otherwise, $L$ is non-boundary, and cutting along $L$ yields a sphere with two boundary components, each containing one marked point.
\qed\end{lem}
Based on these two lemmas, the first statement is given by \cite[Proposition~4.11]{CS}.

The second statement is shown by Chang--Haiden--Schroll~\cite[Corollary~5.3]{CHS}.
What is crucial is that the braid group $B_3=\rho_{\FEC}(B(\D))$ is identified with the mapping class group of $\SS^1_1$, 
the orientation preserving diffeomorphisms $\SS^1_1\longrightarrow \SS^1_1$ fixing $\p\SS^1_1$ modulo isotopy.
Also, that the left action of $\rho_{\FEC}(B(\D))$ and the right action of $B_3$ are identified (recall Proposition~\ref{prop: left right B_3}) plays an important role.
\qed\end{pf}

\begin{cor}
An exceptional object $E$ is uniquely determined by $[E]\in K_0(\D)$ up to the actions of $\S_\D^2$ and $[2]$.
\qed\end{cor}

\subsection{Space of Bridgeland stability conditions}

The result of Haiden--Katzarkov--Kontsevich, when applied to the present setting, yields the following
\begin{prop}[{\cite[Theorem~5.3]{HKK}}]
There exists an injective holomorphic map from the universal covering $UCov(M)$ of $M$ to $\Stab(\D)$ 
whose image is open and closed, thus a union of connected components of $\Stab(\D)$.
\qed\end{prop}

By Takeda~\cite[Lemma~62 and~Lemma 67]{Tak} (see also \cite[Remark~4.6]{HW}), 
it turns out that the injective map $UCov(M)\longrightarrow \Stab(\D)$ is indeed an isomorphism of complex manifolds.
The stability function is given by the exponential period integrals associated to $F$ with the primitive form given in Corollary \ref{primitive_form} (see Section~\ref{sec:Gamma} and put $u=1$),
on which the auto-equivalences $\S_\D^2$ and $[2]$ act trivially.

\subsection{Serre dimension and the global dimension function}

\begin{defn}[{\cite[Definition~2.4]{KOT} (see also \cite[Proposition~6.13]{EL})}]
The {\it upper Serre dimension} $\overline{\S{\rm dim}}\ \D$ and the {\it lower Serre dimension} $\underline{\S{\rm dim}}\ \D$ of $\D$ are given by
\begin{equation}
\overline{\S{\rm dim}\ }\D:=\lim_{t\to\infty}{\frac{h_t(\S_\D)}{t}},\quad 
\underline{\S{\rm dim}\ }\D:=\lim_{t\to-\infty}{\frac{h_t(\S_\D)}{t}}.
\end{equation}
where 
\begin{equation}
h_t(\S_\D)=\lim_{n\rightarrow\infty}\frac{1}{n}\log\sum_{m\in\ZZ}\dim_{\CC} {\rm Hom}_\D(G,\S_\D^n(G)[m]) e^{-mt},
\end{equation}
and $G$ is a split generator of $\D$.
\end{defn}
\begin{rem}
The entropy function $h_t(\S_\D)$ does not depend on the choice of a split generator. 
\end{rem}
\begin{defn}[\cite{IQ,Q1}]
For a stability condition $\sigma=(Z,\P)$ on $\D$, the global dimension ${\rm gldim}\ \sigma$ of $\sigma$ is given by
\begin{equation}
{\rm gldim}\ \sigma:=\sup\left\{\phi_2-\phi_1|~\Hom_\D(A_1,A_2)\neq0\text{ for }A_i\in\P(\phi_i)\right\}\in\RR_{\geq0}\cup\{\infty\}. 
\end{equation}
\end{defn}
\begin{prop}
\[
\overline{\S{\rm dim}}\ \D=\inf_{\sigma\in{\rm Stab}(\D)}{\rm gldim}\ \sigma =2.
\]
\end{prop}
\begin{pf}
There exists $\sigma $ satisfying ${\rm gldim}\ \sigma =2$ \cite[Corollary~5.13]{Q2}.
If we take $P(3)\oplus P(2)\oplus P(1)\oplus E_+\oplus E_-$ as $G$,
we easily obtain the inequality $\overline{\S{\rm dim}}\ \D\ge 2$ since $\S_\D(E_\pm)=E_\mp[2]$.
\qed\end{pf}
\begin{cor}
There is no $\S_\D$-invariant stability condition on $\D$.
\end{cor}
\begin{pf}
Taking $P(3)\oplus P(2)\oplus P(1)\oplus S_{1,c_1}$ as $G$, it is easy to show that $\underline{\S{\rm dim}\ }\D\le 1$ since $\S_\D(S_{1,c_1})\cong S_{1,c_1}[1]$.
Therefore, we have $\underline{\S{\rm dim}\ }\D\le 1< \overline{\S{\rm dim}\ }\D=2$.
The statement follows from \cite[Proposition~6.17]{KP}.
\qed\end{pf}
\begin{rem}
Consider the set ${\rm SOD}_{\S}(\D)$ of all semiorthogonal decompositions for $\D$ 
consisting of admissible subcategories equipped with Serre-invariant stability conditions.
We expect $\Stab(\D)$ to have a stratification indexed by a set ${\rm SOD}_{\S}(\D)$ 
(this will be discussed in a separate paper).
Note that there is a natural inclusion $\FEC(\D)/\ZZ^3\subset {\rm SOD}_\S(\D)$.
Since the derived category of the Kronecker quiver, the semi-orthogonal complement of an exceptional object $E$ with $E \ncong E_\pm [k]$ for any $k \in \ZZ$, 
does not admit any Serre-invariant stability condition and neither does $E_\pm^\perp$ by \cite[Corollary~3.4]{Su}, it turns out that $\FEC(\D)/\ZZ^3={\rm SOD}_\S(\D)$. 
Therefore, it is naturally explained in a categorical way that the universal covering of the configuration space ${\rm Conf}(\CC,3)$, which is $\Stab(\D)$, admits a stratification indexed by $\FEC(\D)/\ZZ^3$, 
called the ``Stokes region'' (see Hertling--Roucairol~\cite[Definition~3.10]{HR}), together with the discussion on the Lyashko--Looijenga map explained later. 
\end{rem}

For a real number $0\le a\le 2$, set
\[
\chi_a:=\begin{pmatrix}
1 & a & a\\
0 & 1 & a\\
0 & 0 & 1
\end{pmatrix}.
\]
With respect to the $\ZZ$-basis $\{[P(3)],[P(2)],[P(1)]\}$ of $K_0(\D)$,  we have $[\S_\D]=\displaystyle\lim_{a\to 2}\chi_a^{-1}\chi_a^T$.
Since the characteristic polynomial of $\chi_a^{-1}\chi_a^T$ is
\[
(t-1)(t^2 -(a^3-3a^2+2) t + 1),
\]
the eigenvalues are 
\[
t=1, \ \frac{ (a^2-2a-2)(a-1) \pm a(a-2)\sqrt{(a+1)(a-3)}}{2}.
\]
Examining the behavior of the phases of the eigenvalues as the parameter $a$
varies from $0$ to $2$ shows that they change monotonically from $0$ to $0, \pm 1/2$.
Therefore, based on the idea of Cecotti--Vafa \cite{CV} (see also Balnojan--Hertling \cite{BH}), 
it is natural to consider the ``exponents" for $\D$ to be $0$, $1/2$ and $1$  (shifted so that the smallest one is 0), 
and hence the ``dimension" of $\D$ to be $1$.
The "exponents" here are rational numbers, analogues of those $q$ with $H^{p,q}(V)\ne 0$ for a smooth proper Deligne--Mumford stack $V$. 

These "exponents" and "dimension" are  consistent with those of the Frobenius structure discussed later, 
but the "dimension" differs from $2$, the upper Serre dimension $\overline{\S{\rm dim}}\ \D$ 
(and hence the infimum $\displaystyle\inf_{\sigma\in{\rm Stab}(\D)} {\rm gldim}\ \sigma$ of the global dimension function) computed above, 
although they coincide in nice cases.
On the other hand, if we regard $\overline{\S{\rm dim}}\ \D$ as the ``dimension",
it is consistent with the variance formula  (see Hertling \cite{He1}, Ebeling-T \cite{ET} and references therein) for the ``exponents"
under the assumption that the ``first Chern class of $\D$" is numerically zero as in the discussion of the gamma integral structure later:
\[
\left(-\frac{1}{2}\right)^2+0^2+\left(\frac{1}{2}\right)^2=\frac{1}{12}\cdot 3\cdot 2,
\]
where $3$ in the RHS is the dimension of the total Hochschild homology group $H\!H_\bullet (\D)$ of $\D$ which is isomorphic to $K_0(\D)\otimes_\ZZ \CC$.
Understanding this discrepancy remains a challenge for future work.

%%%%%%%%%%%%%%%%%%%%%%%%%%%%%%%%%%%%%%%%%%%%%%%%%%%%%%%%%%%%%%%%
\section{Frobenius structure}\label{sec:Frob_Str}
Recall that
$M=\CC\times \CC^\ast\times\HH$ with the coordinate system $(s_1,s_2,\tau)$, 
$\CC\times M$ with the coordinate system $(z;s_1,s_2,\tau)$,  
$\X$ is the quotient of $\CC\times M$ under the group action 
$(z;s_1,s_2,\tau)\mapsto (z+m+n\tau;s_1,s_2,\tau)$, $m,n\in\ZZ$, and 
$p$ the natural projection $\X\longrightarrow M$;
\[
p: \X\longrightarrow M,\quad [(z;s_1,s_2,\tau)]\mapsto (s_1,s_2,\tau).
\]
Here $\X^o$ is the space obtained by removing from $\X$ the zero section of $p$, namely,
\[
\X^o:=\X\setminus \{[(0;s_1,s_2,\tau)]\in \X~|~(s_1,s_2,\tau)\in M\},
\]
and $F: \X^o \longrightarrow \CC$ is a holomorphic map given by
\[
F=F(z;{\bf s})=F(z;s_1,s_2,\tau):=s_2^2\widetilde\wp(z;\tau)+s_1.
\]
It is important that  the critical set $\C$ of $F$ relative to $M$ is given by 
\[
\C:=\left\{(z;{\bf s})\in\X^o\,\left|\, \frac{\p F}{\p z}=0\right.\right\}=\left\{(z;{\bf s})\in\X^o\,\left|\, z=\frac{1}{2},\ \frac{1+\tau}{2},\ \frac{\tau}{2}\right.\right\},
\]
and hence the map $p|_{\C}:\C\longrightarrow M$ is an unramified covering map of degree three.
Therefore, the $\O_M$-algebra 
\[
p_*\O_\C:=p_*\O_{\X^o}\left/\left(\frac{\p F}{\p z}\right)\right.
\]
is semi-simple at {\it any} point ${\bf s}\in M$. See Proposition~\ref{prop: ss} below for a more precise statement.
\subsection{Kodaira--Spencer map}
\begin{prop}
The $\O_M$-linear map
\[
{\rm KS}: \T_M\longrightarrow p_*\O_\C, \quad \frac{\p}{\p s_i}\mapsto \left[\frac{\p F}{\p s_i}\right],
\]
is an isomorphism.
\end{prop}
\begin{pf}
Since $p_*\O_{\X^o}$ is an $\O_M$-algebra generated by $\widetilde\wp(z;{\bf s})$ and $\p\widetilde\wp(z;{\bf s})/\p z$,
its quotient by the ideal $(\p F/\p z)=(\p\widetilde\wp(z;{\bf s})/\p z)$ is a free $\O_M$-module of rank $3$ generated by 
$[1]$, $[\widetilde\wp(z;{\bf s})]$ and $[\widetilde\wp(z;{\bf s})^2]$. 
The statement follows from the equation \eqref{eq: key identity-1}.
\qed
\end{pf}

\begin{defn}
Define an $\O_M$-bilinear product $\circ:\T_M\times\T_M\longrightarrow \T_M$ as 
\[
\frac{\p}{\p s_i}\circ \frac{\p}{\p s_j}:={\rm KS}^{-1}\left(\left[\frac{\p F}{\p s_i}\cdot\frac{\p F}{\p s_j}\right]\right).
\]
\end{defn}

\begin{defn}
The holomorphic vector fields $e, E\in\Gamma(M,\T_M)$ defined as 
\[
e:={\rm KS}^{-1}([1]),\quad E:={\rm KS}^{-1}([F]),
\]
are called the {\it unit vector field} and the {\it Euler vector field}, respectively.
\end{defn}
It is obvious but important that 
\[
\frac{\p F(z;{\bf s})}{\p s_1}=1,\quad F(z;{\bf s})=s_1\frac{\p F(z;{\bf s})}{\p s_1}+\frac{1}{2}s_2\frac{\p F(z;{\bf s})}{\p s_2},
\]
and hence 
\[
e:=\frac{\p}{\p s_1},\quad E:=s_1\frac{\p}{\p s_1}+\frac{1}{2}s_2\frac{\p}{\p s_2}.
\]
\begin{prop}\label{prop: ss}
The holomorphic functions $u_1({\bf s}), u_2({\bf s}), u_3({\bf s})\in\Gamma(M,\O_M)$ defined as 
\[
u_i({\bf s})=s_1+s_2^2\widetilde{e}_i(\tau),\quad
 (\text{recall \eqref{eq: e} for the definition of $\widetilde{e}_i(\tau)$ })
\]
form a coordinate system on $M$. Moreover, we have
\[
e=\frac{\p}{\p u_1}+\frac{\p}{\p u_2}+\frac{\p}{\p u_3},\quad E=u_1\frac{\p}{\p u_1}+u_2\frac{\p}{\p u_2}+u_3\frac{\p}{\p u_3},\quad 
\frac{\p}{\p u_i}\circ\frac{\p }{\p u_j}=\delta_{ij}\frac{\p}{\p u_j},\ i,j=1,2,3,
\]
where $\delta_{ij}$ is the Kronecker delta.
\end{prop}
\begin{pf}
We have $\C=\C_{1}\sqcup  \C_{2} \sqcup  \C_{3}$ where
\begin{equation*}
\displaystyle
\C_{1}:=\left\{\left(\frac{1}{2}; {\bf s}\right)\in \X^{o}\right\}, \quad \C_{2}:=\left\{\left(\frac{1+\tau}{2}; {\bf s}\right)\in \X^{o}\right\}, \quad \C_{3}:=\left\{\left(\frac{\tau}{2}; {\bf s}\right)\in \X^{o}\right\}.
\end{equation*}
Therefore, we have $F|_{\C_{i}}=u_{i}({\bf s})$ and this yields the statement.
\qed
\end{pf}
Note that $u_i\ne u_j$ if $i\ne j$ since $\widetilde{e}_i(\tau)\ne \widetilde{e}_j(\tau)$ if $i\ne j$.
\subsection{Pairing $\eta$ and its flat coordinates}
Define a symmetric $\O_M$-bilinear form $\eta:\T_M\times \T_M\longrightarrow \O_M$ by
\[
\eta\left(\frac{\p}{\p s_i}, \frac{\p}{\p s_j}\right):=\frac{1}{2\pi\sqrt{-1}}\oint_{\frac{\p F(z;{\bf s})}{\p z}=0} \frac{\frac{\p F(z;{\bf s})}{\p s_i}\frac{\p F(z;{\bf s})}{\p s_j}}{\frac{1}{2\pi\sqrt{-1}}\frac{\p F(z;{\bf s})}{\p z}}(2\pi\sqrt{-1}dz)
=-2\pi\sqrt{-1}\oint_{z=0} \frac{\frac{\p F(z;{\bf s})}{\p s_i}\frac{\p F(z;{\bf s})}{\p s_j}}{\frac{\p F(z;{\bf s})}{\p z}}dz.
\]
It is quite important that $\eta$ is {\it not} an invariant of $F$ but depends on 
the particular choice of a holomorphic volume form $2\pi\sqrt{-1}dz$.

By direct calculation using the facts that  $\widetilde\wp(z;\tau)=1/(2\pi\sqrt{-1}z)^2+(\text{holomorphic})$ near $z=0$, 
\eqref{eq: cubic} and \eqref{eq: key identity}, we have
\[
\left(\eta\left(\frac{\p}{\p s_i}, \frac{\p}{\p s_j}\right)\right)=
\begin{pmatrix}
0 & 0 & 1\\
0 & 2 & \frac{1}{6}s_2 E_2(\tau)\\
1 & \frac{1}{6}s_2 E_2(\tau) & \frac{1}{6}s_2^2\frac{1}{2\pi\sqrt{-1}}\frac{\p E_2(\tau)}{\p\tau}
\end{pmatrix}.
\]

Therefore, considering the following coordinate transformation
\begin{equation}
t_1:=s_1+\frac{1}{12}s_2^2E_2(\tau),\quad t_2:=s_2,\quad t_3:=s_3=2\pi\sqrt{-1}\tau,
\end{equation}
we obtain the constant matrix
\[
\left(\eta_{ij}\right):=\left(\eta\left(\frac{\p}{\p t_i}, \frac{\p}{\p t_j}\right)\right)=
\begin{pmatrix}
0 & 0 & 1\\
0 & 2 & 0\\
1 & 0 & 0
\end{pmatrix}.
\]
In terms of this coordinate system $(t_1,t_2,\tau)$ of $M$, $F$ is rewritten as 
\[
F(z;{\bf t})=t_2^2\left(\widetilde\wp(z;\tau)-\frac{1}{12}E_2(\tau)\right)+t_1.
\]
We have
\begin{equation}\label{eq: exponent 0}
\frac{\p F(z;{\bf t})}{\p t_1}=1,\quad F(z;{\bf t})=t_1\frac{\p F(z;{\bf t})}{\p t_1}+\frac{1}{2}t_2\frac{\p F(z;{\bf t})}{\p t_2},
\end{equation}
and the unit vector field $e$ and the Euler vector field $E$ are given by 
\[
e:=\frac{\p}{\p t_1},\quad E:=t_1\frac{\p}{\p t_1}+\frac{1}{2}t_2\frac{\p}{\p t_2}.
\]

In order to define the notion of a primitive form, it is necessary to define a Saito structure associated to $F$, 
namely, a tuple consisting of the filtered de Rham cohomology group $\H_F$ with the increasing filtration $\H_F^{(k)}$, $k\in\ZZ$,
the Gau\ss--Manin connection $\nabla$ on $\H_F$ and the higher residue pairing $K_F$ on $\H_F$. 
In this paper, we omit the details about those objects and refer the interested reader to \cite{SaTa}.
In the setting of this paper, Milanov constructs the triple $(\H_F, \nabla, K_F)$ and shows that 
any ``prime form" in the sense of Dubrovin \cite[Lecture~5]{D} is a primitive form \cite{M}.
 
As in \cite[Lemma~5.2]{IST}, the equation \eqref{eq: exponent 0} and the following proposition imply that $2\pi\sqrt{-1}dz$ gives a primitive form in an elementary way.

\begin{prop}\label{prop: primitive form}
We have
{\Small
\begin{eqnarray*}
\frac{\p F(z;{\bf t})}{\p t_2}\cdot\frac{\p F(z;{\bf t})}{\p t_2}&=&2\frac{\p F(z;{\bf t})}{\p t_3}
-\frac{t_2}{2}E_2(\tau)\frac{\p F(z;{\bf t})}{\p t_2}-\frac{t_2^2}{2}\frac{1}{2\pi\sqrt{-1}}\frac{\p E_2(\tau)}{\p\tau}\frac{\p F(z;{\bf t})}{\p t_1}+\phi_{2,2}(z;{\bf t})\frac{\p F(z;{\bf t})}{\p z},\\
\frac{\p F(z;{\bf t})}{\p t_2}\cdot\frac{\p F(z;{\bf t})}{\p t_3}&=&-\frac{t_2^2}{4}\frac{1}{2\pi\sqrt{-1}}\frac{\p E_2(\tau)}{\p\tau}\frac{\p F(z;{\bf t})}{\p t_2}
-\frac{t_2^3}{6}\frac{1}{(2\pi\sqrt{-1})^2}\frac{\p^2 E_2(\tau)}{\p \tau^2}\frac{\p F(z;{\bf t})}{\p t_1}+\phi_{2,3}(z;{\bf t})\frac{\p F(z;{\bf t})}{\p z},\\
\frac{\p F(z;{\bf t})}{\p t_3}\cdot\frac{\p F(z;{\bf t})}{\p t_3}&=&-\frac{t_2^3}{12}\frac{1}{(2\pi\sqrt{-1})^2}\frac{\p^2 E_2(\tau)}{\p\tau^2} \frac{\p F(z;{\bf t})}{\p t_2}
-\frac{t_2^4}{24}\frac{1}{(2\pi\sqrt{-1})^3}\frac{\p^3 E_2(\tau)}{\p \tau^3}\frac{\p F(z;{\bf t})}{\p t_1}+\phi_{3,3}(z;{\bf t})\frac{\p F(z;{\bf t})}{\p z},
\end{eqnarray*}}
where
{\Small
\begin{gather*}
\phi_{2,2}(z;{\bf t})=2\left(-\widetilde\zeta(z;\tau)-\frac{1}{12}E_2(\tau)z\right),\ 
\phi_{2,3}(z;{\bf t})=2t_2\frac{1}{2\pi\sqrt{-1}}\frac{\p}{\p \tau}\left(-\widetilde\zeta(z;\tau)-\frac{1}{12}E_2(\tau)z\right),\\
\phi_{3,3}(z;{\bf t})=t_2^2\frac{1}{(2\pi\sqrt{-1})^2}\frac{\p^2}{\p \tau^2}\left(-\widetilde\zeta(z;\tau)-\frac{1}{12}E_2(\tau)z\right),
\end{gather*}}
which satisfy
{\Small
\[
\frac{\p \phi_{2,2}(z;{\bf t})}{\p z}=\frac{\p^2 F(z;{\bf t})}{\p t_2 \p t_2},\quad \frac{\p^2 \phi_{2,3}(z;{\bf t})}{\p z}=\frac{\p^2 F(z;{\bf t})}{\p t_2 \p t_3},
\quad \frac{\p \phi_{3,3}(z;{\bf t})}{\p z}=\frac{\p^2 F(z;{\bf t})}{\p t_3 \p t_3}.
\]}
\end{prop}
\begin{pf}
It can be shown by direct computation using equations \eqref{eq: cubic-tilde a}, \eqref{eq: cubic-tilde b}, \eqref{eq: key identity-tilde-1}, \eqref{eq: key identity-tilde}
and \eqref{lem: E_2}. 
Although it is not difficult, it is technical, so we will summarize the computation in Section~\ref{subsec: proof-prop}. 
\qed\end{pf}

\begin{cor}\label{primitive_form}
The class $\zeta_F:=[2\pi\sqrt{-1}dz]$ of the differential form $2\pi\sqrt{-1}dz$  
in the $0$-th filtered de Rham cohomology group $K_F^{(0)}$ is a primitive form with the minimal exponent $0$.
In particular, the tuple $(\circ, \eta, e, E)$ defines a (globally) semi-simple Frobenius structure on $M$ of rank $3$ and dimension $1$,
whose Frobenius potential $\F$ is
\[
\F(t_1,t_2,2\pi\sqrt{-1}\tau):=\frac{1}{2}t_1^2(2\pi\sqrt{-1}\tau)+t_1t_2^2-\frac{1}{24}t_2^4E_2(\tau).
\]
\qed
\end{cor}

The intersection form $g$ of the Frobenius structure is a non-degenerate symmetric $\O_M$-bilinear form 
$g:\Omega^1_M\times \Omega^1_M\longrightarrow \O_M$, which is defined, in terms of the Frobenius potential and flat coordinates, by
\begin{equation*}
g(dt_{i}, dt_{j}):=\sum_{a,b=1}^{3}\eta^{ia}\eta^{jb}E\left(\frac{\p^2 \F}{\p t_{a}\p t_{b}}\right),\quad \left(\eta^{ij}\right):=\left(\eta_{ij}\right)^{-1}.
\end{equation*}
It is calculated as
\[
(g(dt_i , dt_j))=
\begin{pmatrix}
-\frac{t_2^4}{12}\frac{1}{(2\pi\sqrt{-1})^2}\frac{\p^2 E_2(\tau)}{\p \tau^2} & -\frac{t_2^3}{8}\frac{1}{2\pi\sqrt{-1}}\frac{\p E_2(\tau)}{\p \tau}   & t_1\\
-\frac{t_2^3}{8}\frac{1}{2\pi\sqrt{-1}}\frac{\p E_2(\tau)}{\p \tau}   & \frac{1}{2}t_1-\frac{1}{8}t_2^2E_2(\tau) & \frac{1}{2}t_2\\
t_1 & \frac{1}{2}t_2 & 0
\end{pmatrix}.
\]

There is a uniqueness of the Frobenius structure, which may be known for experts.
See \cite[Theorem~2.1]{DZ}, \cite[Proposition~5.2]{Sat1} and \cite[Theorem 7.1]{ShTa} for the relevant statements.
\begin{prop}\label{unique}
A Frobenius structure on $M$ of rank three and of dimension one with the following $e$ and $E$
is uniquely determined by the intersection form $g$:
\[
e=\frac{\p}{\p t_1},\quad E=t_1\frac{\p}{\p t_1}+\frac{1}{2}t_2\frac{\p}{\p t_2}.
\]
\end{prop}
\begin{pf}
We use the following relation between the product $\circ$ and the intersection form $g$.
\begin{lem}
Denote by $\Gamma^{ij}_{k}$ the contravariant component of the Levi--Civita connection $\widetilde \nabla$ 
for the intersection form $g$ on $\Omega^1_M$. Namely, 
\[
\Gamma^{ij}_{k}:=g\left(d t_i,\widetilde\nabla_{\frac{\p}{\p t_k}} dt_j\right).
\]
Then, we have
\begin{equation}
\Gamma^{ij}_{k}=d_{j}C^{ij}_{k},\quad C^{ij}_{k}:=\sum_{a,b=1}^{3}\eta^{ia}\eta^{jb}\frac{\p^3 \F}{\p t_{k}\p t_{a}\p t_{b}},
\end{equation}
where $d_{j}$ is a rational number defined by $E(t_j)=d_j t_j$, namely, $d_1=1$, $d_2=1/2$ and $d_3=0$.
\end{lem}
\begin{pf}
See \cite[Lemma~3.4]{D} and apply $d=1$.
\qed
\end{pf}
It is obvious that $C^{ij}_{k}$ can be reconstructed from the intersection form if $j\neq 3$.
However, 
\begin{equation}
C^{i3}_{k}=\sum_{a,b=1}^3 \eta^{ia}\eta^{3 b}\frac{\p^3 \F}{\p t_{k}\p t_{a}\p t_{b}}=\delta_{ik},
\quad (\text{Kronecker's delta}).
\end{equation}
Therefore, $\p^3\F/\p t_i\p t_j\p t_k$ for all $i,j,k=1,\dots, 3$ can be uniquely reconstructed from the intersection form $g$.
\qed
\end{pf}

The discriminant locus $D_F$ is the subset of $M$ defined as 
\[
D_F:=\left\{{\bf s}\in M\,\left|\,F(z;{\bf s})=\frac{\p F(z;{\bf s})}{\p z}=0\text{ for some }z\right.\right\},
\]
which is explicitly given in terms of the Euler vector field as the zero locus of the holomorphic function ${\rm det}(E\circ -)$: 
{\small
\begin{eqnarray*}
& &{\rm det}(E\circ -)\\
&=&t_1^3-\frac{t_1^2t_2^2}{4}E_2(\tau)+\frac{t_1t_2^4}{4}\frac{1}{2\pi\sqrt{-1}}\frac{\p E_2(\tau)}{\p \tau}-\frac{t_2^6}{24}\frac{1}{(2\pi\sqrt{-1})^2}\frac{\p^2 E_2(\tau)}{\p \tau^2} \\
&=&\left(t_1+t_2^2\left(\widetilde{e}_1(\tau)-\frac{1}{12}E_2(\tau)\right)\right)\left(t_1+t_2^2\left(\widetilde{e}_2(\tau)-\frac{1}{12}E_2(\tau)\right)\right)\left(t_1+t_2^2\left(\widetilde{e}_3(\tau)-\frac{1}{12}E_2(\tau)\right)\right)\\
&=&u_1u_2u_3.
\end{eqnarray*}}
Up to a non-zero constant multiple, it is equal to ${\rm det}(g(dt_i , dt_j))$.

Let ${\rm Conf}(\CC,3)$ be the configuration space of three distinct points on $\CC$ and denote by $[(u_1,u_2,u_3)]$ the point of ${\rm Conf}(\CC,3)$ consisting of $\{u_1,u_2,u_3\}\subset \CC$.
\begin{defn}
The {\it Lyashko--Looijenga map} ${\rm LL} : M \longrightarrow {\rm Conf}(\CC,3)$ is a holomorphic map defined as
\begin{equation*}
{\rm LL} ({\bf s}):= [(u_1, u_2,u_3)],
\end{equation*}
where $u_i$'s are canonical coordinates defined in Proposition \ref{prop: ss}, namely,
\[
u_i=s_1+s_2^2\widetilde{e}_i(\tau)=t_1+t_2^2\left(\widetilde{e}_i(\tau)-\frac{1}{12}E_2(\tau)\right),\quad i=1,2,3.
\]
\end{defn}
For the general theory of the Lyashko--Looijenga map for Frobenius manifolds (or F-manifolds), we refer the reader to \cite[Section~3.5]{He2}.

\begin{prop}\label{deg_LL_map}
The Lyashko--Looijenga map ${\rm LL}$ induces the isomorphism
\[
{\rm LL}^{alg}: M/{\rm SL}(2;\ZZ)\longrightarrow {\rm Conf}(\CC,3).
\]
In particular, we have 
\[
{\rm deg}\ {\rm LL}^{alg}=1=\#(B_3\backslash(\FEC(\D)/\ZZ^3))=\#((\FEC(\D)/\ZZ^3)/B_{3})
\]
where $B_3\backslash(\FEC(\D)/\ZZ^3)$ means the quotient with respect to the left $B_3=\rho_{\rm FEC}(B(\D))$-action
and  $(\FEC(\D)/\ZZ^3)/B_{3}$ does the quotient with respect to the right $B_3$-action.

\end{prop}
\begin{pf}
By Proposition~\ref{prop: left right B_3}, we have $\#((\FEC(\D)/\ZZ^3)/B_{3})=\#(B_{3}\backslash(\FEC(\D)/\ZZ^3))$.
The isomorphism ${\rm LL}^{alg}$ is nothing but the classical theory of elliptic functions. 
In particular, our case is the Jacobi inversion problem of elliptic curves with a quadratic term in the cubic polynomial 
which gives rise to the translation by $s_{1}$. 
Hence ${\rm LL}^{alg}$ is an isomorphism between $3$-dimensional complex manifolds.
See \cite{Zac} and references therein for example.
\qed
\end{pf}

Let $M^{reg}:=M \backslash D_F$ and $F_{\bf s}^{-1}(0):= \{(z; {\bf s}) \in \X^o~|~F(z; {\bf s})=0\}$ for a point ${\bf s}\in M$.
Fix a point ${\bf s}=(s_1,s_2,2\pi\sqrt{-1}\tau)\in M^{reg}$, then the relative homology long exact sequence yields the following short exact sequence:
\begin{equation}
0\longrightarrow H_1 (p^{-1}({\bf s});\ZZ)\longrightarrow H_1 (p^{-1}({\bf s}), F_{\bf s}^{-1}(0);\ZZ)
\stackrel{\p}{\longrightarrow} \widetilde{H}_0 (F_{\bf s}^{-1}(0);\ZZ)\longrightarrow 0,
\end{equation}
where $\widetilde{H}_0$ denotes the $0$-th reduced homology group.
Since $\p F(z;{\bf s})/\p z\ne 0$ for ${\bf s}\in M^{reg}$, $F^{-1}_{\bf s}(0)=\{(x;{\bf s}),(-x;{\bf s})\}$ for some 
$x\ne 1/2, 1/2+\tau/2,\tau/2$ $({\rm mod}\ \ZZ+\ZZ\tau)$. 

For ${\bf s}={\bf s}_0$, by moving along the real axis from $0$ to $+\infty$, 
we construct a homeomorphism between $\overline{F}_{{\bf s}_0}^{-1}(0)$ and $\overline{F}_{{\bf s}_0}^{-1}(1\cdot\infty)$. 
Thus, we obtain an isomorphism 
\[
H_1 (p^{-1}({\bf s}_0), F_{{\bf s}_0}^{-1}(0);\ZZ)= H_1 (\overline{p}^{-1}({\bf s}_0), \overline{F}_{{\bf s}_0}^{-1}(0);\ZZ)\cong
H_1 (\overline{p}^{-1}({\bf s}_0), \overline{F}_{{\bf s}_0}^{-1}(1\cdot\infty);\ZZ)\cong H_1(\SS_1^1,\p\SS_1^1\setminus\MM;\ZZ).
\]
In particular, the class of the path in $\CC$ from $-x$ to $x$ that encircles the origin $0$ in a clockwise direction
is identified with $\alpha$.
Together with the isomorphism $H_1 (p^{-1}({\bf s}_0);\ZZ)\cong H_1(\SS_1^1;\ZZ)$, we have
\[
H_1 (p^{-1}({\bf s}_0), F_{{\bf s}_0}^{-1}(0);\ZZ)=\ZZ\alpha\oplus \ZZ\delta_1\oplus \ZZ\delta_2,\quad  
\widetilde{H}_0 (F_{{\bf s}_0}^{-1}(0);\ZZ)= \ZZ([x]-[-x])=\ZZ\p\alpha,
\]
where $\delta_1$ and $\delta_2$ are the homology classes of the longitude and the meridian of $\SS_1^1$.
Here, it is also useful to take another $\ZZ$-basis $\{\alpha_1,\alpha_2,\alpha_3\}$ of $H_1 (p^{-1}({\bf s}_0), F_{{\bf s}_0}^{-1}(0);\ZZ)$ where 
\[
\alpha_1:=\alpha-\delta_1, \quad \alpha_2:=-\alpha+\delta_1+\delta_2, \quad \alpha_3:=[S(3)]=\alpha-\delta_2. 
\]
Under the identification $H_1 (p^{-1}({\bf s}_0), F_{{\bf s}_0}^{-1}(0);\ZZ)\cong H_1(\SS_1^1,\p\SS_1^1\setminus\MM;\ZZ)$, 
this corresponds to the cycles represented by the dashed lines in Figure~\ref{fig: 2} below.
\begin{figure}[htbp]
  \centering
	\begin{tikzpicture}[scale=5,>=Stealth,
		myred/.style={line width=1.2pt,red!70!black},
		myblue/.style={line width=1.2pt,blue!70!black},
		mygreen/.style={line width=1.2pt,green!60!black},
		pole/.style={black, font=\scriptsize},
		tick/.style={fill=black, circle, inner sep=0.8pt}
	]
	% ===== 枠とラベル =====
	%\draw[black, line width=0.7pt] (0,0) rectangle (1,1);
	\draw (0.13,0) -- (0.87,0);
	\draw (0,0.13) -- (0,0.87);
	\draw (0.13,1) -- (0.87,1);
	\draw (1,0.13) -- (1,0.87);
	\coordinate (A) at (0.02,1);
	\coordinate (B) at (0,0);
	\coordinate (C) at (0.98,0);
	\coordinate (D) at (1,1);

	% 周期（極）
	\node[pole, anchor=north east] at (0,0) {};
	\node[pole, anchor=north west] at (1,0) {};
	\node[pole, anchor=south east] at (0,1) {};
	\node[pole, anchor=south west] at (1,1) {};

	% 半周期（臨界点）
	\coordinate (E) at (0.0,0.5);   % 1/2
	\coordinate (F) at (0.5,1);   % i/2
	\coordinate (G) at (1,0.5);   % 1/2
	\coordinate (H) at (0.5,0);   % i/2
	\coordinate (I) at (0.5,0.5);   % (1+i)/2

	\node[tick] at (E) {};
	\node[tick] at (F) {};
	\node[tick] at (G) {};
	\node[tick] at (H) {};
	\node[tick] at (I) {}; 

	% ===== 逆像 =====

	% L_3
	\draw[red, thick]
		(E) .. controls +(0.04,0) and +(0,-0.0) .. (A);
	\draw[red, thick]
		(G) .. controls +(-0.04,0) and +(0,0) .. (C);

	% L_1
	\draw[blue, thick] 
		(F) .. controls +(-0.00,-0.4) and +(0,-0.5) .. (A);
	\draw[blue, thick] 
		(H) .. controls +(0.00,0.4) and +(0,0.5) .. (C);

	% L_2
	\draw[mygreen, thick] 
		(I) .. controls +(-0.53,0.02) and +(0,0) .. (A);
	\draw[mygreen, thick] 
		(I) .. controls +(0.53,-0.02) and +(0,0) .. (C);

	\node[red, anchor=north west]     at (-0.15,0.8) {$L_3$};
	\node[mygreen, anchor=south]         at (0.35,0.4) {$L_2$};
	\node[blue, anchor=south east] at (0.4,0.7) {$L_1$};

	% 四隅ごとにブローアップ（半径 r）
	\def\r{0.13}

	% 左下
	%\fill[white] (0,0) -- (\r,0) arc[start angle=0,end angle=90,radius=\r] -- cycle;
	%\draw (0,0) -- (\r,0) arc[start angle=0,end angle=90,radius=\r] -- cycle;
	\draw[black, very thick] (0.13,0) arc (0:90:0.13);
	% 右下
	\fill[white] (1,0) -- (1-\r,0) arc[start angle=180,end angle=90,radius=\r] -- cycle;
	%\draw (1,0) -- (1-\r,0) arc[start angle=180,end angle=90,radius=\r] -- cycle;
	\draw[black, very thick] (0.87,0) arc (180:90:0.13);

	% 左上
	\fill[white] (0,1) -- (0,1-\r) arc[start angle=270,end angle=365,radius=\r] -- cycle;
	%\draw (0,1) -- (0,1-\r) arc[start angle=270,end angle=360,radius=\r] -- cycle;
	\draw[black, very thick] (0,0.87) arc (270:360:0.13);

	% 右上
	%\fill[white] (1,1) -- (1,1-\r) arc[start angle=270,end angle=180,radius=\r] -- cycle;
	%\draw (1,1) -- (1,1-\r) arc[start angle=270,end angle=180,radius=\r] -- cycle;
	\draw[black, very thick] (0.87,1) arc (180:270:0.13);

	% 正則点
	\coordinate (J) at (0.25,0.25); % x
	\coordinate (K) at (0.75,0.75); % 1+\tau-x

	\node[below left=2pt] at (J) {$x$};
	\node[above right=0.1pt] at (K) {{\small $1+\sqrt{-1}-x$}};

	\node[tick] at (J) {};
	\node[tick] at (K) {};

	% ===== 双対グラフ =====

	% \alpha_3
	\draw[blue, very thick, dashed]
		(J) .. controls +(0.1,-0.1) and +(-0.1,0.1) .. (H);
	\draw[blue, very thick, dashed]
		(K) .. controls +(-0.1,0.1) and +(0.1,-0.1) .. (F);

	% \alpha_1
	\draw[red, very thick, dashed] 
		(J) .. controls +(-0.1,0.1) and +(0.1,-0.1) .. (E);
	\draw[red, very thick, dashed] 
		(K) .. controls +(0.1,-0.1) and +(-0.1,0.1) .. (G);

	% \alpha_2
	\draw[mygreen, very thick, dashed] 
		(J) .. controls +(0.1,0.1) and +(-0.1,-0.1) .. (I);
	\draw[mygreen, very thick, dashed] 
		(K) .. controls +(-0.1,-0.1) and +(0.1,0.1) .. (I);
	
	\node[red, anchor=south west] at (0.12,0.35) {$\alpha_3$};
	\node[mygreen, anchor=south west] at (0.6,0.5) {$\alpha_2$};
	\node[blue, anchor=south west] at (0.35,0.12) {$\alpha_1$};

	\coordinate (M1) at (0.125,0.035);
	\coordinate (M2) at (0.875,0.965);

	\fill (M1) circle (0.5pt);
	\fill (M2) circle (0.5pt);

	\end{tikzpicture}
  \caption{Cycles representing $\alpha_1$, $\alpha_2$, $\alpha_3$.}
  \label{fig: 2}
\end{figure}
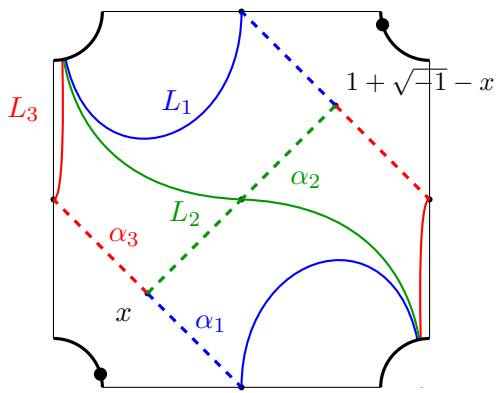

For $\gamma,\gamma'\in H_1 (p^{-1}({\bf s}_0), F_{{\bf s}_0}^{-1}(0);\ZZ)$, 
define a symmetric $\ZZ$-bilinear form $I$ to be the restriction of the intersection form 
$I_{\widetilde{H}_0}$ on $\widetilde{H}_0 (F_{{\bf s}_0}^{-1}(0);\ZZ)$ as 
\[
I(\gamma,\gamma'):=I_{\widetilde{H}_0}(\p\gamma, \p\gamma').
\]

Let $UCov(M^{reg})$ be the universal cover of $M^{reg}$.
Identify $H_1 (p^{-1}({\bf s}_0), F_{{\bf s}_0}^{-1}(0);\ZZ)$ with the space of global horizontal sections of 
$\displaystyle \bigcup_{{\bf s} \in UCov(M^{reg})} H_1 (p^{-1}({\bf s}), F_{\bf s}^{-1}(0);\ZZ)$.

\begin{lem}[{\cite[Section~5 and Section~3.4]{Sa1}}, {\cite[Section 5]{Ta2}} and {\cite[Section 5]{ShTa}}]\label{period mapping}
We have the following:
\begin{enumerate}
\item
For each horizontal section $\gamma$ of the local system 
$\displaystyle \bigcup_{{\bf s} \in M^{reg}} H_1 (p^{-1}({\bf s}), F_{\bf s}^{-1}(0);\ZZ)$, we have
\begin{equation*}
\displaystyle
\frac{1}{2\pi\sqrt{-1}}\int_{\gamma}2\pi\sqrt{-1}dz\in Sol(\widetilde{\nabla}),
\end{equation*}
where $\widetilde{\nabla}$ is the Levi-Civita connection for the intersection form $g$ and 
$\displaystyle Sol(\widetilde{\nabla}):=\{\psi\in \O_{M^{reg}}~|~\widetilde{\nabla} d\psi=0\}$.
Indeed, 
\begin{equation*}
\int_{\alpha}dz= \left[z\right]^x_{-x}=2x, \quad \int_{\delta_{1}}dz=1, \quad 
\quad \int_{\delta_{2}}dz=\tau,
\end{equation*}
which belong to $Sol(\widetilde{\nabla})$.
\item
There exists an $\O_{M^{reg}}$-isomorphism
\begin{equation*}
\Omega^{1}_{M^{reg}} \simeq \O_{M^{reg}}\otimes_{\CC_{M^{reg}}} dSol(\widetilde{\nabla}), 
\end{equation*}
where $dSol(\widetilde{\nabla})$ is the image of $Sol(\widetilde{\nabla})$ in $\Omega^{1}_{M^{reg}}$ under the differential $d$.
\item For $\gamma, \gamma'\in H_1 (p^{-1}({\bf s}_{0}), F_{{\bf s}_{0}}^{-1}(0);\ZZ)$, we have
\begin{equation*}
\displaystyle
g\left(d\left(\int_{\gamma}2\pi\sqrt{-1}dz\right),d\left(\int_{\gamma'}2\pi\sqrt{-1}dz\right)\right)=-I(\gamma, \gamma'),
\end{equation*}
where $\gamma, \gamma'$ in the left-hand side are corresponding global horizontal sections.
\end{enumerate}
In particular, there exists $\phi\in Sol(\widetilde{\nabla})$ such that $\{d\phi, dx,d\tau\}$ forms a $\CC_{M^{reg}}$-basis 
of $dSol(\widetilde{\nabla})$ and $g$ has the following matrix representation with respect to this basis:
\[
g=\frac{1}{(2\pi\sqrt{-1})^2}\begin{pmatrix}
0 & 0 & 1\\
0 & -\frac{1}{2} & 0\\
1 & 0 & 0
\end{pmatrix}.
\]
\end{lem}
\begin{pf}
See {\cite[Section~6.10]{Ta2}}.
\qed
\end{pf}
\begin{rem}
In \cite{ShTa} and \cite{Ta2}, $F$ is a function in three variables, so the coefficient $(2\pi\sqrt{-1})^{-2}$ was necessary 
in front of the integral. 
However, since we are now dealing with a function in one variable as $F$,
as the weight decreases by one, the coefficient becomes $(2\pi\sqrt{-1})^{-1}$.
\end{rem}
By Lemma ~\ref{period mapping}, we have the following map:
\[
UCov(M^{reg})\longrightarrow \CC\times \CC\times \HH,\quad {\bf s}\mapsto (\phi,x,\tau).
\]
and the period mapping:
\begin{equation*}
UCov(M^{reg})\longrightarrow \Hom_{\ZZ}(H_1 (p^{-1}({\bf s}_0), F_{{\bf s}_0}^{-1}(0);\ZZ),\CC).
\end{equation*}
The period domain becomes $\widetilde \EE^o$ introduced later in \eqref{EEo}
since $\alpha$ is not a vanishing cycle and hence the period map cannot be extended to $x\in 2\ZZ+2\ZZ\tau$. 

\subsection{Elliptic Weyl group}\label{Elliptic Weyl group quotient}
In this subsection, we introduce the elliptic Weyl group $\widetilde{W}$ of type $A_{1}^{(1,1)*}$ and its extension $\widehat{W}$ by ${\rm SL}(2; \ZZ)$ and 
describe their actions on suitable domains.

\begin{defn}
A {\it root system} $\R$ of rank $\mu$ is a tuple $(\N,I,\Delta_{re})$ where
$\N$ is a free $\ZZ$-module of rank $\mu$, $I:\N\times \N\longrightarrow \ZZ$ is a symmetric $\ZZ$-bilinear form,
and $\Delta_{re}$ is a subset of $\N$, called the set of real roots,
satisfying the following properties:
\begin{enumerate}
\item
$\N=\ZZ\Delta_{re}$.
\item
For each $\alpha\in\Delta_{re}$, $I(\alpha,\alpha)=2$.
\item
For each $\alpha\in\Delta_{re}$, define a reflection $r_\alpha\in{\rm Aut}_\ZZ(\N,I)$ by
$r_\alpha(\lambda):=\lambda-I(\alpha,\lambda)\alpha$.
Then, $r_\alpha(\Delta_{re})=\Delta_{re}$.
\end{enumerate}
\end{defn}
\begin{defn}
The group
$W(\R):=\langle r_\alpha\,|\, \alpha\in\Delta_{re}\rangle$ is called the {\it Weyl group} of $\R$.
\end{defn}
\begin{defn}
Let $\R=(\N,I,\Delta_{re})$ be a root system of rank $\mu$.
A subset $B=\{\alpha_1,\dots, \alpha_\mu\}$ of $\Delta_{re}$ is a {\it root basis}
if $\Delta_{re}=W_B\cdot B$ where $W_B:=\langle r_{\alpha_1},\dots, r_{\alpha_\mu}\rangle\subset W(\R)$.
\end{defn}
\begin{defn}
Let $\R=(\N,I,\Delta_{re})$ be a root system of rank $\mu$.
An element $c\in W(\R)$ is called a {\it Coxeter element} of $\R$ if
there exists a root basis $B=\{\alpha_1,\dots, \alpha_\mu\}$ such that
$c=r_{\alpha_1}\dots r_{\alpha_\mu}$.
\end{defn}
\begin{defn}[{Saito~\cite[Section~5.3]{Sa3}}]
A pair $(\R,c)$ of a root system $\R$ and a Coxeter element $c$ of $\R$ is called
a {\it generalized root system}.
\end{defn}

Applying \cite[Proposition 2.10]{STW} to our case, we can have the following generalized root system $(\N,I,\Delta_{re},c)$ from $\D$.
Since we relax the condition that any exceptional object can be completed to a full exceptional collection,
the set of real roots can be given by the $B_3 \ltimes (\ZZ [1])^3$-orbit of the full exceptional collection $(S(1), S(2), S(3))$ :
\[
\N:=K_0(\D)=\ZZ\alpha\oplus\ZZ\delta_1\oplus\ZZ\delta_2=\ZZ\alpha_1\oplus\ZZ\alpha_2\oplus\ZZ\alpha_3,
\]
where $\alpha_1:=[S(1)]=\alpha-\delta_1$, $\alpha_2:=[S(2)]=-\alpha+\delta_1+\delta_2$, $\alpha_3:=[S(3)]=\alpha-\delta_2$. 
\[
\left(I(\alpha_i,\alpha_j)\right):=\left(\chi(\alpha_i,\alpha_j)+\chi(\alpha_j,\alpha_i)\right)=\begin{pmatrix}
2 & -2 & 2\\
-2& 2 & -2\\
2 & -2 & 2
\end{pmatrix}.
\]
\[
W:=\langle r_{\alpha_1}, r_{\alpha_2}, r_{\alpha_3}\rangle\subset \Aut_\ZZ(\N,I).
\]
\[
\Delta_{re}:=W\{\alpha_1,\alpha_2,\alpha_3\}=\{\pm \alpha-(p+1)\delta_1-(q+1)\delta_2\,|\,p,q\in\ZZ,\ pq\in 2\ZZ\}.
\]
\[
c:=[\S_\D[-1]]=r_{\alpha_1}r_{\alpha_2}r_{\alpha_3}=r_{\alpha}.
\]
\begin{rem}
As is explained later in \eqref{hyperplane}, reflection hyperplanes corresponding to simple roots $\alpha_{1}, \alpha_{2}, \alpha_{3}$ are $\displaystyle \frac{1}{2}, \frac{1+\tau}{2}, \frac{\tau}{2}$
respectively.  
\end{rem}
\begin{rem}
The matrix $\left(I(\alpha_i,\alpha_j)\right)$ is not a generalized Cartan matrix and the group $W$ is not a Coxeter group.
\end{rem}

The subgroup $H:=\langle \tau_1, \tau_2\rangle$ of $W$ generated by 
\begin{gather}
\tau_1(\alpha,\delta_1,\delta_2):=r_{\alpha_2}r_{\alpha_3}(\alpha,\delta_1,\delta_2)=(\alpha-2\delta_1,\delta_1,\delta_2),\\
\tau_2(\alpha,\delta_1,\delta_2):=r_{\alpha_2}r_{\alpha_1}(\alpha,\delta_1,\delta_2)=(\alpha-2\delta_2,\delta_1,\delta_2).
\end{gather}
is a free abelian group of rank $2$ and we have 
the following short exact sequence:
\[
\{1\}\longrightarrow H\longrightarrow W\longrightarrow \ZZ/2\ZZ\longrightarrow \{1\}.
\]
More precisely, we have
\begin{eqnarray}
W&\cong &\langle r_{\alpha_1},\tau_1,\tau_2\,|\, r_1^2=1,\ r_{\alpha_1}\tau_ir_{\alpha_1}=\tau_i^{-1},\ i=1,2\rangle\\
&\cong & \langle r_1, r_2, r_3\,|\, r_1^2=1,r_2^2=1,\ r_3^2=1,\ (r_1r_2r_3)^2=1\rangle
\end{eqnarray}

Consider the following free abelian group of rank $5$
\[
\widetilde\N:=\N\oplus\ZZ\kappa_1\oplus\ZZ\kappa_2=\ZZ\alpha\oplus\ZZ\delta_1\oplus\ZZ\delta_2\oplus \ZZ\kappa_1\oplus \ZZ\kappa_2,
\]
and the symmetric $\ZZ$-bilinear form $\widetilde I$ defined by 
\begin{gather*}
\widetilde I(\alpha, \alpha):=2,\quad \widetilde I(\alpha, \delta_1):=0,\quad \widetilde I(\alpha, \delta_2):=0,
\quad \widetilde I(\alpha, \kappa_1):=0,\quad \widetilde I(\alpha, \kappa_2):=0,\\
\widetilde I(\delta_1, \delta_1):=0,\quad \widetilde I(\delta_1, \delta_2):=0,
\quad \widetilde I(\delta_1, \kappa_1):=1,\quad \widetilde I(\delta_1, \kappa_2):=0,
\quad \widetilde I(\delta_2, \delta_2):=0,\\
\widetilde I(\delta_2, \kappa_1):=0,\quad \widetilde I(\delta_2, \kappa_2):=-1,\quad \widetilde I(\kappa_1, \kappa_1):=0,
\quad \widetilde I(\kappa_1, \kappa_2):=0,\quad \widetilde I(\kappa_2, \kappa_2):=0.
\end{gather*}
We identify the lattice $(\N,I)$ with the sub-lattice $(\N\oplus 0\oplus 0,\widetilde I|_{\N\oplus 0\oplus 0})$ of $(\widetilde N,\widetilde I)$.
\begin{rem}\label{rem_hyperbolic_ext_I}
When considering the ${\rm SL}(2;\ZZ)$-action later, it is more convenient to consider the matrix representation of 
$\widetilde I$ with respect to the $\ZZ$-basis $\{\kappa_1,\kappa_2,\alpha,\delta_2,\delta_1\}$, which is given by:
\begin{equation}\label{widetilde I}
\widetilde I=\begin{pmatrix}
0 & 0 & 0 & 0 & 1\\
0 & 0 & 0 & -1 & 0\\
0 & 0 & 2 & 0 & 0\\
0 & -1 & 0 & 0 & 0\\
1 & 0 & 0 & 0 & 0
\end{pmatrix}.
\end{equation}
\end{rem}

Define the hyperbolic extension $\widetilde W$ of the Weyl group $W$ by
\[
\widetilde W:=\langle \widetilde r_\beta\,|\,\beta\in \Delta_{re}\rangle\subset \Aut_\ZZ(\widetilde\N,\widetilde I),\quad 
\widetilde r_\beta(\widetilde \lambda):=\widetilde \lambda-\widetilde I(\beta, \widetilde \lambda)\beta,\quad \widetilde \lambda\in \widetilde \N.
\]
Explicitly, we have
\begin{gather*}
\widetilde r_{\alpha_1}(\alpha, \delta_1,\delta_2,\kappa_1,\kappa_2)=(-\alpha+2\delta_1, \delta_1,\delta_2,\kappa_1+\alpha-\delta_1,\kappa_2),\\
\widetilde r_{\alpha_2}(\alpha, \delta_1,\delta_2,\kappa_1,\kappa_2)=(-\alpha+2\delta_1+2\delta_2, \delta_1,\delta_2,\kappa_1+\alpha-\delta_1-\delta_2,\kappa_2-\alpha+\delta_1+\delta_2)\\
\widetilde r_{\alpha_3}(\alpha, \delta_1,\delta_2,\kappa_1,\kappa_2)=(-\alpha+2\delta_2, \delta_1,\delta_2,\kappa_1,\kappa_2-\alpha+\delta_2),
\end{gather*}
It is important that $\widetilde r_{\alpha}\notin \widetilde W$.

We have
\[
\widetilde c(\alpha, \delta_1,\delta_2,\kappa_1,\kappa_2):=\widetilde r_{\alpha_1}\widetilde r_{\alpha_2}\widetilde r_{\alpha_3}(\alpha, \delta_1,\delta_2,\kappa_1,\kappa_2)
=(-\alpha,\delta_1,\delta_2,\kappa_1-\delta_2,\kappa_2-\delta_1)
\]
and hence ${\widetilde c}^2$ belongs to the center of $\widetilde W$.

It is interesting to observe the relationship between this group and a certain fundamental group.
\begin{prop}
Consider the complex manifold 
$\overline{M}:=\CC\times \CC\times\HH$ with the coordinate system $(s_1,s_2,\tau)$, 
and the holomorphic function $\Delta$ given by 
\[
\Delta:=\left(s_1+s_2^2\widetilde{e}_1(\tau)\right)\left(s_1+s_2^2\widetilde{e}_2(\tau)\right)\left(s_1+s_2^2\widetilde{e}_3(\tau)\right).
\]
The groups $\pi_1(\overline{M}\setminus\{\Delta=0\},*)$ and $\widetilde W$ have presentations by generators and relations as follows
\begin{gather}
\pi_1(\overline{M}\setminus\{\Delta=0\},*)\cong \langle a_1,a_2,a_3\,|\, a_i(a_1a_2a_3)^2=(a_1a_2a_3)^2a_i,\ i=1,2,3\rangle,\\
\widetilde W\cong \langle \widetilde r_1,\widetilde r_2,\widetilde r_3\,|\, \widetilde r_i^2=1,\ \widetilde r_i(\widetilde r_1\widetilde r_2\widetilde r_3)^2=(\widetilde r_1\widetilde r_2\widetilde r_3)^2\widetilde r_i,\ i=1,2,3\rangle.
\label{ellipticArtin}
\end{gather}
\end{prop}
\begin{pf}
Consider the projection 
\[
\CC\times \CC \times \HH \longrightarrow \CC \times \HH \quad (s_1,s_2,\tau)\mapsto (s_2,\tau).
\]
and intersections of its fibers with $\Delta=0$. Then the bifurcation set is $s_2=0$. 
A fiber with $s_{2}\ne 0$ contains three counterclockwise loops $a_{i}$ ($i=1,2,3$) with a common generic marked point $*$ (one may take ${\bf s}_0$ as $*$), 
where each $a_{i}$ encircles $-s^{2}_{2}e_{i}(\tau)$. 
Along a rotation around $s_2=0$, three points $-s^{2}_{2}e_{i}(\tau)$ rotate twice around $s_{1}=0$
without any permutations. This rotation induces a double full twist $((\sigma_{1}\sigma_{2})^{3})^{2}$ on $a_{i}$,
where $((\sigma_{1}\sigma_{2})^{3})^{2}(a_{i})=(a_{1}a_{2}a_{3})^{2}a_{i}(a_{1}a_{2}a_{3})^{-2}$.
By the Zariski-Van Kampen method, these $a_{i}$ form generators for and the relations $a_{i}=((\sigma_{1}\sigma_{2})^{3})^{2}(a_{i})=(a_{1}a_{2}a_{3})^{2}a_{i}(a_{1}a_{2}a_{3})^{-2}$ are 
exactly the ones for the fundamental group $\pi_1(\overline{M}\setminus\{\Delta=0\},*)$. 
\qed
\end{pf}
\begin{rem}
The complement $\overline{M}\setminus\{\Delta=0\}$ is a $K(\pi, 1)$-space. The group \eqref{ellipticArtin} is called the elliptic Artin group of type $A^{(1,1)*}_{1}$.
\end{rem}

The restriction to the subspace $\N$ of the action of $\widetilde W$ yields the surjective group homomorphism
$\widetilde W\longrightarrow W$. Moreover, we have the following commutative diagram
\[
\begin{tikzcd}
& & \{1\}\arrow[d] & \{1\}\arrow[d] & \\
\{1\} \arrow[r] & K:=\langle {\widetilde c}^2\rangle\arrow[r]\arrow[d, equal]& \widetilde H:=\langle \widetilde \tau_1, \widetilde \tau_2\rangle\arrow[r]\arrow[d]& H=\langle \tau_1, \tau_2\rangle\arrow[r]\arrow[d]& \{1\}\\
\{1\} \arrow[r] & K:=\langle {\widetilde c}^2\rangle\arrow[r]& \widetilde W\arrow[r]\arrow[d]& W\arrow[r]\arrow[d]& \{1\}\\
& & \ZZ/2\ZZ \arrow[r, equal]\arrow[d] & \ZZ/2\ZZ\arrow[d]& \\
& & \{1\} & \{1\}& 
\end{tikzcd},
\]
where 
\begin{gather*}
\widetilde \tau_1(\alpha, \delta_1,\delta_2,\kappa_1,\kappa_2):=
\widetilde r_{\alpha_2}\widetilde r_{\alpha_3}(\alpha, \delta_1,\delta_2,\kappa_1,\kappa_2)
=(\alpha-2\delta_1, \delta_1,\delta_2,\kappa_1+\alpha-\delta_1-\delta_2,\kappa_2-\delta_1),\\
\widetilde \tau_2(\alpha, \delta_1,\delta_2,\kappa_1,\kappa_2):=
\widetilde r_{\alpha_2}\widetilde r_{\alpha_1}(\alpha, \delta_1,\delta_2,\kappa_1,\kappa_2)
=(\alpha-2\delta_2, \delta_1,\delta_2,\kappa_1+\delta_2,\kappa_2-\alpha+\delta_1+\delta_2),
\end{gather*}
which satisfies
\[
\widetilde \tau_1\widetilde \tau_2{\widetilde \tau}_1^{-1}{\widetilde \tau}_2^{-1}={\widetilde c}^{-2}.
\]

Let $\langle -,-\rangle$ denote the natural pairing between $\widetilde\N\otimes_\ZZ\CC$ and $\Hom_\ZZ(\widetilde\N,\CC)$, 
and $\{\alpha^\vee, d_1, d_2, c_1, c_2\}$ a $\CC$-basis of $\Hom_\ZZ(\widetilde\N,\CC)$ such that 
the values of the pairing between bases are given by
\[
\langle \alpha,\alpha^\vee\rangle=2, \quad \langle \delta_1,d_1\rangle=1,\quad \langle \delta_2,d_2\rangle=1,\quad \langle \kappa_1,c_1\rangle=1, \quad \langle \kappa_2,c_2\rangle=1,
\]
and zero otherwise.
Let $(x, \omega_1,  \omega_2, \phi', \phi)$ denote a system of coordinates associated to the basis above. Set
\[
\widetilde\Pi:=\{[\widetilde h]\in \PP\Hom_\ZZ(\widetilde\N,\CC)\,|\,
\widetilde{I}(\widetilde h,\widetilde h)=0,\ 
\sqrt{-1}(\omega_1\overline{\omega_2}-\overline{\omega_1}\omega_2)>0 \}
\]
where we denote by the same letter $\widetilde I$ the symmetric $\CC$-bilinear form induced by $\widetilde I$ on $\widetilde \N$ 
and $\omega_i:=\langle\delta_i,\widetilde h\rangle$, $i=1,2$, and
\[
\widetilde\EE:=\{x\alpha^\vee+d_1+\tau d_2+\phi c_2\in \Hom_\ZZ(\widetilde\N,\CC)\,|\,(\phi,x,\tau)\in\CC\times \CC\times \HH\},
\]
where $\tau:=\omega_{2}/\omega_{1}$.
\begin{prop}\label{prop: Pi}
The correspondence
\[
x\alpha^\vee+d_1+\tau d_2+\phi c_2 \mapsto 
[x\alpha^\vee+d_1+\tau d_2+(\phi\tau-x^2)c_1+\phi c_2]
\]
gives the isomorphism of complex manifold between $\widetilde \EE$ and $\widetilde\Pi$.
\end{prop}
\begin{pf}
The coefficient $\phi\tau-x^2$ of $c_{1}$ comes from $\widetilde{I}(\widetilde h,\widetilde h)=0$.
The coefficient $1$ of $d_{1}$ fixes the ambiguity of the projection $\widetilde\Pi\rightarrow \widetilde \EE$ and is clearly holomorphic. 
The correspondence is also holomorphic. These yield the proposition.
\qed
\end{pf}

We extend the action on $K_0(\D)$ to $\widetilde{\N}$ so that it is $\CC$-linear and $\widetilde{I}$-invariant. 
Under the expression \eqref{widetilde I} of $\widetilde{I}$, 
it is straightforward to see that we have the following matrix representation of the action:
\begin{equation}
\begin{pmatrix}
0 & 0 & 0 & d & -b\\
0 & 0 & 0 & -c & a\\
0 & 0 & 1 & 0 & 0\\
d &- b & 0 & 0 & 0\\
-c & a & 0 & 0 & 0
\end{pmatrix},\quad 
\begin{pmatrix}
a&b\\
c&d
\end{pmatrix}\in {\rm SL}(2;\ZZ).
\end{equation}
Therefore, the space $\Hom_\ZZ(\widetilde\N,\CC)$ admits the ${\rm SL}(2;\ZZ)$-action defined by
\begin{equation}
\langle \widetilde \lambda, \gamma (h)\rangle:=\langle \gamma^{-1}(\widetilde \lambda), h\rangle,\quad \widetilde \lambda\in \widetilde \N,\ h\in \Hom_\ZZ(\widetilde\N,\CC),\ 
\gamma\in{\rm SL}(2;\ZZ),
\end{equation}
which is explicitly written as 
\begin{multline}\label{eq: SLPi}
\gamma[x\alpha^\vee+\omega_1 d_1+\omega_2 d_2+\phi' c_1+\phi c_2]\\=
[x\alpha^\vee+(c\omega_2+d\omega_1) d_1+(a \omega_2+b \omega_1) d_2+(a \phi'+b \phi) c_1+(c \phi'+d \phi ) c_2],
\end{multline}
for $\gamma=\begin{pmatrix}
a&b\\
c&d
\end{pmatrix}
\in{\rm SL}(2;\ZZ)$.

We shall often identify $\widetilde\EE$ with $\CC\times \CC\times \HH$, namely, identify an element 
$x\alpha^\vee+d_1+\tau d_2+\phi c_2\in \Hom_\ZZ(\widetilde\N,\CC)$ 
with a point $(\phi,x,\tau)\in\CC\times \CC\times \HH$, and write $(\phi,x,\tau)\in\widetilde \EE$ for simplicity.

The natural action of $\widetilde W$ on the space $\widetilde \EE$ is given by
\[
\langle \widetilde \lambda, \widetilde w(h)\rangle:=\langle \widetilde{w}^{-1}(\widetilde \lambda), h\rangle, 
\quad \widetilde \lambda\in \widetilde \N,\ h\in \widetilde \EE,\ \widetilde w\in\widetilde W.
\]

Explicitly, we have
\begin{subequations}\label{eq: action on E-tilde}
\begin{gather}
\widetilde r_{\alpha_1}(\phi,x,\tau)=(\phi, -x+1, \tau),\\
\widetilde r_{\alpha_2}(\phi,x,\tau)=(\phi-2x+1+\tau,-x+1+\tau, \tau),\\
\widetilde r_{\alpha_3}(\phi,x,\tau)=(\phi-2x+\tau, -x+\tau, \tau),
\end{gather}
and hence we have
\begin{gather}
\widetilde \tau_1(\phi,x,\tau)=(\phi+1, x+1, \tau),\\
\widetilde \tau_2(\phi,x,\tau)=(\phi+2x-1+\tau, x+\tau, \tau),\\
\widetilde \tau_1^m\widetilde \tau_2^n(\phi,x,\tau)=(\phi+2nx+n^2\tau+m-n,x+m+n\tau,\tau),\\
\widetilde c(\phi,x,\tau)=(\phi+1,-x,\tau).
\end{gather}
\end{subequations}
By Proposition~\ref{prop: Pi} and \eqref{eq: SLPi}, we obtain the ${\rm SL}(2;\ZZ)$-action on $\widetilde{\EE}$, which is given by
\[
\gamma (\phi,x,\tau)=\left(\phi-\frac{cx^{2}}{c\tau+d}, \frac{x}{c\tau+d}, \frac{a\tau+b}{c\tau+d}\right),
\]

\begin{rem}
Let us summarize the relationships among the groups that have appeared so far.
The group $\widehat W$ which fits into the following commutative diagram of short exact sequences is the natural generalization of 
the extended affine Weyl group introduced by Dubrovin--Zhang~\cite{DZ}:
\[
\begin{tikzcd}
& \{1\}\arrow[d]& \{1\}\arrow[d] &  & \\
\{1\} \arrow[r] & \ZZ=K\arrow[r]\arrow[d]& \widetilde W\arrow[r]\arrow[d]& W\arrow[r]\arrow[d, equal]& \{1\}\\
\{1\} \arrow[r] & B_3\arrow[r]\arrow[d]& \widehat W\arrow[r]\arrow[d]& W\arrow[r]& \{1\}\\
& {\rm SL}(2;\ZZ)\arrow[r, equal]\arrow[d]& {\rm SL}(2;\ZZ)\arrow[d] & & \\
& \{1\}& \{1\} & & 
\end{tikzcd}.
\]
Since it is more convenient to study the $W$-invariants and the ${\rm SL}(2;\ZZ)$-action on them via the variable $e^{\pi\sqrt{-1}\phi}$, 
namely, by passing to the quotient by $\ZZ$ rather than by lifting the action to $B_3$, 
we do not develop the precise definition of $\widehat W$ or the corresponding invariant theory in this paper.
\end{rem}

It is also important that $\widetilde W$ acts trivially on the third component $\tau$, 
and for each root $\alpha-(p+1)\delta_1-(q+1)\delta_2\in\Delta_{re}$ the reflection hyperplane in $\widetilde \EE$ defined by 
\begin{equation*}\label{hyperplane}
\{h\in\widetilde \EE\,|\,\langle \alpha-p\delta_1-q\delta_2,h\rangle=0\}
\end{equation*}
is given by $x=\frac{p}{2}+\frac{q}{2}\tau$.

From these explicit forms of the $\widetilde W$-action, we see the following
\begin{prop}[{\cite[Section~3.5]{Sa2}}]
Let 
\[
\EE:=\{x\alpha^\vee+d_1+\tau d_2\in \Hom_\ZZ(\N,\CC)\,|\,(x,\tau)\in\CC\times \HH\},
\]
and consider the natural projection map $\widetilde \EE\longrightarrow \EE$, $(\phi,x,\tau)\mapsto (x,\tau)$. 
\begin{enumerate}
\item
The following diagrams commute:
\[
\begin{tikzcd}
\widetilde W\times \widetilde \EE\arrow [r]\arrow[d] & W\times \EE \arrow[d]\\
\widetilde \EE\arrow[r] & \EE
\end{tikzcd},\quad
\begin{tikzcd}
\widetilde H\times \widetilde \EE\arrow [r]\arrow[d] & H\times \EE \arrow[d]\\
\widetilde \EE\arrow[r] & \EE
\end{tikzcd}.
\]
\item 
The action of $\widetilde H$ on $\widetilde \EE$ is properly discontinuous and fixed point free.
\item
The complex manifold $\widetilde \EE/K$ is the trivial $\CC^*$-bundle over $\EE$ and 
the complex manifold $\widetilde \EE/\widetilde H$ is a principal $\CC^*$-bundle over $\EE/H$.
\item
The action of $\widetilde W$ on $\widetilde \EE$ is properly discontinuous.
\item
Let $L$ be the line bundle over $\EE/H$ associated to $\widetilde \EE/\widetilde H$, which is as a set, a union $L=\widetilde \EE/\widetilde H \cup \EE/H$.
The finite group $\widetilde W/\widetilde H\cong W/H\cong \ZZ/2\ZZ$ acts equivariantly on both $L$ and $\EE/ H$.
\item
For each $\tau\in\HH$, the degree of the Chern class of the line bundle $L|_{{\rm E}_\tau}$ is $-1$ where 
${\rm E}_\tau:=\CC/(\ZZ+\ZZ\tau)$ is the fiber of the natural projection $\EE/H\longrightarrow \HH$ over $\tau\in\HH$.
\end{enumerate}
\qed\end{prop}

\subsection{Elliptic Weyl group invariant theory and Frobenius structures}\label{Invariant-Isom}
The purpose of this subsection is to show an isomorphism of Frobenius structures on $M$ between the one constructed from the deformation theory of $F$ with the primitive form $[2\pi\sqrt{-1}dz]$ and the one constructed from the invariant theory of an elliptic Weyl group $\widetilde{W}$.
The relation between the Weierstrass elliptic functions and the theta functions will play important roles 
and relate these two constructions.

The holomorphic function $e^{\pi\sqrt{-1}\phi}$ on $\widetilde \EE$ is $K$-invariant and 
the action of $\tau_1^m \tau_2^n\in H=\widetilde H/K$ on $e^{\pi\sqrt{-1}\phi}$ is calculated as 
\[
\tau_1^m \tau_2^n(e^{\pi\sqrt{-1}\phi})=(-1)^{m-n}{\bf e}\left[\frac{n^2}{2}\tau+n x\right]e^{\pi\sqrt{-1}\phi}.
\]
Define a $\ZZ$-graded algebra $S$ over $\Gamma(\HH,\O_\HH)$ by 
\[
S:=\bigoplus_{k=0}^\infty S_k,\quad S_k:=\Gamma(\EE/H,\O(L^{-1})^{\otimes k}).
\]
Note that $S_k$ is the set of $\widetilde H$-invariant holomorphic functions on $\widetilde \EE$
satisfying the degree condition
\[
\frac{1}{\pi\sqrt{-1}}\frac{\p y(\phi,x,\tau)}{\p \phi}=k\cdot y(\phi,x,\tau),\quad y(\phi,x,\tau)\in S_k.
\]
Therefore, an element $y(\phi,x,\tau)\in S_k$ is described as a product $y(\phi,x,\tau)=y(x,\tau)e^{k\pi\sqrt{-1}\phi}$ where
$y(x,\tau)$ is a holomorphic function satisfying 
\[
\tau_1^m \tau_2^n(y(x,\tau))=(-1)^{k(m-n)}{\bf e}\left[-\frac{kn^2}{2}\tau-kn x\right]y(x,\tau).
\]

Since the group $\widetilde W/\widetilde H$ acts equivariantly on both $\widetilde \EE/\widetilde H$ and the line bundle $L$, 
it induces a natural action on $S_k$. 
Denote by $S_k^{\inv}$ $(\text{resp.}\, S_k^{-\inv})$ the set of $\widetilde W/\widetilde H$-invariant (resp. anti-invariant) elements of $S_k$, the set of $\widetilde{W}$-invariant (resp. anti-invariant) holomorphic functions on $\widetilde{\EE}$ 
satisfying the above degree condition. Set
\[
S^{\inv}:=\bigoplus_{k=0}^\infty S_k^{\inv},\quad S^{-\inv}:=\bigoplus_{k=0}^\infty S_k^{-\inv}.
\]
Naturally, $S^{\inv}$ is a $\Gamma(\HH,\O_\HH)$-algebra and $S^{-\inv}$ is an $S^{\inv}$-module.

Consider the theta function
\[
\vartheta_{11}(x;\tau):=\sum_{n\in\ZZ} {\bf e}\left[\frac{1}{2}\left(n+\frac{1}{2}\right)^2\tau +\left(n+\frac{1}{2}\right)\left(x+\frac{1}{2}\right)\right],
\]
which satisfies $\vartheta_{11}(-x;\tau)=-\vartheta_{11}(x;\tau)$ and  
\begin{equation}\label{eq: theta lattice}
\vartheta_{11}(x+m+n\tau;\tau)=(-1)^{m+n}{\bf e}\left[-\frac{n^2}{2}\tau-n x\right]\vartheta_{11}(x;\tau),\quad m,n\in\ZZ.
\end{equation}
Define $y_1\in S_2^{\inv}$, $y_2\in S_1^{\inv}$, $y_3\in S_0^{\inv}=\Gamma(\HH,\O_\HH)$ and $\J\in S_3^{-\inv}$ by
\begin{gather*}
y_1:=\frac{\vartheta_{11}(x;\tau)^2}{\vartheta_{11}'(0;\tau)^2} \widetilde\wp(x;\tau){\bf e}[\phi],\quad 
y_2:=\frac{\vartheta_{11}(x;\tau)}{\vartheta_{11}'(0;\tau)}{\bf e}\left[\frac{1}{2}\phi\right],\quad 
y_3:=2\pi\sqrt{-1}\tau,\\
\J:=-\frac{1}{2}\frac{1}{2\pi\sqrt{-1}}\frac{\vartheta_{11}(x;\tau)^3}{\vartheta_{11}'(0;\tau)^3} \frac{\p \widetilde\wp(x;\tau)}{\p x}{\bf e}\left[\frac{3}{2}\phi\right],
\end{gather*}
where 
\[
\vartheta_{11}'(0;\tau):=\left.\frac{\p\vartheta_{11}(x;\tau)}{\p x}\right|_{x=0}.
\]
Then we have
\[
dy_1\wedge dy_2\wedge dy_3=(2\pi\sqrt{-1})^3 \J d\phi\wedge dx\wedge d\tau.
\]

The following proposition highlights the importance of considering invariants with respect to the elliptic Weyl group rather than the Jacobi group. 
Otherwise, the discriminant does not become a cubic polynomial in the invariant with the highest degree.
\begin{prop}[Chevalley type theorem, see {\cite[Section~4.5]{Sa2}} and references therein]\label{prop: Chevalley}
The following hold:
\begin{enumerate}
\item
The $\ZZ$-graded algebra $S$ is an algebra over $\Gamma(\HH,\O_\HH)$ generated by $y_1, y_2, \J$ 
with the unique relation
\begin{equation}\label{eq: ell to P123}
\J^2=y_1^3-\frac{1}{48}E_4(\tau)y_1y_2^2+\frac{1}{864}E_6(\tau)  y_2^3.
\end{equation}
\item
The $\ZZ$-graded algebra $S^{\inv}$ is a polynomial algebra over $\Gamma(\HH,\O_\HH)$ generated by $y_1,y_2$.

\item
The $S^{\inv}$-module $S^{-\inv}$ is free of rank one generated by $\J$.

\item
The zero locus of $\J$ on $\widetilde \EE$ is equal to the union of the reflection hyperplanes $H_\beta$, $\beta\in \Delta_{re}$. Moreover, 
\[
\J=\frac{1}{(2\pi\sqrt{-1})^3}\frac{\vartheta_{11}(2x;\tau)}{2\vartheta_{11}(x;\tau)} {\bf e}\left[\frac{3}{2}\phi\right]
\]
and it has the following product formula
\begin{eqnarray*}
\J&=&\frac{1}{(2\pi\sqrt{-1})^3}\cdot e^{3\pi\sqrt{-1}\phi}\cdot \prod_{n=1}^\infty \left(1+{\bf e}\left[n\tau+x\right]\right)\left(1+{\bf e}\left[n\tau-x\right]\right)\\
& &\cdot \prod_{n=1}^\infty\left(1-{\bf e}\left[\left(n+\frac{1}{2}\right)\tau+x\right]\right)
\left(1-{\bf e}\left[\left(n+\frac{1}{2}\right)\tau-x\right]\right)\\
& &\cdot \frac{1}{2}\left({\bf e}\left[\frac{x}{2}\right]+{\bf e}\left[-\frac{x}{2}\right]\right)
\prod_{n=1}^\infty\left(1+{\bf e}\left[\left(n+\frac{1}{2}\right)\tau+x\right]\right)\left(1+{\bf e}\left[\left(n+\frac{1}{2}\right)\tau-x\right]\right).
\end{eqnarray*}
\end{enumerate}
\end{prop}
\begin{rem}
Up to a non-zero constant multiple the function $\vartheta_{11}(2x;\tau)$ is the denominator of the affine $A_1$-type 
with $\alpha_1=\alpha-\delta_1$ and $\alpha_2=-\alpha+\delta_1+\delta_2$ as the simple roots.
Since $\alpha\notin\Delta_{re}$, the infinite product obtained by dividing $\vartheta_{11}(2x;\tau)$ by $\vartheta_{11}(x;\tau)$ appears as $\J$.
\end{rem}
\begin{rem}
Since the line bundle $L$ is negative relative to $\HH$, the zero section of the line bundle $L$ can be blow down to $\HH$,  
The blow-down space $\LL=\widetilde \EE/\widetilde H\cup \HH$ is isomorphic to ${\rm Spec}(S)$.
Note also that for each $\tau\in\HH$ the $\ZZ$-graded $\CC$-algebra $S$ gives an embedding of an elliptic curve ${\rm E}_\tau$ 
into the weighted projective space $\PP(1,2,3)$.
\end{rem}
\begin{rem}
The reason why we choose $y_{1}, y_{2}$ as in \eqref{eq: ell to P123} will be explained around \eqref{normalized}.
\end{rem}
\begin{pf}
For (i),  see \cite[Section~A.11]{Sa2}. (ii) and (iii) are direct consequence of (i). 
As for the proof of (iv), we have
\[
\frac{\vartheta_{11}(x;\tau)^3}{\vartheta_{11}'(0;\tau)^3} \left(-\frac{\p}{\p x}\widetilde\wp(x;\tau)\right)=\frac{1}{(2\pi\sqrt{-1})^2}\frac{\vartheta_{11}(2x;\tau)}{\vartheta_{11}(x;\tau)},
\]
which can be checked by using \eqref{eq: theta lattice}.
The last part follows from the 
Jacobi's triple product formula
\begin{equation}\label{eq: theta product}
\vartheta_{11}(x;\tau)=\sqrt{-1}{\bf e}\left[\frac{\tau}{8}\right]\left({\bf e}\left[\frac{x}{2}\right]-{\bf e}\left[-\frac{x}{2}\right]\right)
\prod_{n=1}^\infty \left(1-{\bf e}\left[n\tau\right]\right) \left(1-{\bf e}\left[n\tau+x\right]\right)\left(1-{\bf e}\left[n\tau-x\right]\right).
\end{equation}
\qed
\end{pf}
Recall briefly the relationship between the Weierstrass's elliptic functions and the theta function $\vartheta_{11}(x;\tau)$,
which will be used later. The Weierstrass's $\sigma$-function satisfies the identity
\[
\sigma(x;1,\tau):=e^{\frac{\pi^2}{6}E_2(\tau)x^2}\frac{\vartheta_{11}(x;\tau)}{\vartheta_{11}'(0;\tau)}
\]
which can be checked by using \eqref{eq: theta lattice} together with the Legendre's relation.
Then 
\[
-\widetilde\zeta(x;\tau)-\frac{1}{12}E_2(\tau)x=-\frac{1}{(2\pi\sqrt{-1})^2}\frac{\p}{\p x}\log\vartheta_{11}(x;\tau),
\]
and hence
\[
\widetilde\wp(x;\tau)-\frac{1}{12}E_2(\tau)=-\frac{1}{(2\pi\sqrt{-1})^2}\frac{\p^2}{\p x^2}\log\vartheta_{11}(x;\tau).
\]
Moreover, it follows from the heat equation
\begin{equation}\label{eq: theta heat}
\frac{1}{2\pi\sqrt{-1}}\frac{\p}{\p\tau}\vartheta_{11}(x;\tau)=\frac{1}{2}\frac{1}{(2\pi\sqrt{-1})^2}\frac{\p^2}{\p x^2}\vartheta_{11}(x;\tau),
\end{equation}
that  
\begin{equation}\label{eq: theta E2}
E_2(\tau)=\frac{4}{(2\pi\sqrt{-1})^2}\frac{\vartheta_{11}'''(0;\tau)}{\vartheta_{11}'(0;\tau)},\quad
\vartheta_{11}'''(0;\tau):=\left.\frac{\p^3\vartheta_{11}(x;\tau)}{\p x^3}\right|_{x=0},
\end{equation}
where we used the fact that in terms of the Dedekind's eta function one has 
\[
\frac{1}{24}E_2(\tau)=\frac{1}{2\pi\sqrt{-1}}\frac{\p}{\p \tau}\log \eta(\tau),\quad 
\eta(\tau):={\bf e}\left[\frac{\tau}{24}\right]\prod_{n=1}^\infty \left(1-{\bf e}\left[n\tau\right]\right),
\]
and $\vartheta_{11}'(0;\tau)=-2\pi\eta(\tau)^3$ by \eqref{eq: theta product}.
It is also important that 
\begin{equation}\label{normalized}
\frac{\vartheta_{11}\left(\frac{x}{c\tau+d};\frac{a\tau+b}{c\tau+d}\right)}{\vartheta_{11}'\left(0;\frac{a\tau+b}{c\tau+d}\right)}
{\bf e}\left[\frac{1}{2}\left(\phi-\frac{cx^2}{c\tau+d}\right)\right]
=\frac{1}{c\tau+d}\frac{\vartheta_{11}\left(x;\tau\right)}{\vartheta_{11}'\left(0;\tau\right)}
{\bf e}\left[\frac{1}{2}\phi\right],\quad 
\begin{pmatrix}
a & b\\
c & d
\end{pmatrix}\in {\rm SL}(2;\ZZ),
\end{equation}
which is one of the reasons why we normalized $y_1$, $y_2$ and $\J$ in terms of $\vartheta_{11}'\left(0;\tau\right)$.
The equation \eqref{eq: key identity-tilde-1} is also a direct consequence of the heat equation \eqref{eq: theta heat}.

Let 
\begin{eqnarray}\label{EEo}
\widetilde \EE^o&:=&\widetilde \EE\setminus \{h\in \widetilde \EE\,|\, \langle \alpha+2k\delta_1+2l\delta_2,h\rangle =0,\ k,l\in\ZZ\}\\
&=&\{x\alpha^\vee+d_1+\tau d_2+\phi c_2\in \Hom_\ZZ(\widetilde\N,\CC)\,|\,(\phi,x,\tau)\in
\CC\times \left(\CC\setminus \left(\ZZ+\ZZ\tau\right)\right)\times \HH\}. 
\end{eqnarray}

Define a complex manifold $\widetilde \EE/\!/\widetilde W$ as 
\[
\widetilde \EE/\!/\widetilde W:={\rm Specan}(S^{\inv}).
\]
Due to Proposition~\ref{prop: Chevalley}, $\widetilde \EE/\!/\widetilde W\cong \CC\times \CC\times \HH$.
Then we have
\[
\CC\times\CC^*\times\HH\cong \widetilde \EE^o/\widetilde W\subsetneq \widetilde \EE/\widetilde W \subsetneq \widetilde \EE/\!/\widetilde W\cong \CC\times \CC\times \HH.
\]
Note that $(0,0,\tau)\notin \widetilde \EE/\widetilde W$.

We have $M=\CC\times\CC^*\times\HH\cong \widetilde \EE^o/\widetilde W$ as complex manifolds.
From now on, we shall construct a Frobenius structure on the tangent bundle of $\widetilde \EE/\!/\widetilde W$ via the invariant theory of $\widetilde W$ and the intersection form
and compare this Frobenius structure restricted to $M$ with the one in Corollary~\ref{primitive_form}.

Consider the $\widetilde W$-invariant holomorphic functions $t_1, t_2, t_3$ on $\widetilde \EE$ defined by 
\begin{eqnarray*}
t_1(\phi,x,\tau)&:=&-y_1+\frac{1}{12}E_2(\tau)y_2^2
=-\frac{\vartheta_{11}(x;\tau)^2}{\vartheta_{11}'(0;\tau)^2} \left(\widetilde\wp(x;\tau)-\frac{1}{12}E_2(\tau)\right)e^{2\pi\sqrt{-1}\phi}\\
&=&\frac{1}{(2\pi\sqrt{-1})^2}
\frac{\vartheta_{11}''(x;\tau)\vartheta_{11}(x;\tau)-\vartheta_{11}'(x;\tau)^2}{\vartheta_{11}'(0;\tau)^2}e^{2\pi\sqrt{-1}\phi},
\end{eqnarray*}
\[
t_2(\phi,x,\tau):=y_2=\frac{\vartheta_{11}(x;\tau)}{\vartheta_{11}'(0;\tau)}e^{\pi\sqrt{-1}\phi},\quad 
t_3(\phi,x,\tau):=y_3=2\pi\sqrt{-1}\tau,
\]
which forms a global coordinate system on $\widetilde \EE/\!/\widetilde W$. 

\begin{rem}
The holomorphic functions $t_{1}$, $t_{2}$ are good invariants in the sense of \cite{Sat2}, namely, 
their higher derivatives with respect to $x$ and $\phi$ vanish at $x=0$.
\end{rem}

Since $\widetilde \N$ can be identified with the cotangent space of $\Hom_\ZZ(\widetilde\N,\CC)$ at $0$, 
it turns out that $\widetilde I$ induces the $\widetilde W$-invariant non-degenerate symmetric bilinear form on 
the cotangent sheaf $\Omega_{\widetilde \EE}$ on $\widetilde \EE$ in the following way:
\begin{prop}
Define a symmetric $\O_{\widetilde \EE}$-bilinear form $g_{\widetilde W}:\Omega^1_{\widetilde \EE}\times \Omega^1_{\widetilde \EE}\longrightarrow \O_{\widetilde \EE}$ by the following matrix representation with respect to the $\O_{\widetilde \EE}$-basis 
$\{d\phi, dx,d\tau\}$:
\begin{equation}
g_{\widetilde W}:=\frac{1}{(2\pi\sqrt{-1})^2}\begin{pmatrix}
0 & 0 & 1\\
0 & -\frac{1}{2} & 0\\
1 & 0 & 0
\end{pmatrix}.
\end{equation}

Then $g_{\widetilde W}$ is $\widetilde W$-invariant.
\end{prop}
\begin{pf}
This is clear from \eqref{eq: action on E-tilde}.
\qed\end{pf}

Therefore, $g_{\widetilde W}$ naturally induces the non-degenerate symmetric $\O_{\widetilde \EE/\!/\widetilde W}$-bilinear form on $\Omega^1_{\widetilde \EE/\!/\widetilde W}$, 
which we shall denote by the same symbol $g_{\widetilde W}$.
\begin{prop}\label{prop: invariant g}
The non-degenerate symmetric $\O_{\widetilde \EE/\!/\widetilde W}$-bilinear form $g_{\widetilde W}$ on $\Omega^1_{\widetilde \EE/\!/\widetilde W}$ is given by
\[
(g_{\widetilde W}(dt_i , dt_j))=
\begin{pmatrix}
-\frac{t_2^4}{12}\frac{1}{(2\pi\sqrt{-1})^2}\frac{\p^2 E_2(\tau)}{\p \tau^2} & -\frac{t_2^3}{8}\frac{1}{2\pi\sqrt{-1}}\frac{\p E_2(\tau)}{\p \tau}   & t_1\\
-\frac{t_2^3}{8}\frac{1}{2\pi\sqrt{-1}}\frac{\p E_2(\tau)}{\p \tau}   & \frac{1}{2}t_1-\frac{1}{8}t_2^2E_2(\tau) & \frac{1}{2}t_2\\
t_1 & \frac{1}{2}t_2 & 0
\end{pmatrix}.
\]
In particular, each entry is at most linear in $t_1$ and 
the new symmetric $\O_{\widetilde \EE/\!/\widetilde W}$-bilinear form $\eta_{\widetilde W}$ on 
$\Omega^1_{\widetilde \EE/\!/\widetilde W}$ defined by
\[
(\eta_{\widetilde W}(dt_i , dt_j)):=\left(\frac{\p g_{\widetilde W}(dt_i , dt_j)}{\p t_1}\right)=
\begin{pmatrix}
0 & 0 & 1\\
0 & \frac{1}{2} & 0\\
1 & 0 & 0
\end{pmatrix}
\]
is non-degenerate.
\end{prop}
\begin{pf}
The proof is done by direct calculation using identities of theta functions and $\wp$-functions.
{\tiny
\[
g_{\widetilde W}(dt_1,dt_3)=\frac{1}{(2\pi\sqrt{-1})^2}\frac{\p t_1}{\p\phi}(2\pi\sqrt{-1})=t_1,\quad g_{\widetilde W}(dt_2,dt_3)=\frac{1}{(2\pi\sqrt{-1})^2}\frac{\p t_2}{\p\phi}(2\pi\sqrt{-1})=\frac{1}{2}t_2,\quad 
g_{\widetilde W}(dt_3,dt_3)=0,
\]}
{\tiny
\begin{eqnarray*}
g_{\widetilde W}(dt_2,dt_2) &=& -\frac{1}{2}\frac{1}{(2\pi\sqrt{-1})^2}\frac{\p t_2}{\p x}\frac{\p t_2}{\p x}+2\frac{1}{(2\pi\sqrt{-1})^2}\frac{\p t_2}{\p \tau}\frac{\p t_2}{\p \phi}\\
&=&-\frac{1}{2}\frac{1}{(2\pi\sqrt{-1})^2}\frac{\vartheta_{11}'(x;\tau)^2}{\vartheta_{11}'(0;\tau)^2}e^{2\pi\sqrt{-1}\phi}+\frac{1}{2}\frac{1}{(2\pi\sqrt{-1})^2}\left(\frac{\vartheta_{11}''(x;\tau)\vartheta_{11}(x;\tau)}{\vartheta_{11}'(0;\tau)^2}-\frac{\vartheta_{11}'''(0;\tau)\vartheta_{11}(x;\tau)^2}{\vartheta_{11}'(0;\tau)^3}\right)e^{2\pi\sqrt{-1}\phi}\\
&=&\frac{1}{2}t_1-\frac{1}{8}t_2^2E_2(\tau),
\end{eqnarray*}}
{\tiny
\begin{eqnarray*}
g_{\widetilde W}(dt_1,dt_2) &=& -\frac{1}{2}\frac{1}{(2\pi\sqrt{-1})^2}\frac{\p t_1}{\p x}\frac{\p t_2}{\p x}+\frac{1}{(2\pi\sqrt{-1})^2}\frac{\p t_1}{\p \tau}\frac{\p t_2}{\p \phi}
+\frac{1}{(2\pi\sqrt{-1})^2}\frac{\p t_1}{\p \phi}\frac{\p t_2}{\p \tau}\\
&=&\frac{1}{(2\pi\sqrt{-1})^2}t_2\frac{\p t_2}{\p x}\frac{\p t_2}{\p x}\left(\widetilde\wp(x;\tau)-\frac{1}{12}E_2(\tau)\right)
+\frac{1}{2}\frac{1}{(2\pi\sqrt{-1})^2}t_2^2\frac{\p t_2}{\p x}\frac{\p }{\p x}\left(\widetilde\wp(x;\tau)-\frac{1}{12}E_2(\tau)\right)\\
& &-2\frac{1}{2\pi\sqrt{-1}}t_2^2\frac{\p t_2}{\p \tau}\left(\widetilde\wp(x;\tau)-\frac{1}{12}E_2(\tau)\right)
-\frac{1}{2}\frac{1}{2\pi\sqrt{-1}}t_2^3\frac{\p}{\p \tau}\left(\widetilde\wp(x;\tau)-\frac{1}{12}E_2(\tau)\right)\\
&=&t_2^3\left(\widetilde\wp(x;\tau)-\frac{1}{12}E_2(\tau)\right)^2+\frac{t_2^3}{4}E_2(\tau)\left(\widetilde\wp(x;\tau)-\frac{1}{12}E_2(\tau)\right)\\
& &-\frac{1}{2}t_2^3\left(-\widetilde\zeta(x;\tau)-\frac{1}{12}E_2(\tau)x\right)\frac{\p }{\p x}\left(\widetilde\wp(x;\tau)-\frac{1}{12}E_2(\tau)\right)
-\frac{t_2^3}{2}\frac{1}{2\pi\sqrt{-1}}\frac{\p}{\p \tau}\left(\widetilde\wp(x;\tau)-\frac{1}{12}E_2(\tau)\right)\\
&=&-\frac{1}{8}t_2^3\frac{1}{2\pi\sqrt{-1}}\frac{\p E_2(\tau)}{\p \tau}.
\end{eqnarray*}}
{\tiny
\begin{eqnarray*}
g_{\widetilde W}(dt_1,dt_1) &=& -\frac{1}{2}\frac{1}{(2\pi\sqrt{-1})^2}\frac{\p t_1}{\p x}\frac{\p t_1}{\p x}+2\frac{1}{(2\pi\sqrt{-1})^2}\frac{\p t_1}{\p \tau}\frac{\p t_1}{\p \phi}\\
&=& -\frac{1}{2}\frac{1}{(2\pi\sqrt{-1})^2}\left(-2t_2\frac{\p t_2}{\p x}\left(\widetilde\wp(x;\tau)-\frac{1}{12}E_2(\tau)\right)-t_2^2\frac{\p }{\p x}\left(\widetilde\wp(x;\tau)-\frac{1}{12}E_2(\tau)\right)\right)^2\\
& &+4t_2^3\frac{1}{2\pi\sqrt{-1}}\frac{\p t_2}{\p \tau}\left(\widetilde\wp(x;\tau)-\frac{1}{12}E_2(\tau)\right)^2
+2t_2^4\left(\widetilde\wp(x;\tau)-\frac{1}{12}E_2(\tau)\right)\frac{1}{2\pi\sqrt{-1}}\frac{\p}{\p \tau}\left(\widetilde\wp(x;\tau)-\frac{1}{12}E_2(\tau)\right)\\
&=& -2t_2^2\left(\widetilde\wp(x;\tau)-\frac{1}{12}E_2(\tau)\right)^3-\frac{t_2^2}{2}E_2(\tau)\left(\widetilde\wp(x;\tau)-\frac{1}{12}E_2(\tau)\right)^2\\
& & +2t_2^4\frac{1}{(2\pi\sqrt{-1})^2}\left(-\widetilde\zeta(x;\tau)-\frac{1}{12}E_2(\tau)x\right)\left(\widetilde\wp(x;\tau)-\frac{1}{12}E_2(\tau)\right)\frac{\p }{\p x}\left(\widetilde\wp(x;\tau)-\frac{1}{12}E_2(\tau)\right)\\
& &-\frac{1}{2}\frac{1}{(2\pi\sqrt{-1})^2}t_2^4\left(\frac{\p }{\p x}\left(\widetilde\wp(x;\tau)-\frac{1}{12}E_2(\tau)\right)\right)^2
+2t_2^4\left(\widetilde\wp(x;\tau)-\frac{1}{12}E_2(\tau)\right)\frac{1}{2\pi\sqrt{-1}}\frac{\p}{\p \tau}\left(\widetilde\wp(x;\tau)-\frac{1}{12}E_2(\tau)\right)\\
&=&-\frac{t_2^4}{12}\frac{1}{(2\pi\sqrt{-1})^2}\frac{\p^2 E_2(\tau)}{\p\tau^2}.
\end{eqnarray*}}
\qed\end{pf}

Let $\widetilde\nabla$ be the Levi--Civita connection for $g_{\widetilde W}$ on $\Omega^1_{\widetilde \EE}$ and 
consider the contravariant components of the Christoffel symbols 
\[
\Gamma^{ij}_{k}:=g_{\widetilde W}\left(dt_i,\widetilde\nabla_{\frac{\p}{\p t_k}}dt_j\right),\quad i,j,k=1,2,3.
\] 
Let $g^{ij}:=g_{\widetilde W}(dt_i,dt_j)$. It follows from the definition that for $i,j,k=1,2,3$ 
\begin{subequations}
\begin{gather}
\frac{\p g^{ij}}{\p t_k}=\Gamma^{ij}_k+\Gamma^{ji}_k\label{eq: contra a},\\ 
\sum_{a=1}^3g^{ia}\Gamma^{jk}_a=\sum_{a=1}^3g^{ja}\Gamma^{ik}_a.\label{eq: contra b}
\end{gather}
\end{subequations}
and $\Gamma^{i3}_k=0$ since $t_3=2\pi\sqrt{-1}\tau$ is flat with respect to $\widetilde\nabla$.
\begin{prop}
We have 
\[
(\Gamma^{ij}_{1})=
\begin{pmatrix}
0 & 0 & 0\\
0 & \frac{1}{4} & 0\\
1 & 0 & 0
\end{pmatrix},\quad 
(\Gamma^{ij}_{2})=
\begin{pmatrix}
-\frac{1}{6}t_2^3\frac{1}{(2\pi\sqrt{-1})^2}\frac{\p^2 E_2(\tau)}{\p\tau^2} & -\frac{1}{8}t_2^2\frac{1}{2\pi\sqrt{-1}}\frac{\p E_2(\tau)}{\p\tau}  & 0\\
-\frac{1}{4}t_2^2\frac{1}{2\pi\sqrt{-1}}\frac{\p E_2(\tau)}{\p\tau} & -\frac{1}{8}t_2E_2(\tau)  & 0\\
0 & \frac{1}{2} & 0
\end{pmatrix},
\]
\[
(\Gamma^{ij}_{3})=
\begin{pmatrix}
-\frac{1}{24}t_2^4\frac{1}{(2\pi\sqrt{-1})^3}\frac{\p^3 E_2(\tau)}{\p\tau^3} & -\frac{1}{24}t_2^3\frac{1}{(2\pi\sqrt{-1})^2}\frac{\p^2 E_2(\tau)}{\p\tau^2} & 0\\
-\frac{1}{12}t_2^3\frac{1}{(2\pi\sqrt{-1})^2}\frac{\p^2 E_2(\tau)}{\p\tau^2} & -\frac{1}{16}t_2^2\frac{1}{2\pi\sqrt{-1}}\frac{\p E_2(\tau)}{\p\tau}  & 0\\
0 & 0 & 0
\end{pmatrix},
\]
In particular, for all $i,j,k=1,2,3$,
\[
\frac{\p \Gamma^{ij}_{k}}{\p t_1}=0.
\]
\end{prop}
\begin{pf}
Except for $\Gamma^{12}_{k}$ and $\Gamma^{21}_{k}$, the calculations can be easily carried out using the result of Proposition~\ref{prop: invariant g}
and the equation \eqref{eq: contra a}. Since we have
\[
\frac{\p}{\p x}\left(\frac{\vartheta_{11}(x;\tau)}{\vartheta_{11}'(0;\tau)}\right)=\frac{\vartheta_{11}(x;\tau)}{\vartheta_{11}'(0;\tau)}\cdot \frac{\vartheta'_{11}(x;\tau)}{\vartheta_{11}'(0;\tau)}=(2\pi\sqrt{-1})^2
\frac{\vartheta_{11}(x;\tau)}{\vartheta_{11}'(0;\tau)}\left(\widetilde\zeta(x;\tau)+\frac{1}{12}E_2(\tau)x\right),
\]
we obtain $\Gamma^{12}_1=0$ and $\Gamma^{12}_2=-\frac{1}{8}t_2^2\frac{1}{2\pi\sqrt{-1}}\frac{\p E_2(\tau)}{\p\tau}$ by using \eqref{eq: key identity-tilde-1} and the result of Proposition~\ref{prop: invariant g}.
Hence $\Gamma^{21}_1=0$ and $\Gamma^{21}_2=-\frac{1}{4}t_2^2\frac{1}{2\pi\sqrt{-1}}\frac{\p E_2(\tau)}{\p\tau}$ by \eqref{eq: contra a}.
The remaining $\Gamma^{12}_3$ and $\Gamma^{21}_3$ can be derived from \eqref{eq: contra b}.
\qed
\end{pf}

Define the unit vector field and the Euler vector field as 
the holomorphic vector fields $e, E\in\Gamma(\widetilde \EE/\!/\widetilde W,\T_{\widetilde \EE/\!/\widetilde W})$ given by
\[
e:=\frac{\p}{\p y_1}=\frac{\p}{\p t_1},\quad E:=\frac{1}{2\pi\sqrt{-1}}\frac{\p}{\p \phi}=t_1\frac{\p}{\p t_1}+\frac{1}{2}t_2\frac{\p}{\p t_2},
\]
The non-degenerate symmetric $\O_{\widetilde \EE/\!/\widetilde W}$-bilinear form $\eta_{\widetilde W}$ on $\Omega^1_{\widetilde \EE/\!/\widetilde W}$
given in Proposition~\ref{prop: invariant g} induces the one on $\T_{\widetilde \EE/\!/\widetilde W}$, which we denote by the same symbol $\eta_{\widetilde W}$.
Define the product $\circ$ on $\T_{\widetilde \EE/\!/\widetilde W}$ by 
\[
\frac{\p}{\p t_i}\circ \frac{\p}{\p t_j}:=\sum_{k=1}^3 C_{ij}^k \frac{\p}{\p t_k},
\]
where 
\[
C_{ij}^k=\sum_{a=1}^3 \eta_{\widetilde W}\left(\frac{\p}{\p t_i},\frac{\p}{\p t_a}\right)C^{ka}_j,\quad 
C^{ka}_j =\begin{cases} 
\Gamma^{k1}_j & a= 1\\
2\Gamma^{k2}_j & a= 2\\
\delta_{kj} & a=3
\end{cases}.
\]

\begin{cor}\label{inv_Frob}
The tuple $(\circ, \eta_{\widetilde W}, e, E)$ defines a Frobenius structure on $\EE/\!/\widetilde W\cong \CC\times \CC\times \HH$ 
of rank $3$ and dimension $1$,
whose intersection form $g$ on $\Omega^1_{\EE/\!/\widetilde W}$ is 
\[
(g(dt_i , dt_j))=
\begin{pmatrix}
-\frac{t_2^4}{12}\frac{1}{(2\pi\sqrt{-1})^2}\frac{\p^2 E_2(\tau)}{\p \tau^2} & -\frac{t_2^3}{8}\frac{1}{2\pi\sqrt{-1}}\frac{\p E_2(\tau)}{\p \tau}   & t_1\\
-\frac{t_2^3}{8}\frac{1}{2\pi\sqrt{-1}}\frac{\p E_2(\tau)}{\p \tau}   & \frac{1}{2}t_1-\frac{1}{8}t_2^2E_2(\tau) & \frac{1}{2}t_2\\
t_1 & \frac{1}{2}t_2 & 0
\end{pmatrix}.
\]
Therefore, the Frobenius manifold restricted on $M$ is isomorphic to the one in Corollary \ref{primitive_form}.
\qed
\end{cor}

%%%%%%%%%%%%%%%%%%%%%%%%%%%%%%%%%%%%%%%%%%%%%%%%%%%%%%%%%%%%%%%%
\section{Gamma integral structure}\label{sec:Gamma}

Based on the gamma integral structure of Iritani \cite[Section~2.4]{I}, we define a map
\[
{\rm ch}_\Gamma: K_0(\D)\longrightarrow H\!H_*(\D)
\] 
by the following matrix representation with respect to the basis $\{[P(3)], [P(2)],[P(1)]\}$
{\small 
\[
{\rm ch}_\Gamma:=
\begin{pmatrix}
1 & 0 & 0\\
0 & \Gamma\left(\frac{1}{2}\right) & 0\\
0 & 0 & 1
\end{pmatrix}\cdot
\begin{pmatrix}
1 & 0 & 0\\
0 & -1 & 0\\
0 & 0 & 2\pi\sqrt{-1}
\end{pmatrix}\cdot
\begin{pmatrix}
1 & 1 & 0\\
-1 & -1 & -1\\
0 & 1 & 1
\end{pmatrix}
=\begin{pmatrix}
1 & 1 & 0\\
\sqrt{\pi} & \sqrt{\pi} & \sqrt{\pi}\\
0 & 2\pi\sqrt{-1} & 2\pi\sqrt{-1}
\end{pmatrix}.
\]}
and the identification $H\!H_*(\D)\cong H^{0,0}(\PP^1;\CC)\oplus \CC\oplus H^{1,1}(\PP^1;\CC)$,
where the (normalized) degree operator $Q$, defined by 
\[
Q:=
\begin{pmatrix}
-\frac{1}{2} & 0 & 0\\
0 & 0 & 0\\
0 & 0 & \frac{1}{2}
\end{pmatrix},
\]
acts diagonally on the right hand side.
Here the middle summand $\CC$ can be considered as a ``twisted sector'' (of the form $H^{\frac{1}{2},\frac{1}{2}}$) coming from the non-commutative resolution of the nodal cubic.

It is important that we have
\[
\left(\frac{1}{(2\pi)^{\frac{1}{2}}}{\rm ch}_\Gamma\right)^{-1}\cdot {\bf e}[Q]\cdot\left(\frac{1}{(2\pi)^{\frac{1}{2}}}{\rm ch}_\Gamma\right)
=\begin{pmatrix}
1 & 2 & 2\\
-2 & -3 & -2\\
2 & 2 & 1
\end{pmatrix}
=\chi_\D^{-1}\chi_\D^T,
\]
\[
\left(\frac{1}{(2\pi)^{\frac{1}{2}}}{\rm ch}_\Gamma\right)^T\cdot {\bf e}\left[\frac{1}{2}Q\right]\cdot\eta\cdot\left(\frac{1}{(2\pi)^{\frac{1}{2}}}{\rm ch}_\Gamma\right)
=\begin{pmatrix}
1 & 2 & 2\\
0 & 1 & 2\\
0 & 0 & 1
\end{pmatrix}
=\chi_\D.
\]

Note that the bifurcation set $\B\subset M$ is empty. In particular, the Frobenius structure on $M$ is semi-simple at any point. 

By mirror symmetry, it is expected that ${\rm Stab}(\D)$ should be equipped with a Frobenius structure.
As explained above, just as there are three types of descriptions of derived categories, 
there are also three different constructions of Frobenius structures, 
from the Gromov--Witten theory, the deformation theory with a choice of a primitive form, and 
the invariant theory of Weyl groups. 
Under this expectation, especially in the case of the smooth algebraic varieties, 
Dubrovin~\cite{D3} conjectured that the Frobenius structure above is semi-simple if and only if 
the bounded derived category of coherent sheaves has a full exceptional collection,
and that the monodromy data of the semi-simple Frobenius structure is compatible with numerical data of the full exceptional collection.
Dubrovin's conjecture has since been refined to Gamma Conjecture II (Galkin--Golyshev--Iritani~\cite{GGI}) and
the refined Dubrovin conjecture (Cotti--Dubrovin--Guzzetti~\cite{CDG}). 

Unfortunately, the Gromov--Witten theory for non-commutative algebraic varieties has not been constructed yet.
However, if it is made, it is also very natural from the perspective of mirror symmetry to expect that the Frobenius structure in this paper 
should be isomorphic to the one obtained from the ``Gromov--Witten theory for the non-commutative algebraic variety 
whose category of coherent sheaves is ${\rm coh}(\A_{\bf E})$''. 
Indeed, we have the following analogy for the above conjectures:
\begin{prop}\label{Gamma}
Choose a point on $M$ and an admissible line $l$, 
\begin{enumerate}
\item
There is a full exceptional collection $(E_1,E_2,E_3)$ such that 
the Stokes matrix $S$ of the Frobenius structure is given by $S=(\chi_\D([E_i],[E_j]))$.
\item
The central connection matrix of the Frobenius structure is given by $(2\pi)^{-1/2}{\rm ch }_\Gamma$.
\end{enumerate}
\end{prop}
\begin{rem}
The statement does not depend on the choice of a point on $M$ and an admissible line $l$ 
due to \cite[Isomonodromicity Theorem (second part)]{D2} and the compatibility of the two braid group actions: 
the one on the set of pairs of Stokes and central connection matrices and the one on the set of full exceptional collections.
\end{rem}
\begin{pf}
Fix a point ${\bf t}\in M$ such that $u\in \RR_{>0}$, $t_2\in\RR_{>0}$ and $\tau\in\sqrt{-1}\RR_{>0}$. 
Choose as a full exceptional collection $(P(3),P(2),P(1))$, or equivalently $(L_1,L_2,L_3)$, and denote by $(\gamma_1, \gamma_2,\gamma_3)$ the corresponding elements in the Grothendieck group, identified with $H_1(\overline{p}^{-1}({\bf t}), \overline{F}_{\bf t}^{-1}(+\infty); \ZZ)$.

Since $2\pi\sqrt{-1}dz$ is a primitive form, for any $\gamma\in H_1(\overline{p}^{-1}({\bf t}), \overline{F}_{\bf t}^{-1}(+\infty); \ZZ)$, 
the exterior derivative of the exponential period
\[
d\left(u^{-\frac{1}{2}}\int_{\gamma}  e^{-\frac{F(z;{\bf t})}{u}} (2\pi\sqrt{-1}dz)\right)
\]
is flat with respect to the first structure connection.
Therefore, by computing the behavior as $u\to +\infty$ of 
\[
\sum_{a=1}^3\eta^{ia}\cdot (-u)\frac{\p}{\p t_a}\left(u^{-\frac{1}{2}}\int_{\gamma_j}  e^{-\frac{F(z;{\bf t})}{u}} (2\pi\sqrt{-1}dz)\right),\quad i,j=1,2,3,
\]
and comparing it with the fundamental matrix at $u=\infty$ in the sense of \cite[Theorem~2.1]{D}, we obtain ${\rm ch}_\Gamma$.
Indeed, we have 
\begin{eqnarray*}
-u^{-\frac{1}{2}}\int_{\frac{1}{2}}^{\frac{1}{2}+\tau} e^{-\frac{F(z;{\bf t})}{u}} (2\pi\sqrt{-1}dz)&=&-2\pi\sqrt{-1}\tau\cdot u^{-\frac{1}{2}}-(2\pi\sqrt{-1}\tau t_1+t_2^2)\cdot u^{-\frac{3}{2}}+O(u^{-\frac{5}{2}}),\\
u^{-\frac{1}{2}}\left(\int_0^\tau -\int_{\frac{1}{2}}^{\frac{1}{2}+\tau}\right)e^{-\frac{F(z;{\bf t})}{u}} (2\pi\sqrt{-1}dz)&=&2\sqrt{\pi}t_2\cdot u^{-1}-2\sqrt{\pi}t_2t_1\cdot u^{-2}+O(u^{-3}),\\
u^{-\frac{1}{2}}\int_{\frac{\tau}{2}}^{1+\frac{\tau}{2}} e^{-\frac{F(z;{\bf t})}{u}} (2\pi\sqrt{-1}dz)&=&2\pi\sqrt{-1}\cdot u^{-\frac{1}{2}}-2\pi\sqrt{-1}t_1\cdot u^{-\frac{3}{2}}+O(u^{-\frac{5}{2}}),
\end{eqnarray*}
and hence
\[
-\delta_2=\begin{pmatrix}0\\ 1\\ -1\end{pmatrix}\mapsto \begin{pmatrix}1\\ 0\\ 0\end{pmatrix},\quad 
-\alpha=\begin{pmatrix}-1\\ 1\\ -1\end{pmatrix}\mapsto \begin{pmatrix}0\\ -\sqrt{\pi}\\ 0\end{pmatrix},\quad 
\delta_1=\begin{pmatrix}-1\\ 1\\ 0\end{pmatrix}\mapsto \begin{pmatrix}0\\ 0\\ 2\pi\sqrt{-1}\end{pmatrix}
\]
\qed\end{pf}
%%%%%%%%%%%%%%%%%%%%%%%%%%%%%%%%%%%%%%%%%%%%%%%%%%%%%%%%%%%%%%%%
\section{Appendix}
\subsection{Proof of Proposition~\ref{prop: primitive form}}\label{subsec: proof-prop}
It follows that
{\tiny 
\begin{eqnarray*}
& &\frac{\p F(z;{\bf t})}{\p t_2}\frac{\p F(z;{\bf t})}{\p t_2}\\
&=&4t_2^2\left(\widetilde\wp(z;\tau)-\frac{1}{12}E_2(\tau)\right)^2\\
&=&2t_2^2\frac{1}{2\pi\sqrt{-1}}\frac{\p}{\p \tau}\left(\widetilde\wp(z;\tau)-\frac{1}{12}E_2(\tau)\right)
-t_2^2E_2(\tau)\left(\widetilde\wp(z;\tau)-\frac{1}{12}E_2(\tau)\right)
-\frac{t_2^2}{2}\frac{1}{2\pi\sqrt{-1}}\frac{\p E_2(\tau)}{\p\tau}\\
& &+2t_2^2\left(-\widetilde\zeta(z;\tau)-\frac{1}{12}E_2(\tau)z\right)\frac{\p }{\p z}\left(\widetilde\wp(z;\tau)-\frac{1}{12}E_2(\tau)\right),\\
&=&2\frac{\p F(z;{\bf t})}{\p t_3}
-\frac{t_2}{2}E_2(\tau)\frac{\p F(z;{\bf t})}{\p t_2}-\frac{t_2^2}{2}\frac{1}{2\pi\sqrt{-1}}\frac{\p E_2(\tau)}{\p\tau}\frac{\p F(z;{\bf t})}{\p t_1}+2\left(-\widetilde\zeta(z;\tau)-\frac{1}{12}E_2(\tau)z\right)\frac{\p F(z;{\bf t})}{\p z},
\end{eqnarray*}}
and
{\tiny
\[
2\frac{\p}{\p z}\left(-\widetilde\zeta(z;\tau)-\frac{1}{12}E_2(\tau)z\right)=2\left(\widetilde\wp(z;\tau)-\frac{1}{12}E_2(\tau)\right)=\frac{\p^2 F(z;{\bf t})}{\p t_2 \p t_2}.
\]}

It follows that
{\tiny
\begin{eqnarray*}
& &\frac{\p F(z;{\bf t})}{\p t_2}\frac{\p F(z;{\bf t})}{\p t_3}\\
&=&2t_2^3\left(\widetilde\wp(z;\tau)-\frac{1}{12}E_2(\tau)\right)\cdot \frac{1}{2\pi\sqrt{-1}}\frac{\p}{\p \tau}\left(\widetilde\wp(z;\tau)-\frac{1}{12}E_2(\tau)\right)\\
&=&4t_2^3\left(\widetilde\wp(z;\tau)-\frac{1}{12}E_2(\tau)\right)^3+t_2^3E_2(\tau)\left(\widetilde\wp(z;\tau)-\frac{1}{12}E_2(\tau)\right)^2
+\frac{t_2^3}{2}\frac{1}{2\pi\sqrt{-1}}\frac{\p E_2(\tau)}{\p\tau}\left(\widetilde\wp(z;\tau)-\frac{1}{12}E_2(\tau)\right)\\
& &-2t_2^3\left(\widetilde\wp(z;\tau)-\frac{1}{12}E_2(\tau)\right)\left(-\widetilde\zeta(z;\tau)-\frac{1}{12}E_2(\tau)z\right)\frac{\p }{\p z}\left(\widetilde\wp(z;\tau)-\frac{1}{12}E_2(\tau)\right),\\
&=&
-\frac{t_2^3}{2}\frac{1}{2\pi\sqrt{-1}}\frac{\p E_2(\tau)}{\p\tau} \left(\widetilde\wp(z;\tau)-\frac{1}{12}E_2(\tau)\right)
-\frac{t_2^3}{6}\frac{1}{(2\pi\sqrt{-1})^2}\frac{\p^2 E_2(\tau)}{\p \tau^2} +t_2^3\frac{1}{(2\pi\sqrt{-1})^2}\left(\frac{\p}{\p z}\left(\widetilde\wp(z;\tau)-\frac{1}{12}E_2(\tau)\right)\right)^2\\
& &-2t_2^3\left(\widetilde\wp(z;\tau)-\frac{1}{12}E_2(\tau)\right)\left(-\widetilde\zeta(z;\tau)-\frac{1}{12}E_2(\tau)z\right)\frac{\p }{\p z}\left(\widetilde\wp(z;\tau)-\frac{1}{12}E_2(\tau)\right),\\
&=&-\frac{t_2^2}{4}\frac{1}{2\pi\sqrt{-1}}\frac{\p E_2(\tau)}{\p\tau}\frac{\p F(z;{\bf t})}{\p t_2}
-\frac{t_2^3}{6}\frac{1}{(2\pi\sqrt{-1})^2}\frac{\p^2 E_2(\tau)}{\p \tau^2}\frac{\p F(z;{\bf t})}{\p t_1} \\
& &+t_2\left(\frac{1}{(2\pi\sqrt{-1})^2}\frac{\p}{\p z}\left(\widetilde\wp(z;\tau)-\frac{1}{12}E_2(\tau)\right)-2\left(\widetilde\wp(z;\tau)-\frac{1}{12}E_2(\tau)\right)\left(-\widetilde\zeta(z;\tau)-\frac{1}{12}E_2(\tau)z\right)\right)\frac{\p F(z;{\bf t})}{\p z}\\
&=&-\frac{t_2^2}{4}\frac{1}{2\pi\sqrt{-1}}\frac{\p E_2(\tau)}{\p\tau}\frac{\p F(z;{\bf t})}{\p t_2}
-\frac{t_2^3}{6}\frac{1}{(2\pi\sqrt{-1})^2}\frac{\p^2 E_2(\tau)}{\p \tau^2}\frac{\p F(z;{\bf t})}{\p t_1}\\ 
& & +2t_2\frac{1}{2\pi\sqrt{-1}}\frac{\p}{\p \tau}\left(-\widetilde\zeta(z;\tau)-\frac{1}{12}E_2(\tau)z\right)\frac{\p F(z;{\bf t})}{\p z},
\end{eqnarray*}}
and 
{\tiny
\[
\frac{\p}{\p z}\left(2t_2\frac{1}{2\pi\sqrt{-1}}\frac{\p}{\p \tau}\left(-\widetilde\zeta(z;\tau)-\frac{1}{12}E_2(\tau)z\right)\right)=2t_2\frac{1}{2\pi\sqrt{-1}}\frac{\p}{\p \tau}\left(\widetilde\wp(z;\tau)-\frac{1}{12}E_2(\tau)\right)=\frac{\p^2 F(z;{\bf t})}{\p t_2 \p t_3}.
\]}

It follows that
{\tiny
\begin{eqnarray*}
& &\frac{\p F(z;{\bf t})}{\p t_3}\frac{\p F(z;{\bf t})}{\p t_3}\\
&=&t_2^4\left(\frac{1}{2\pi\sqrt{-1}}\frac{\p}{\p \tau}\left(\widetilde\wp(z;\tau)-\frac{1}{12}E_2(\tau)\right)\right)^2\\
&=& 4t_2^4\left(\widetilde\wp(z;\tau)-\frac{1}{12}E_2(\tau)\right)^4+2t_2^4E_2(\tau)\left(\widetilde\wp(z;\tau)-\frac{1}{12}E_2(\tau)\right)^3+t_2^4\left(\frac{1}{2\pi\sqrt{-1}}\frac{\p E_2(\tau)}{\p\tau} +\frac{1}{2}E_2(\tau)^2\right)\left(\widetilde\wp(z;\tau)-\frac{1}{12}E_2(\tau)\right)^2\\
& &+\frac{t_2^4}{4}\frac{1}{2\pi\sqrt{-1}}\frac{\p E_2(\tau)}{\p\tau}\cdot E_2(\tau)+\frac{t_2^4}{16}\left(\frac{1}{2\pi\sqrt{-1}}\frac{\p E_2(\tau)}{\p\tau}\right)^2
-t_2^4\left(-\widetilde\zeta(z;\tau)-\frac{1}{12}E_2(\tau)z\right)\frac{\p}{\p z}\left(\widetilde\wp(z;\tau)-\frac{1}{12}E_2(\tau)\right)\\
& &\cdot \left(4\left(\widetilde\wp(z;\tau)-\frac{1}{12}E_2(\tau)\right)^2
+E_2(\tau)\left(\widetilde\wp(z;\tau)-\frac{1}{12}E_2(\tau)\right)
+\frac{1}{2}\frac{1}{2\pi\sqrt{-1}}\frac{\p E_2(\tau)}{\p\tau}\right)\\
& & +t_2^4\left(-\widetilde\zeta(z;\tau)-\frac{1}{12}E_2(\tau)z\right)^2\left(\frac{\p}{\p z}\left(\widetilde\wp(z;\tau)-\frac{1}{12}E_2(\tau)\right)\right)^2\\
&=& -\frac{t_2^4}{6}\frac{1}{(2\pi\sqrt{-1})^2}\frac{\p^2 E_2(\tau)}{\p\tau^2}\left(\widetilde\wp(z;\tau)-\frac{1}{12}E_2(\tau)\right)-\frac{t_2^4}{24}\frac{1}{(2\pi\sqrt{-1})^3}\frac{\p^3 E_2(\tau)}{\p \tau^3}\\
& & +t_2^4\frac{1}{(2\pi\sqrt{-1})^2}\left(\widetilde\wp(z;\tau)-\frac{1}{12}E_2(\tau)\right)\left(\frac{\p}{\p z}\left(\widetilde\wp(z;\tau)-\frac{1}{12}E_2(\tau)\right)\right)^2
+\frac{t_2^4}{4}E_2(\tau)\left(\frac{\p}{\p z}\left(\widetilde\wp(z;\tau)-\frac{1}{12}E_2(\tau)\right)\right)^2\\
& &+2t_2^4\left(-\widetilde\zeta(z;\tau)-\frac{1}{12}E_2(\tau)z\right)\left(\widetilde\wp(z;\tau)-\frac{1}{12}E_2(\tau)\right)^2\frac{\p}{\p z}\left(\widetilde\wp(z;\tau)-\frac{1}{12}E_2(\tau)\right)\\
& &
-t_2^4\frac{1}{(2\pi\sqrt{-1})^2}\frac{\p^2}{\p z^2}\left(\widetilde\wp(z;\tau)-\frac{1}{12}E_2(\tau)\right)\frac{\p}{\p z}\left(\widetilde\wp(z;\tau)-\frac{1}{12}E_2(\tau)\right)\\
& & +t_2^4\left(-\widetilde\zeta(z;\tau)-\frac{1}{12}E_2(\tau)z\right)^2\left(\frac{\p}{\p z}\left(\widetilde\wp(z;\tau)-\frac{1}{12}E_2(\tau)\right)\right)^2\\
&=&-\frac{t_2^3}{12}\frac{1}{(2\pi\sqrt{-1})^2}\frac{\p^2 E_2(\tau)}{\p\tau^2} \frac{\p F(z;{\bf t})}{\p t_2}
-\frac{t_2^4}{24}\frac{1}{(2\pi\sqrt{-1})^3}\frac{\p^3 E_2(\tau)}{\p \tau^3}\frac{\p F(z;{\bf t})}{\p t_1}+\phi_{3,3}(z;{\bf t})\frac{\p F(z;{\bf t})}{\p z},
\end{eqnarray*}}
where, since \eqref{eq: key identity-tilde} yields
{\tiny
\begin{eqnarray*}
& &\frac{1}{(2\pi\sqrt{-1})^3}\frac{\p^2}{\p \tau\p z}\left(\widetilde\wp(z;\tau)-\frac{1}{12}E_2(\tau)\right)\\
&=&-\frac{1}{(2\pi\sqrt{-1})^2}\left(-\widetilde\zeta(z;\tau)-\frac{1}{12}E_2(\tau)z\right)\cdot\frac{\p^2 }{\p z^2}\left(\widetilde\wp(z;\tau)-\frac{1}{12}E_2(\tau)\right)
+\left(3\left(\widetilde\wp(z;\tau)-\frac{1}{12}E_2(\tau)\right)+\frac{1}{2}E_2(\tau)\right)\cdot\frac{\p }{\p z}\left(\widetilde\wp(z;\tau)-\frac{1}{12}E_2(\tau)\right),
\end{eqnarray*}}
$\phi_{3,3}(z;{\bf t})$ is given by
{\tiny
\begin{eqnarray*}
& &\phi_{3,3}(z;{\bf t})\\
&:=&t_2^2\frac{1}{(2\pi\sqrt{-1})^2}\left(\widetilde\wp(z;\tau)-\frac{1}{12}E_2(\tau)\right)\frac{\p}{\p z}\left(\widetilde\wp(z;\tau)-\frac{1}{12}E_2(\tau)\right)
+\frac{t_2^2}{4}E_2(\tau)\frac{\p}{\p z}\left(\widetilde\wp(z;\tau)-\frac{1}{12}E_2(\tau)\right)\\
& &+2t_2^2\left(-\widetilde\zeta(z;\tau)-\frac{1}{12}E_2(\tau)z\right)\left(\widetilde\wp(z;\tau)-\frac{1}{12}E_2(\tau)\right)^2
-t_2^2\frac{1}{(2\pi\sqrt{-1})^2}\frac{\p^2}{\p z^2}\left(\widetilde\wp(z;\tau)-\frac{1}{12}E_2(\tau)\right)\\
& & +t_2^2\left(-\widetilde\zeta(z;\tau)-\frac{1}{12}E_2(\tau)z\right)^2\frac{\p}{\p z}\left(\widetilde\wp(z;\tau)-\frac{1}{12}E_2(\tau)\right)\\
&=&-t_2^2\left(-\widetilde\zeta(z;\tau)-\frac{1}{12}E_2(\tau)z\right)\left(\widetilde\wp(z;\tau)-\frac{1}{12}E_2(\tau)\right)^2-\frac{t_2^2}{2}\frac{1}{(2\pi\sqrt{-1})^2}\left(\widetilde\wp(z;\tau)-\frac{1}{12}E_2(\tau)\right)\frac{\p}{\p z}\left(\widetilde\wp(z;\tau)-\frac{1}{12}E_2(\tau)\right)\\
& &+t_2^2\left(-\widetilde\zeta(z;\tau)-\frac{1}{12}E_2(\tau)z\right)^2\frac{\p }{\p z}\left(\widetilde\wp(z;\tau)-\frac{1}{12}E_2(\tau)\right)\\
& &-t_2^2\left(-\widetilde\zeta(z;\tau)-\frac{1}{12}E_2(\tau)z\right)\cdot\left(\frac{1}{2}E_2(\tau)\left(\widetilde\wp(z;\tau)-\frac{1}{12}E_2(\tau)\right)+\frac{1}{4}\frac{1}{2\pi\sqrt{-1}}\frac{\p E_2(\tau)}{\p\tau}\right)\\
&  &-\frac{t_2^2}{2}\frac{1}{(2\pi\sqrt{-1})^2}\left(-\widetilde\zeta(z;\tau)-\frac{1}{12}E_2(\tau)z\right)\cdot\frac{\p^2 }{\p z^2}\left(\widetilde\wp(z;\tau)-\frac{1}{12}E_2(\tau)\right)\\
& &+\frac{t_2^2}{2}\left(3\left(\widetilde\wp(z;\tau)-\frac{1}{12}E_2(\tau)\right)+\frac{1}{2}E_2(\tau)\right)\cdot\frac{\p }{\p z}\left(\widetilde\wp(z;\tau)-\frac{1}{12}E_2(\tau)\right)\\
&=& t_2^2\left(-\widetilde\zeta(z;\tau)-\frac{1}{12}E_2(\tau)z\right)\left(\widetilde\wp(z;\tau)-\frac{1}{12}E_2(\tau)\right)^2
-\frac{1}{2}\frac{1}{(2\pi\sqrt{-1})^2}\frac{\p}{\p z}\left(\widetilde\wp(z;\tau)-\frac{1}{12}E_2(\tau)\right)\cdot\left(\widetilde\wp(z;\tau)-\frac{1}{12}E_2(\tau)\right)\\
& &-t_2^2\left(-\widetilde\zeta(z;\tau)-\frac{1}{12}E_2(\tau)z\right)\cdot 2\left(\widetilde\wp(z;\tau)-\frac{1}{12}E_2(\tau)\right)^2\\
& &-t_2^2\left(-\widetilde\zeta(z;\tau)-\frac{1}{12}E_2(\tau)z\right)\cdot\left(\frac{1}{2}E_2(\tau)\left(\widetilde\wp(z;\tau)-\frac{1}{12}E_2(\tau)\right)+\frac{1}{4}\frac{1}{2\pi\sqrt{-1}}\frac{\p E_2(\tau)}{\p\tau}\right)\\
& &+t_2^2\left(-\widetilde\zeta(z;\tau)-\frac{1}{12}E_2(\tau)z\right)^2\frac{\p }{\p z}\left(\widetilde\wp(z;\tau)-\frac{1}{12}E_2(\tau)\right)\\
& &+\frac{t_2^2}{2}\frac{1}{(2\pi\sqrt{-1})^3}\frac{\p^2}{\p \tau\p z}\left(\widetilde\wp(z;\tau)-\frac{1}{12}E_2(\tau)\right)\\
&=&-t_2^2\frac{1}{2\pi\sqrt{-1}}\frac{\p}{\p\tau}\left(-\widetilde\zeta(z;\tau)-\frac{1}{12}E_2(\tau)z\right)\cdot\left(\widetilde\wp(z;\tau)-\frac{1}{12}E_2(\tau)\right)\\
& &-t_2^2\left(-\widetilde\zeta(z;\tau)-\frac{1}{12}E_2(\tau)z\right)\cdot\frac{1}{2\pi\sqrt{-1}}\frac{\p}{\p\tau}\left(\widetilde\wp(z;\tau)-\frac{1}{12}E_2(\tau)\right)\\
&  &+\frac{t_2^2}{2}\frac{1}{(2\pi\sqrt{-1})^3}\frac{\p^2}{\p \tau\p z}\left(\widetilde\wp(z;\tau)-\frac{1}{12}E_2(\tau)\right)\\
&=&t_2^2\frac{1}{(2\pi\sqrt{-1})^2}\frac{\p^2}{\p\tau^2}\left(-\widetilde\zeta(z;\tau)-\frac{1}{12}E_2(\tau)z\right).
\end{eqnarray*}}
and
{\tiny
\[
\frac{\p}{\p z}\left(
t_2^2\frac{1}{(2\pi\sqrt{-1})^2}\frac{\p^2}{\p \tau^2}\left(-\widetilde\zeta(z;\tau)-\frac{1}{12}E_2(\tau)z\right)\right)=\frac{\p^2 F(z;{\bf t})}{\p t_3^2}.
\]}

%%%%%%%%%%%%%%%%%%%%%%%%%%%%%%%%%%%%%%%%%%%%%%%%%%%%%%%%%%%%%%%%%%%%%%%%%%%%%%%%%%%%%%


\begin{thebibliography}{99}

\bibitem[A]{A}
G.~F.~Almeida,
\textit{The differential geometry of the orbit space of extended affine Jacobi group $A_1$},
SIGMA {\bf 17} (2021), 022, 39~pp.

\bibitem[AS]{AS}
I.~Assem and A.~Skowro\'{n}ski,
\textit{Iterated tilted algebras of type $\widetilde{\mathbb{A}}_n$},
Math. Z. {\bf 195} (1987) 269--290.

\bibitem[BH]{BH}
S.~Balnojan and C.~Hertling,
\textit{Conjectures on spectral numbers for upper triangular matrices and for singularities},
Math. Phys. Anal. Geom. \textbf{23} (2020), 5, 49~pp.

\bibitem[Be1]{Be1}
M.~Bertola,
\textit{Frobenius manifold structure on orbit space of Jacobi groups. I.}
Differential Geom. Appl. \textbf{13} (2000), 19--41.

\bibitem[Be2]{Be2}
M.~Bertola,
\textit{Frobenius manifold structure on orbit space of Jacobi groups. II.}
Differential Geom. Appl. \textbf{13} (2000), 213--233.

\bibitem[Bo]{Bo}
A.~I.~Bondal,
\textit{Representations of associative algebras and coherent sheaves},
Math. USSR-Izv. \textbf{34} (1990), 23--42.

\bibitem[BK]{BK}
A.~I.~Bondal and M.~M.~Kapranov,
\textit{Representable functors, Serre functors, and reconstructions},
Math. USSR-Izv. \textbf{35} (1990), 519--541.

\bibitem[BP]{BP}
A.~I.~Bondal and A.~E.~Polishchuk,
\textit{Homological properties of associative algebras: the method of helices},
Russian Acad. Sci. Izv. Math. \textbf{42} (1994), 219--260.

\bibitem[B]{Br}
T.~Bridgeland,
\textit{Stability conditions on triangulated categories}, 
Ann. of Math. {\bf 166} (2007), 317--345.

\bibitem[BPP]{BPP}
N.~Broomhead, D.~Pauksztello, and D.~Ploog,
\textit{Discrete derived categories I: homomorphisms, autoequivalences and $t$-structures},
Math. Z. \textbf{285} (2017), 39--89.

\bibitem[BD]{BD}
I.~Burban and Y.~Drozd,
{\it Tilting on non-commutative rational projective curves},
Math. Ann. {\bf 351} (2011), 665--709.

\bibitem[BuKr]{BuKr}
I.~Burban and B.~Kreu\ss ler,
\textit{Derived categories of irreducible projective curves of arithmetic genus one},
Compos. Math. \textbf{142} (2006), 1231--1262.

\bibitem[CDG]{CDG}
G.~Cotti, B.~Dubrovin and D.~Guzzetti, 
\textit{Helix structures in quantum cohomology of Fano varieties}
arXiv:1811.09235.

\bibitem[CV]{CV}
S.~Cecotti and C.~Vafa,
\textit{On classification of $N=2$ supersymmetric theories},
Comm. Math. Phys. \textbf{158} (1993), 569--644.

\bibitem[CHS]{CHS}
W.~Chang, F.~Haiden, and S.~Schroll,
\textit{Braid group actions on branched coverings and full exceptional sequences},
Adv. Math. \textbf{472} (2025), Paper No.~110284, 24~pp.

\bibitem[CS]{CS}
W.~Chang and S.~Schroll,
\textit{Exceptional sequences in the derived category of a gentle algebra},
Selecta Math. (N.S.) \textbf{29} (2023), Paper No.~33, 41~pp.

\bibitem[Ch]{Ch}
X.~W.~Chen,
\textit{A note on standard equivalences},
Bull. Lond. Math. Soc. \textbf{48} (2016), 797--801.

\bibitem[Du1]{D}
B.~Dubrovin,
\textit{Geometry of $2$D topological field theories},
Integrable systems and quantum groups (Montecatini Terme, 1993), 120--348,
Lecture Notes in Math., \textbf{1620}, Springer, Berlin, 1996.

\bibitem[Du2]{D2}
B.~Dubrovin,
\textit{Painlevé transcendents in two-dimensional topological field theory},
The Painlevé property, 287--412,
CRM Ser. Math. Phys., Springer, New York, 1999.

\bibitem[Du3]{D3}
B.~Dubrovin, 
\textit{Geometry and Analytic theory of Frobenius manifolds}, 
Proceedings of the International Congress of Mathematicians, Vol. II (Berlin, 1998).
Doc. Math. 1998, Extra Vol. II, 315--326.

\bibitem[DZ]{DZ}
B.~Dubrovin and Y.~Zhang,
\textit{Extended affine Weyl groups and Frobenius manifolds},
Compos. Math. \textbf{111} (1998), 167--219.

\bibitem[ET]{ET}
W.~Ebeling and A.~Takahashi,
\textit{Variance of the exponents of orbifold Landau-Ginzburg models},
Math. Res. Lett. \textbf{20} (2013), 51--65.

\bibitem[EL]{EL}
A.~Elagin and V.~A.~Lunts,
\textit{Three notions of dimension for triangulated categories},
J. Algebra \textbf{569} (2021), 334--376.

\bibitem[GGI]{GGI}
S.~Galkin, V.~Golyshev, H.~Iritani, 
\textit{Gamma classes and quantum cohomology of Fano manifolds: Gamma conjectures}, 
Duke Math. J. 165. (2016), 2005--2077.

\bibitem[HKK]{HKK}
F.~Haiden, L.~Katzarkov, M.~Kontsevich,
\textit{Flat surfaces and stability structures},
Publ. Math. Inst. Hautes \'{E}tudes Sci. {\bf 126} (2017), 247--318.

\bibitem[HW]{HW}
F.~Haiden and D.~Wu,
\textit{A counterexample to the Jordan-Hölder property for polarizable semiorthogonal decompositions},
arXiv:2502.12075.

\bibitem[He1]{He1}
C.~Hertling,
\textit{Frobenius manifolds and variance of the spectral numbers},
New developments in singularity theory (Cambridge, 2000), 235--255,
NATO Sci. Ser. II Math. Phys. Chem., \textbf{21}, Kluwer Acad. Publ., Dordrecht, 2001.

\bibitem[He2]{He2}
C.~Hertling,
\textit{Frobenius manifolds and moduli spaces for singularities},
Cambridge Tracts in Math., \textbf{151}, Cambridge Univ. Press, Cambridge, 2002.

\bibitem[HR]{HR}
C.~Hertling and C.~Roucairol,
\textit{Distinguished bases and Stokes regions for the simple and the simple elliptic singularities},
Moduli spaces and locally symmetric spaces, 39--106,
Surv. Mod. Math., \textbf{16}, Int. Press, Somerville, MA, 2021.

\bibitem[IQ]{IQ}
A.~Ikeda and Y.~Qiu,
\textit{$q$-Stability conditions on Calabi--Yau-$\XX$ categories}, 
Compos. Math. \textbf{159} (2023), 1347--1386.

\bibitem[I]{I}
H.~Iritani, 
\textit{An integral structure in quantum cohomology and mirror symmetry for toric orbifolds}, 
Adv. Math. \textbf{222} (2009), 1016--1079.

\bibitem[IST]{IST}
Y.~Ishibashi, Y.~Shiraishi, and A.~Takahashi,
\textit{Primitive forms for affine cusp polynomials},
Tohoku Math. J. (2) \textbf{71} (2019), 437--464.

\bibitem[KOT]{KOT}
K.~Kikuta, G.~Ouchi, and A.~Takahashi,
\textit{Serre dimension and stability conditions},
Math. Z. \textbf{299} (2021), 997--1013.

\bibitem[Kon]{Kon}
M.~Kontsevich
\textit{Homological algebra of mirror symmetry},
Proceedings of the International Congress of Mathematicians (1994),
Birkhäuser, Basel, 1995, 120--139.

\bibitem[Ku]{Ku}
A.~Kuznetsov,
\textit{Derived categories of cubic fourfolds},
Cohomological and geometric approaches to rationality problems, 219--243,
Progr. Math., \textbf{282}, Birkhäuser Boston, Boston, MA, 2010.

\bibitem[KL]{KL}
A.~Kuznetsov and V.~A.~Lunts,
\textit{Categorical resolutions of irrational singularities},
Int. Math. Res. Not. IMRN \textbf{2015} (2015), 4536--4625.

\bibitem[KP]{KP}
A.~Kuznetsov and A.~Perry,
\textit{Serre functors and dimensions of residual categories},
arXiv:2109.02026.

\bibitem[LP]{LP}
Y.~Lekili and A.~E.~Polishchuk,
\textit{Auslander orders over nodal stacky curves and partially wrapped Fukaya categories},
J. Topol. \textbf{11} (2018), 615--644.

\bibitem[M]{M}
T.~Milanov, 
\textit{Primitive forms and Frobenius structures on the Hurwitz spaces},
arXiv:1701.00393.

\bibitem[O1]{O1}
S.~Opper,
\textit{On auto-equivalences and complete derived invariants of gentle algebras},
arXiv:1904.04859.

\bibitem[O2]{O2}
S.~Opper,
\textit{Spherical objects, transitivity and auto-equivalences of Kodaira cycles via gentle algebras},
J. Eur. Math. Soc. (2024), doi:10.4171/JEMS/1473.

\bibitem[OST]{OST}
T.~Otani, Y.~Shiraishi and A.~Takahashi
\textit{The number of full exceptional collections modulo spherical twists for extended Dynkin quivers},
arXiv:2308.04031.
  
\bibitem[Q1]{Q1}
Y.~Qiu,
\textit{Global dimension function on stability conditions and Gepner equations},
Math. Z. \textbf{303} (2023), Paper No.~11, 24~pp.

\bibitem[Q2]{Q2}
Y.~Qiu,
\textit{Contractible flow of stability conditions via global dimension function},
J. Differential Geom. \textbf{129} (2025), 491--521.

\bibitem[Sa1]{Sa1}
K.~Saito,
\textit{Period mapping associated to a primitive form},
Publ. Res. Inst. Math. Sci. \textbf{19} (1983), 1231--1264.

\bibitem[Sa2]{Sa2}
K.~Saito,
\textit{Extended affine root systems. II. Flat invariants},
Publ. Res. Inst. Math. Sci. \textbf{26} (1990), 15--78.

\bibitem[Sa3]{Sa3}
K.~Saito,
\textit{Around the theory of the generalized weight system: relations with singularity theory, the generalized Weyl group and its invariant theory, etc.},
Selected papers on harmonic analysis, groups, and invariants,
Amer. Math. Soc. Transl. \textbf{183}, 
Amer. Math. Soc., Providence, RI, 1998, pp.~101--143.

\bibitem[Sa4]{Sa4}
K.~Saito,
{\it On a linear structure of the quotient variety by a finite reflection group},
Publ. RIMS 1993 Volume 29 Issue 4 Pages 535--579.


\bibitem[SaTa]{SaTa}
K.~Saito and A.~Takahashi,
\textit{From primitive forms to Frobenius manifolds},
From Hodge theory to integrability and TQFT tt*-geometry, 31--48,
Proc. Sympos. Pure Math. \textbf{78}, Amer. Math. Soc., Providence, RI, 2008.

\bibitem[SYS]{SYS}
K.~Saito, T.~Yano and J.~Sekigichi,
\textit{On a certain generator system of the ring of invariants of a finite reflection group},
Comm. Algebra \textbf{8} (1980), 373--408.

\bibitem[Sat1]{Sat1}
I.~Satake,
\textit{Frobenius manifolds for elliptic root systems},
Osaka J. Math. \textbf{47} (2010), 301--330.

\bibitem[Sat2]{Sat2}
I.~Satake, 
\textit{Good Basic Invariants for Elliptic Weyl Groups and Frobenius Structures},
arXiv:2004.03587.

\bibitem[SatTa]{SatTa}
I.~Satake, A.~Takahashi,  
\textit{Gromov--Witten invariants for mirror orbifolds of simple elliptic singularities}
Annales de l'Institut Fourier, Volume 61 (2011) no. 7, pp. 2885--2907.

\bibitem[SeTh]{ST}
P.~Seidel, R.~Thomas, 
\textit{Braid group actions on derived categories of coherent sheaves},
Duke Math. J. {\bf 108} (2001), 37--108.

\bibitem[ShTa]{ShTa}
Y.~Shiraishi, A.~Takahashi,
\textit{On the Frobenius manifolds for cusp singularities}
Adv. Math. {\bf 273} (2015), 485--522.

\bibitem[STW]{STW}
Y.~Shiraishi, A.~Takahashi, K.~Wada,
\textit{On Weyl groups and Artin groups associated to orbifold projective lines},
J. Algebra {\bf 453} (2016), 249--290.

\bibitem[Su]{Su}
B.~Sung,
\textit{Remarks on the nodal quiver},
J. Algebra \textbf{664} (2025), 410--451.

\bibitem[SZ]{SZ}
J.~Schröer, A.~Zimmermann,
\textit{Stable endomorphism algebras of modules over special biserial algebras},
Math. Z. \textbf{244} (2003), 515--530.

\bibitem[Ta1]{Ta1}
A.~Takahashi,
\textit{Mirror Symmetry Between Orbifold Projective Lines and Cusp Singularities},
Singularities in geometry and topology 2011, 257--282,
Adv. Stud. Pure Math. \textbf{66}, Math. Soc. Japan, Tokyo, 2015.

\bibitem[Ta2]{Ta2}
A.~Takahashi,
\textit{Introduction to Primitive Forms and Mirror Symmetry},
Iwanami Studies in Advances Mathematics (Japanese), 
Iwanami Shoten, Tokyo, 2021.

\bibitem[TZ]{TZ}
A.~Takahashi and H.~Zhang,
\textit{Number of full exceptional collections modulo spherical twists for elliptic orbifolds},
J. Algebra \textbf{667} (2025), 570-586.

\bibitem[Tak]{Tak}
A.~Takeda,
\textit{Relative stability conditions on Fukaya categories of surfaces},
Math. Z. \textbf{301} (2022), 3019--3070.

\bibitem[Zac]{Zac}
J.~Juan-Zacarías
\textit{Some remarks on the correspondence between elliptic curves and four points in the Riemann sphere},
arXiv:1810.08742.



\end{thebibliography}
\end{document}